\documentclass[]{elsarticle}
\usepackage{lineno,hyperref}
\usepackage{algorithm}
\usepackage[noend]{algpseudocode}
\modulolinenumbers[5]
\makeatletter
\def\ps@pprintTitle{%
 \let\@oddhead\@empty
 \let\@evenhead\@empty
 \def\@oddfoot{}%
 \let\@evenfoot\@oddfoot}
\makeatother
\usepackage[a4paper, total={6.4in, 9.4in}]{geometry}
\usepackage[a4paper]{geometry}
 
\usepackage{color}
\usepackage[x11names]{xcolor}
\usepackage[utf8]{inputenc}
\usepackage[english]{babel}
\usepackage{amsmath}
\usepackage{pgfplots}
\usepackage{tabularx,ragged2e,booktabs,caption}

\usepackage[labelfont=bf]{caption}
\usepackage[font=footnotesize]{caption}
\usepackage{amssymb,latexsym}
\usepackage{graphicx}
\graphicspath{{./images/}}
\usepackage{booktabs}
\usepackage{amsfonts}
\usepackage{tikz}
\usepackage{array}
\usepackage{mathtools}
 \usepackage{relsize}
\usepackage[T1]{fontenc}
\usepackage{float}
\usepackage{cases}
\usepackage{multirow}
\usepackage{amsthm}
\usepackage{subcaption}
\pgfplotsset{compat=newest}
\usepackage[toc,page]{appendix}
\newcolumntype{Y}{>{\centering\arraybackslash}X}

\newtheorem{remark}{Remark}[section]
\numberwithin{equation}{section}
\usepackage{tabularx,ragged2e,booktabs,caption}

\usepackage[toc]{appendix}
\DeclareMathAlphabet{\mathpzc}{OT1}{pzc}{m}{it}
\newcommand{\lnorm}{\left|\left|}

\newcommand{\rnorm}{\right|\right|}

\newcommand{\Ld}{L^2(\mathcal{D})}
\newcommand{\Lomd}{L^2(\Omega,\mathcal{D})}
\title{A parametric acceleration of multilevel Monte Carlo convergence for nonlinear variably saturated flow}
\begin{document}
\begin{frontmatter}
\address[cwi]{CWI - Centrum Wiskunde \& Informatica, Amsterdam, The Netherlands}
\address[aero]{Faculty of Aerospace Engineering, Delft University of Technology, Delft, The Netherlands.}
\address[diam]{DIAM, Delft University of Technology, Delft, The Netherlands.}
\address[zara]{IUMA and Applied Mathematics Department, University of Zaragoza}
\author[cwi,aero]{Prashant~Kumar}
\ead{pkumar@cwi.nl}
\author[zara]{Carmen~Rodrigo}
\ead{carmenr@unizar.es}
\author[zara]{Francisco~J.~Gaspar}
\ead{fjgaspar@unizar.es}
\author[cwi,diam]{Cornelis~W.~Oosterlee}
\ead{c.w.oosterlee@cwi.nl}
\begin{abstract}
We present a multilevel Monte Carlo (MLMC) method for the uncertainty quantification of variably saturated porous media flow that are modeled using the Richards' equation. We propose a stochastic extension for the empirical models that are typically employed to close the Richards' equations. This is achieved by treating the soil parameters in these models as spatially correlated random fields with appropriately defined marginal distributions. As some of these parameters can only take values in a specific range, non-Gaussian models are utilized. The randomness in these parameters may result in path-wise highly nonlinear systems, so that a robust solver with respect to the random input is required. For this purpose, a solution method based on a combination of the modified Picard iteration and a cell-centered multigrid method for heterogeneous diffusion coefficients is utilized. Moreover, we propose a non-standard MLMC estimator to solve the resulting high-dimensional stochastic Richards' equation. The improved efficiency of this multilevel estimator is achieved by parametric continuation that allows us to incorporate simpler nonlinear problems on coarser levels for variance reduction while the target strongly nonlinear problem is solved only on the finest level. Several numerical experiments are presented showing computational savings obtained by the new estimator compared to the original MC estimator.
\end{abstract}
\end{frontmatter}
\section{Introduction}
Mass transport through a variably saturated porous medium can be accurately predicted using the Richards' equation \cite{richards}. This modelling approach is of critical importance for several physics and engineering problems, for instance, when studying aquifer recharge via rainfall infiltration, or for understanding the environmental impact of mining operations. When reliable measurements of the hydraulic properties are available, numerical solutions originating from the Richards' equation have been reasonably successful for transport prediction in a broad range of soil types.

Different formulations for the Richards' equation are available in the literature, along with well established mathematical theory, such as the pressure head, the water content or a mixed formulation, see e.g. \cite{Celia1990,FORSYTH199525,Arbogast1996,Eymard1999,Pop2002,Woodward2000}. The aforementioned formulations contain nonlinearities due to a parametric dependence of the pressure head on the saturation and the relative hydraulic conductivity. Depending on the soil parameters, these nonlinearities range from mild to strong. The extreme sensitivity of the soil parameters on the output necessitates accurate measurements of the hydraulic properties. For many realistic problems, complete information of these quantities is however not available. In such scenarios, these parameters may be modeled in a probabilistic framework and the solution output may be expressed by means of a prediction interval (with mean and variance), rather than as a single value. Such approaches are nowadays common in the case of saturated groundwater flow, where uncertainties are included when modeling the hydraulic conductivity as a spatially correlated lognormal random field \cite{WRCR:WRCR2456,WRCR:WRCR1821,WRCR:WRCR3746}. The purpose of the present work is to develop and analyze a stochastic extension of the Richards' equation, along with an efficient numerical method to solve the resulting nonlinear partial differential equation with random coefficients. 

Previous work on the uncertainty quantification (UQ) of unsaturated flows was often based on an uncertain hydraulic conductivity \cite{Mantoglou1987,Russo1997,Zhang1999,Iden2007}. In addition to that, in the present work, we introduce stochasticity in the so-called van Genuchten and Mualem model \cite{vanG, Mualem}, which is typically utilized to close the Richards' equation. This model provides a closed-form analytic expression for the unsaturated hydraulic conductivity based on a \emph{sigmoid} type function for the soil-water retention curve. This curve is defined by four independent parameters that are estimated by curve-fitting, based on field measurements. Typically, these parameters are fixed throughout the domain during numerical simulations, assuming the soil to be homogeneous. Realistic models should however also incorporate the intrinsic heterogeneity in the soil. Therefore, we model these soil parameters as random variables with a certain, specified probability distribution and spatial correlations. To assure the well-posedness of the Richards' equation, these parameters should be within a certain range. Thus, the probability distributions for these parameters are chosen such that the random samples will be in the domain of validity for these parameters. A practical choice is to employ non-Gaussian random fields with marginal distributions from expert knowledge or from field measurements.

With the stochastic Richards' equation formulated, an appropriate UQ technique is required to compute the statistical moments of the desired quantities of interest (QoIs). This choice primarily  depends on the number of uncertainty dimensions. Other practical factors, such as ease of implementation and availability of an iterative solver which is robust with respect to the random input, also play a role in the selection of a suitable UQ technique. The proposed stochastic extension of the Richards' equation results in a very high-dimensional problem, and the use of deterministic sampling approaches such as polynomial chaos expansion, stochastic collocation or stochastic Galerkin is therefore limited.  For these UQ methods, the cost grows exponentially with the number of random inputs. Furthermore, a deterministic sampling approach may not adequately represent those regions in the stochastic space where strong nonlinearity may be encountered. In previous works of Zhang \cite{Zhang1999,Zhang2002}, the moments method was applied for the uncertainty quantification of solutions of the Richards' equation. The main disadvantage of a moment-based method is that it can only be reliably employed when the effect of uncertain inputs is mild and largely linear. For the proposed stochastic formulation of the Richards' equation, Monte Carlo (MC) based sampling approaches are the \emph{methods of choice}, due to their dimension independent convergence. Moreover, these MC-type methods can accurately represent the entire stochastic space given a sufficiently large number of samples. A well-known drawback of the plain MC method is its slow convergence of the sampling error, with $\mathcal{O}(1/\sqrt{N})$, where $N$ is the number of samples, making it intractable for problems with a large cost per sample. Recently, efficient MC estimators based on the multilevel Monte Carlo (MLMC) method have been developed for a large class of problems, see, e.g. \cite{MLMC1,MLMC2,Mishra20123365,mishra2016multi}. The efficiency of the MLMC estimation comes from solving the problem of interest on a coarse grid and subsequently adding corrections based on finer mesh resolutions. As these correction terms have smaller variances, they can be computed accurately using only a few samples. The estimates at different levels are then combined using a telescopic sum. The standard practice is to solve the PDE with random coefficients on a hierarchy of grids. 

The original, grid-based MLMC estimator may be utilized to solve the stochastic Richards' equation, however, this approach may not be the most efficient, especially not when strongly nonlinear problems need to be solved. Such problems require a very fine spatio-temporal mesh thereby restricting the use of coarse grids to improve the efficiency of the MLMC estimator. In this article, we utilize a non-standard MLMC estimator based on the  \emph{parametric continuation} technique. Continuation methods for solving nonlinear PDEs are very popular in engineering applications \cite{Bra77,brandt2003multigrid,Mittelmann1,arclength,dinar,rheinboldt}. Within continuation methods, a nonlinearity dictating parameter $\Theta$  is introduced in the interval $\Theta_0\leq\Theta\leq\Theta_*$ where the solution $p(\Theta_0)$ corresponds to a linear (or mildly nonlinear) problem and $p(\Theta_*)$ to the target strongly nonlinear. The key idea is to march from $p(\Theta_0)$ to $p(\Theta_*)$ in small steps of size $\delta \Theta$, where at each step we use the solution from the previous step as an initial guess. Usually, $\Theta$ is some physical parameter, for e.g. the Reynolds number, the Mach number, etc. In the current work, we use parametric continuation to obtain variance reduction within the multilevel Monte Carlo framework. This is achieved by solving simpler nonlinear problems on coarser levels and the target strongly nonlinear problem is only solved on the finest level. This new estimator allows us to incorporate comparatively coarser spatio-temporal grids in the MLMC hierarchy and, as such, the computational cost of each estimator in the telescopic sum is greatly reduced. 

We furthermore propose a solution method for Richards' equation based on a combination of the modified Picard method \cite{Celia1990} and a cell-centered multigrid, as proposed in \cite{prashant3}. We benchmark the performance of this combined solver in a probabilistic framework. A number of tests for a wide range of soil parameters and for hydraulic conductivities with different heterogeneity levels is performed. 

The rest of the article is organized as follows. In Section \ref{deterRichards}, we briefly discuss the deterministic Richards' equation along with the van Genuchten-Mualem parameterization. Section \ref{stochRE} describes the stochastic Richards' equation as well as the modeling of various uncertain soil parameters. The description of the modified Picard method in combination with the cell-centered multigrid method is provided in Section \ref{MPCCMG}. Also, in this section we present some numerical experiments to assess the performance of the combined solver for an infiltration problem. The non-standard MLMC estimator is explained in Section \ref{nonMLMC} and its performance is analysed in Section \ref{numexp2}. Finally, some conclusions are drawn in Section \ref{conclu6}.

\section{Deterministic Richards' equation}\label{deterRichards}

We begin by describing the deterministic version of the problem. The governing equations are defined in a bounded domain $\mathcal{D}\subset \mathbb{R}^n$, with the boundary $\partial\mathcal{D}$ and a finite time interval $\mathcal{T}= (0,T_{final}]$, with $T_{final}<\infty$. The classical Richards' equation is a result of coupling the mass conservation equation of the water-phase and the Darcian flow, i.e.,
\begin{eqnarray}
\phi\dfrac{\partial S_w}{\partial t}  + \nabla\cdot\mathbf{q} = {f}\quad\text{in}\quad\mathcal{D}\times\mathcal{T},\label{continuity}\\
\mathbf{q} = -{K}_sK_{rw}(\nabla p+z)\quad\text{in}\quad\mathcal{D}\times\mathcal{T},\label{darcy}
\end{eqnarray}
respectively, subject to boundary and initial conditions:
\begin{align}\label{BCs}
p &= p_0\quad\text{in}\quad\mathcal{D},\quad t=0,\\
p &= g_D\quad\text{in}\quad\Gamma_D\times\mathcal{T},\\
\mathbf{q}\cdot \mathbf{n} &= g_N\quad\text{in}\quad\Gamma_N\times\mathcal{T},
\end{align}
where $\phi [L^3/L^3]$ is the porosity, $S_w [L^3/L^3]$ is the water-phase saturation; ${q}$ is the Darcy flux, which depends on the pressure head, $p [L]$, and the depth $z [L]$ in the vertical direction;  ${K}_s[L/T]$ represents the saturated hydraulic conductivity field at saturation; $K_{rw}$ is the relative conductivity of the water phase with respect to air and ${f}$ is the source/sink term. The initial pressure head value is given by $p_0$. The quantities $g_D$ and $g_N$ denote, respectively, the Dirichlet and Neumann boundary conditions that are imposed at the boundaries $\Gamma_D$ and $\Gamma_N$, respectively, with $\mathbf{n}$ the unit normal vector to $\Gamma_N$. 

The coupling of \eqref{continuity} and \eqref{darcy} may result in different variants of the Richards' equation, such as the pressure head, the moisture content and the mixed formulation. The mixed formulation of the Richards' equation is given by:
\begin{equation}\label{mixedForm}
\dfrac{\partial \theta(p)}{\partial t} - \nabla\cdot\left({K}_sK_{rw}(\nabla p+z)\right) = {f}\quad\text{in}\quad\mathcal{D}\times\mathcal{T}.
\end{equation}
% where $C(p):=\dfrac{\partial \theta}{\partial p}$ is the specific moisture capacity and $\theta:=\phi S_w(p)$ is the moisture content. 
It is obtained by substituting the moisture content, i.e,  $\theta =\phi S_w(p)$. By using $$\frac{\partial \theta(p)}{\partial t}  = C(p) \frac{\partial p}{\partial t},$$ the above PDE can be reformulated into the pressure head formulation:
\begin{equation}\label{pressurehead}
 C(p) \dfrac{\partial p}{\partial t} - \nabla\cdot\left({K}_sK_{rw}(\nabla p+z)\right) = {f}\quad\text{in}\quad\mathcal{D}\times\mathcal{T},
\end{equation}
where $C(p)= \frac{\partial \theta}{\partial p}$ is the specific moisture capacity. It is well-known that numerical solutions originating from the pressure head formulation may give rise to a significant mass balance error, resulting  in an inaccurate prediction of the infiltration depth. 

Numerical methods based on the mixed form (using finite differences or mass-lumped finite elements) are popular as they result in mass conservation schemes \cite{Celia1990}. Therefore, we will work with the mixed form \eqref{mixedForm} of the Richards' equation. 
 \subsection{Van Genuchten-Mualem model}\label{VGMmodel}
To complete  the PDE formulations, \eqref{mixedForm} or \eqref{pressurehead}, closure models for approximating $K_{rw}$ and $\theta$ are required. A number of models have been presented in the literature and the most popular ones are by Brooks-Corey \cite{brooks1964hydrau} and van Genuchten-Mualem \cite{vanG,Mualem}. These two models employ nonlinear constitutive relations for $K_{rw}$ and $p$, and for $\theta$ and $p$, respectively. We consider the parameterization introduced by van Genuchten and Mualem here. For the saturation, van Genuchten \cite{vanG} proposed the following analytic formula:
\begin{equation}\label{Sw}
S_w(p) =\frac{\theta(p) - \theta_r}{\theta_s - \theta_r}=
 \begin{cases}{}
(1+(|\alpha p|)^n)^{-m},\quad p<0,\\
1,\hspace{2.6cm} p\geq0,
\end{cases}
\end{equation}
where $\theta_s$ and $\theta_r$ are the saturated and residual water contents, respectively, and $\alpha [L^{-1}],n$ and $m = 1 - n^{-1}$ are obtained by fitting data characterizing the statistics of the soil.  Specifically, the parameter $\alpha$ provides a measure of the average pore-size in the soil matrix and $n$ is related to the pore-size distribution of the soil \cite{Miller1998}. %These parameters are typically fixed throughout the computational domains for simulation. 

We may derive the specific moisture content, $C(p)$, analytically from \eqref{Sw}, as
\begin{equation}
\label{sc_van}
C(p) = \left\{
\begin{array}{l}
\displaystyle (\theta_s - \theta_r) \alpha m n (1 + | \alpha p |^n)^{-(m+1)} | \alpha p |^{n-1}, \quad p < 0, \\
0, \hspace{6cm} p \geq 0.
\end{array}
\right.
\end{equation}
In previous work, Mualem \cite{Mualem} derived a closed-form expression for $K_{rw}$, which is given by:
\begin{equation}\label{Krw1}
K_{rw} = S_w^{1/2}\left[\int_0^{S_w} {dS_w}/{p} \middle/\int_0^{1} {dS_w}/{p}\right]^2.
\end{equation}
Using \eqref{Sw}, the above integral equation reduces to the following analytic expression:
 \begin{equation}\label{Krw}
K_{rw}(p) =
 \begin{cases}{}
S_w(p)^{1/2}\bigg(1 - \big(1-S_w(p)^{1/m}\big)^m\bigg)^2,\quad p<0,\\
1,\hspace{5.1cm} p\geq0.
\end{cases}
\end{equation}
The complexity of the numerical solution of the Richards' equation depends on the values of the parameters $n$ and $\alpha$. For $n\in(1,2)$ and $p\rightarrow0$, the relative hydraulic conductivity $K_{rw}(p)$ is not Lipschitz continuous and the derivative $K_{rw}'(p)$ becomes infinite as $p$ approaches zero \cite{Miller1998,Ippisch2006}. Moreover, for small values of $n$, a sharp $K_{rw}$ vs. $p$ profile is encountered.  Similarly, for large values of the parameter $\alpha$, the pressure head exhibits a transition behaviour with a steep gradient from the saturated to the unsaturated region. In general, for a small $n$ or for large $\alpha$ strong nonlinearities are encountered, thus implying convergence issues for nonlinear iterative techniques such as the Newton or Picard methods.

\section{Stochastic Richards' model}\label{stochRE}
Here, we describe a stochastic extension of the van Genuchten model. We assume that the unknown soil parameters belong to the probability space $(\Omega,\mathbb{F},\mathbb{P})$, where $\Omega$ is the sample space with $\sigma$-field $\mathbb{F}\subset 2^{\Omega}$ as a set of events and the probability measure $\mathbb{P}:\mathbb{F}\rightarrow[0,1]$. 

The stochastic extension is based on modeling the soil parameters as spatially correlated random fields in order to incorporate spatial heterogeneity. For  the saturated hydraulic conductivity, $K_s$, it is standard practice to model it as a lognormal random field, as follows,

\begin{equation}\label{RandK}
K_s(\mathbf{x},\omega) = K_s^{(bl)}(\mathbf{x}) \exp(Z(\mathbf{x},\omega)),\qquad \mathbf{x}\in \mathcal{D},\omega\in\Omega,
\end{equation}
where $K_s^{(bl)}(\mathbf{x})$ is the baseline hydraulic conductivity and $Z(\mathbf{x},\omega)$ is a zero mean Gaussian random field with a specified covariance kernel. So,
\begin{align}\label{eq:GRF}
\mathbb{E}[Z(\mathbf{x},\cdot)]  &=0,\\
\text{Cov}(Z(\mathbf{x_1},\cdot),Z(\mathbf{x_2},\cdot))&= \mathbb{E}[ Z(\mathbf{x_1},\cdot)Z(\mathbf{x_2},\cdot)],\quad \mathbf{x_1},\mathbf{x_2}\in \mathcal{D}.
\end{align}
In the present work, we consider an anisotropic Mat\'ern covariance function, $C_{\Phi}$, defined as
 \begin{equation}\label{eq:nonIsotropic}
 \begin{cases}
C_\Phi(\mathbf{x_1},\mathbf{x_2}) = \sigma_c^2\dfrac{2^{1-\nu_c}}{\Gamma(\nu_c)}\left( 2\sqrt{\nu_c}\tilde{r}\right)^{\nu_c} K_{\nu_c}\left( 2\sqrt{\nu_c}\tilde{r}\right),\\[1.5ex] 
\tilde{r}  = \sqrt{\dfrac{(x_1-x_2)^2}{\lambda^2_{cx}}+ \dfrac{(z_1-z_2)^2}{\lambda^2_{cz}}},\quad\text{with}\quad \mathbf{x_1} = (x_1,z_1),\mathbf{x_2} = (x_2,z_2).
\end{cases}
\end{equation}
Here, we denote the gamma function by $\Gamma$ and by $K_{\nu_c}$ the modified Bessel function of the second kind. The Mat\'ern function is characterized by the parameter set $\Phi=\{\nu_c,{\lambda_{cx}},{\lambda_{cz}},\sigma_c^2\}$. Parameter $\nu_c\geq 0$ defines the differentiability of $Z$, $\sigma_c^2>0$ is the marginal variance and  $\lambda_{cx}$ and $\lambda_{cz}$ are correlation lengths along x- and z-coordinates, respectively. When $\nu_c=1/2$, the Mat\'ern function corresponds to an exponential covariance function and for $\nu_c\rightarrow\infty$ to a squared exponential covariance model.  Simulating a Gaussian random field can be based on the Karhunen-Lo\'eve (KL) decomposition \cite{RF1} of $Z(\mathbf{x},\omega)$, i.e.,
\begin{equation}\label{KLexp_ch6}
Z(\mathbf{x},\omega) = \sum_{j=1}^\infty\sqrt{\lambda_j}\Psi_j(\mathbf{x})\xi_j,\qquad \xi_j\sim\mathcal{N}(0,1).
\end{equation}
Here, $\lambda_j$ and $\Psi_j$ are eigenvalues and eigenfunctions of the covariance kernel ${C}_\Phi(\mathbf{x_1},\mathbf{x_2})$, obtained from the solution of the Fredholm integral,
\begin{equation}\label{eq:Fred}
\int_\mathcal{D} C_\Phi(\mathbf{x_1},\mathbf{x_2})\Psi(\mathbf{x_1}) d\mathbf{x_1} =\lambda \Psi(\mathbf{x_2}).
\end{equation}
The sum \eqref{KLexp_ch6} represents an infinite-dimensional uncertain field with a decaying contribution of the eigenmodes. The rate of decay typically depends on the smoothness and correlation length of the covariance function. The sum is truncated at a finite number of terms, $M_{KL}$, which is usually decided by balancing the KL-truncation error with other sources of error, like the discretization and sampling errors. For Gaussian processes with small correlation lengths and large variances, typically a large number of terms is needed to include the critical eigenmodes \cite{RF1}. The evaluation of the eigenmodes in the KL-expansion is expensive as it requires solving the integral equation \eqref{eq:Fred} for each mode. In the case of stationary covariance models, fast sampling of random fields can be achieved via spectral generators that employ the FFT (Fast Fourier Transform) \cite{RF2,RF4} for the factorization of the covariance matrix. Another advantage of using these spectral methods is that they are able to simulate random fields on the sampling mesh without any bias (for example, in the case of the KL-expansion). In this article, we use the  Fast Fourier Transform moving average (FFT-MA) algorithm from \cite{Ravalec2000} to sample Gaussian random fields, see \ref{appendix_A1} for details.

\subsection{Sampling of non-Gaussian random fields}\label{nonGaussian}
For sampling the van Genuchten parameters, $\alpha(\mathbf{x},\omega),n(\mathbf{x},\omega),\theta_s(\mathbf{x},\omega),\theta_r(\mathbf{x},\omega)$ in Section \ref{VGMmodel}, we employ random fields with \emph{non-Gaussian marginal} distributions. This choice of distributions is practical as these parameters can only take values in a certain range, see e.g \cite{Ippisch2006}. We introduce stochasticity in the parameters via an additive noise,
\begin{equation}\label{alphaRF}
\alpha(\mathbf{x},\omega) = \alpha^{(bl)}(\mathbf{x}) + \varepsilon_\alpha(\mathbf{x},\omega),
\end{equation}
where $\alpha^{(bl)}(\mathbf{x}) $ is the deterministic baseline value and $\varepsilon_\alpha(\mathbf{x},\omega)$ is a random field with a non-Gaussian marginal distribution and covariance $C_\Phi$. Notations are analogously for the other three van Genuchten parameters. Next, we describe a technique proposed in \cite{sakamoto} for the point-wise transformation of a standard Gaussian random field to a non-Gaussian random field. 

Non-Gaussian random fields are difficult to simulate as they are not uniquely determined by their mean and variance. There are however different techniques available for simulating non-Gaussian fields, see e.g. \cite{sakamoto,PHOON2005}. In this work, we will follow a basic approach based on a generalized Polynomial Chaos (gPC) expansion \cite{sakamoto}, which approximates the non-Gaussian field in terms of a weighted combination of Hermite orthogonal polynomials of the standard Gaussian field, 
\begin{equation}\label{gPC}
%Y_i(\bfx,\omega)\approx\sum^{N_{PC}}_{n=1} \text{w}_n(\bfx) \mathcal{H}_n(Z(\bfx,\omega)) ,
Y(\mathbf{x},\omega)\approx\sum^{N_{PC}}_{j=0} \text{w}_j \mathcal{H}_j(Z(\mathbf{x},\omega)) ,
\end{equation}
where $Y(\mathbf{x},\omega)$ is the non-Gaussian random field (with a marginal distribution, e.g. the uniform distribution, gamma distribution, truncated normal, etc). $\mathcal{H}_j(Z)$ is the Hermite polynomial in $Z$ of order $j$ with weight $\text{w}_j$ and $N_{PC}$ is the order of the expansion. Hermite polynomials can be expressed as:
\begin{eqnarray}\label{hermite}
&\mathcal{H}_0(Z) =1, &\mathcal{H}_j(Z) = (-1)^j\exp(Z^2/2)\frac{d^j}{dx^j}\exp(Z^2/2),\quad j\in \mathbb{N}.
\end{eqnarray}
As Hermite polynomials are orthogonal with respect to the Gaussian measure, the weights can be evaluated using
\begin{equation}\label{weight}
\text{w}_j  = \frac{\mathbb{E}[Y\mathcal{H}_j(Z)]}{\mathbb{E}[\mathcal{H}_j(Z)^2]}.
\end{equation}
Here, the denominator is basically an expectation of a polynomial of the Gaussian random variable, which has an analytic expression. As the dependence between $Y$ and $Z$ is unknown, the expectation in the numerator is not well-defined. Since the cumulative distribution for $Y$, defined as $F_{Y}(y) = \mathbb{P}\text{rob}(Y\leq y)$, is however known, one can utilize the relation $Y = F_Y^{-1}(F_Z(Z))$ to reformulate \eqref{weight} as

\begin{equation}\label{weight2}
\text{w}_j = \frac{1}{\mathbb{E}[\mathcal{H}_j(Z)^2]}\int_{I_Z} F_{Y}^{-1}[F_{Z}(z)]\mathcal{H}_j(z) \text{d}F_Z(z),
\end{equation}
where $I_Z$ is the support of $Z$ in the range $(-\infty,\infty)$ and $F_Y^{-1}$ representing the inverse of the distribution $F_Y$. Similarly, $F_{Z}(z) = \mathbb{P}\text{rob}(Z\leq z)$ is the cumulative distribution for standard Gaussian random variable $Z$. The integral \eqref{weight2} can be numerically computed using any conventional integration technique, for example, by using Monte Carlo quadrature. %With above weights, the gPC expansion in \eqref{gPC} converges to $Y$ in \emph{weak sense} (or convergence in probability distribution) \cite{Dxiu2,Xiu:2010}. 
The weights only need to be computed once, so that the cost of sampling a non-Gaussian random field with a stationary covariance function is of the same order as that of a Gaussian random field. 

\begin{figure}[hbt]
\begin{subfigure}[b]{0.49\textwidth}
\includegraphics[scale = 0.32]{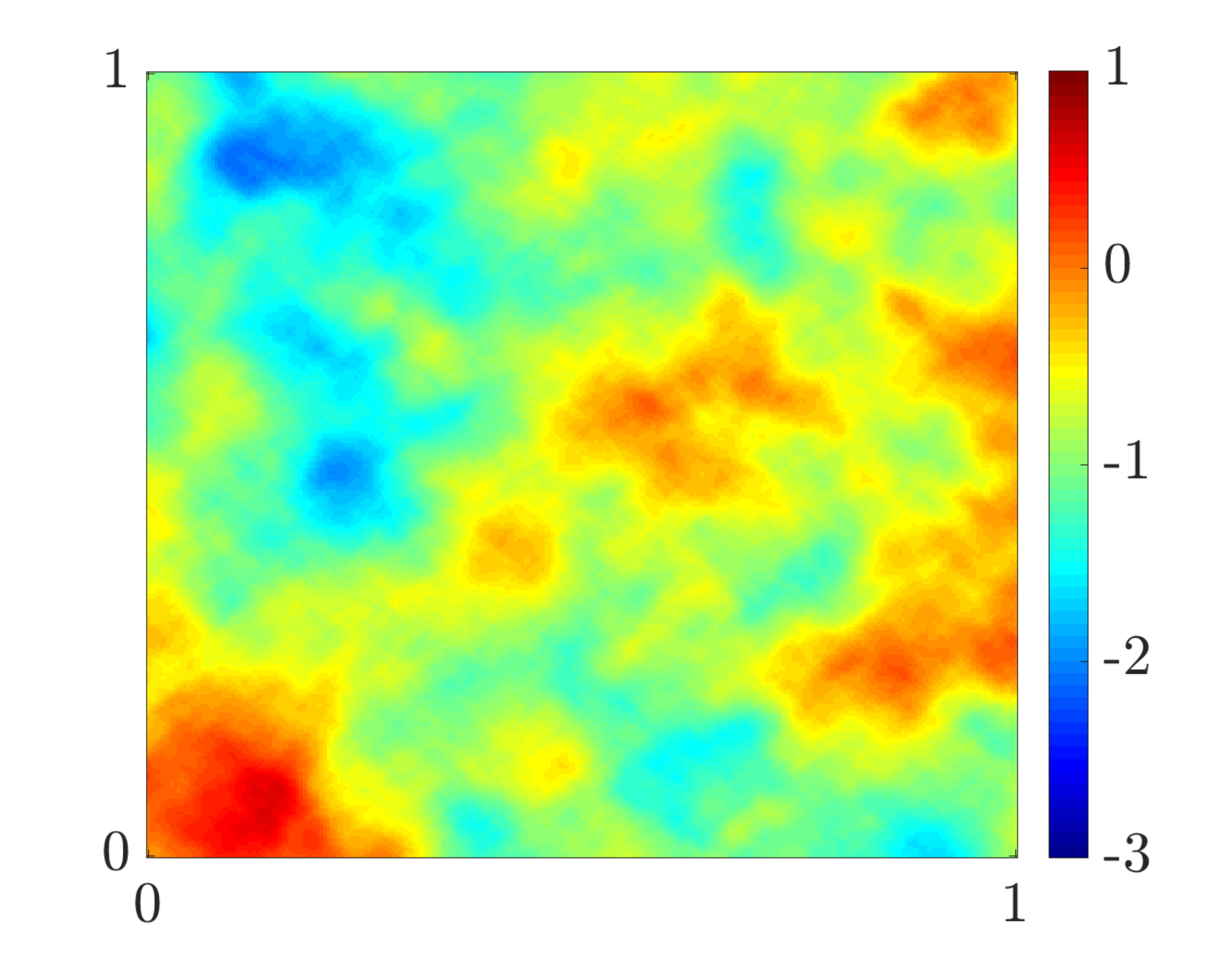}
\subcaption{$\mathcal{N}(0,C_{\Phi_1})$}
\end{subfigure}
\begin{subfigure}[b]{0.49\textwidth}
\includegraphics[scale = 0.32]{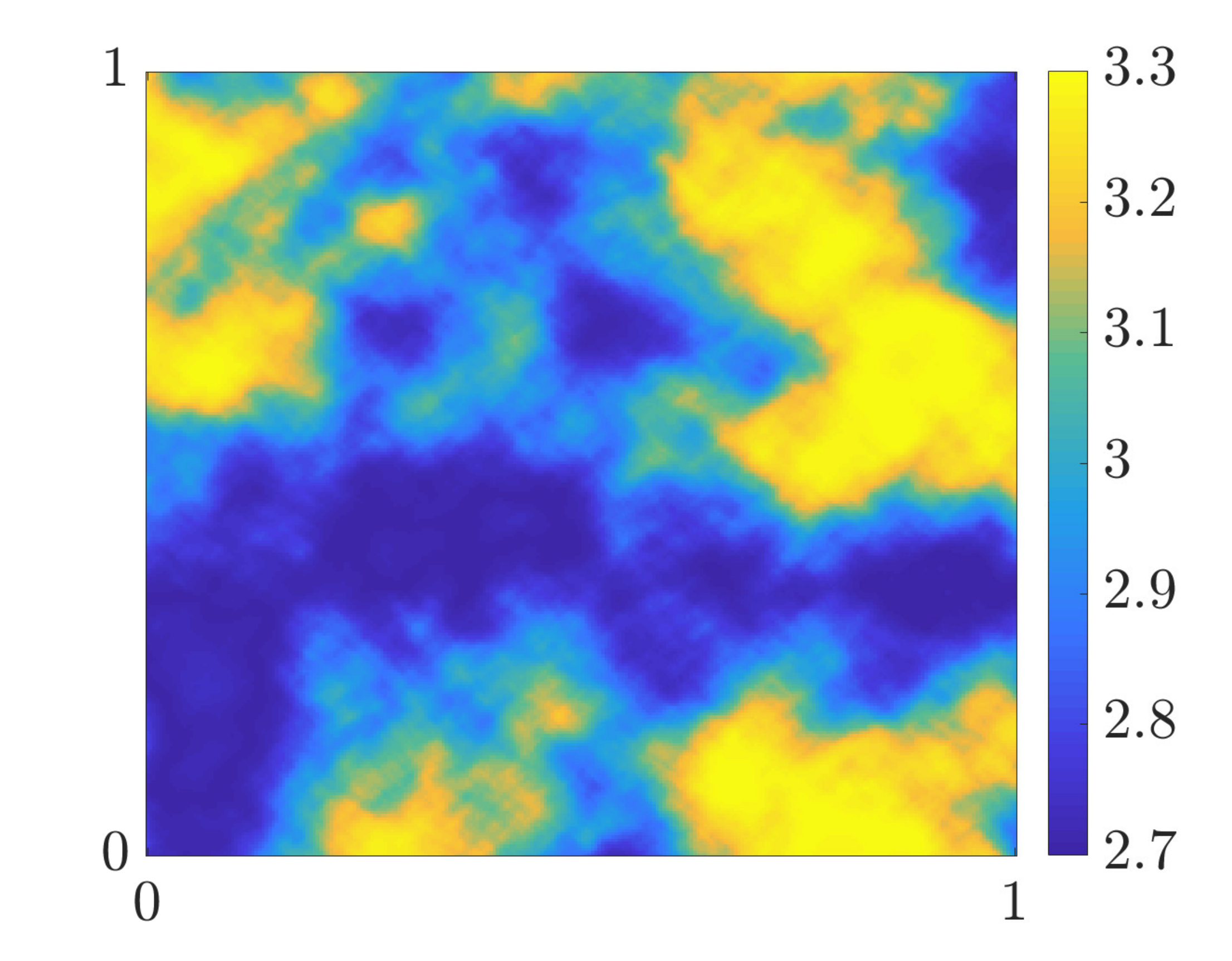}
\subcaption{$\mathcal{U}(2.7 ,3.3,C_{\Phi_1})$}
\end{subfigure}
\begin{subfigure}[b]{0.49\textwidth}
\includegraphics[scale = 0.32]{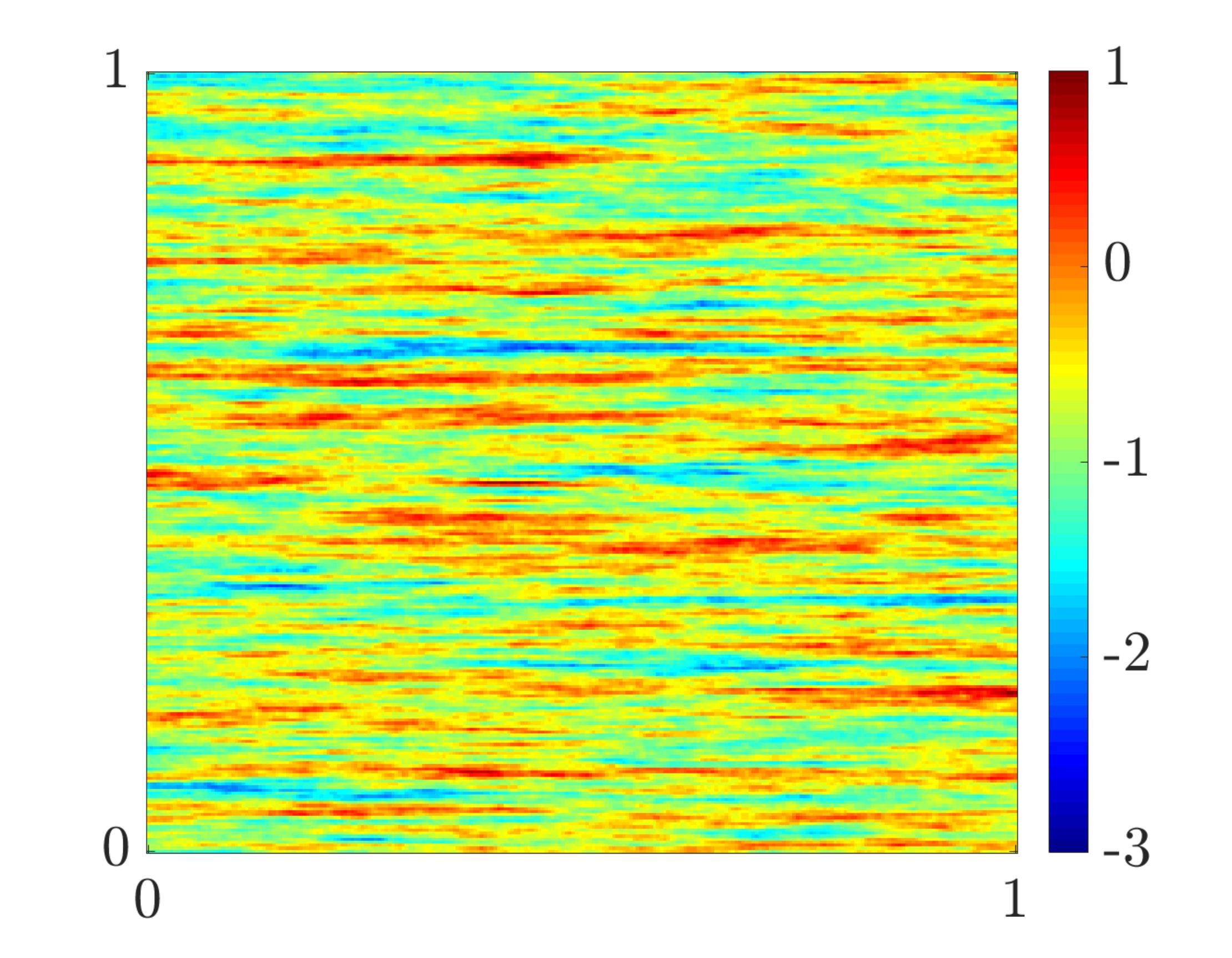}
\subcaption{$\mathcal{N}(0,C_{\Phi_2})$}
\end{subfigure}
\begin{subfigure}[b]{0.49\textwidth}
\hbox{\hspace{0.2cm}\includegraphics[scale = 0.32]{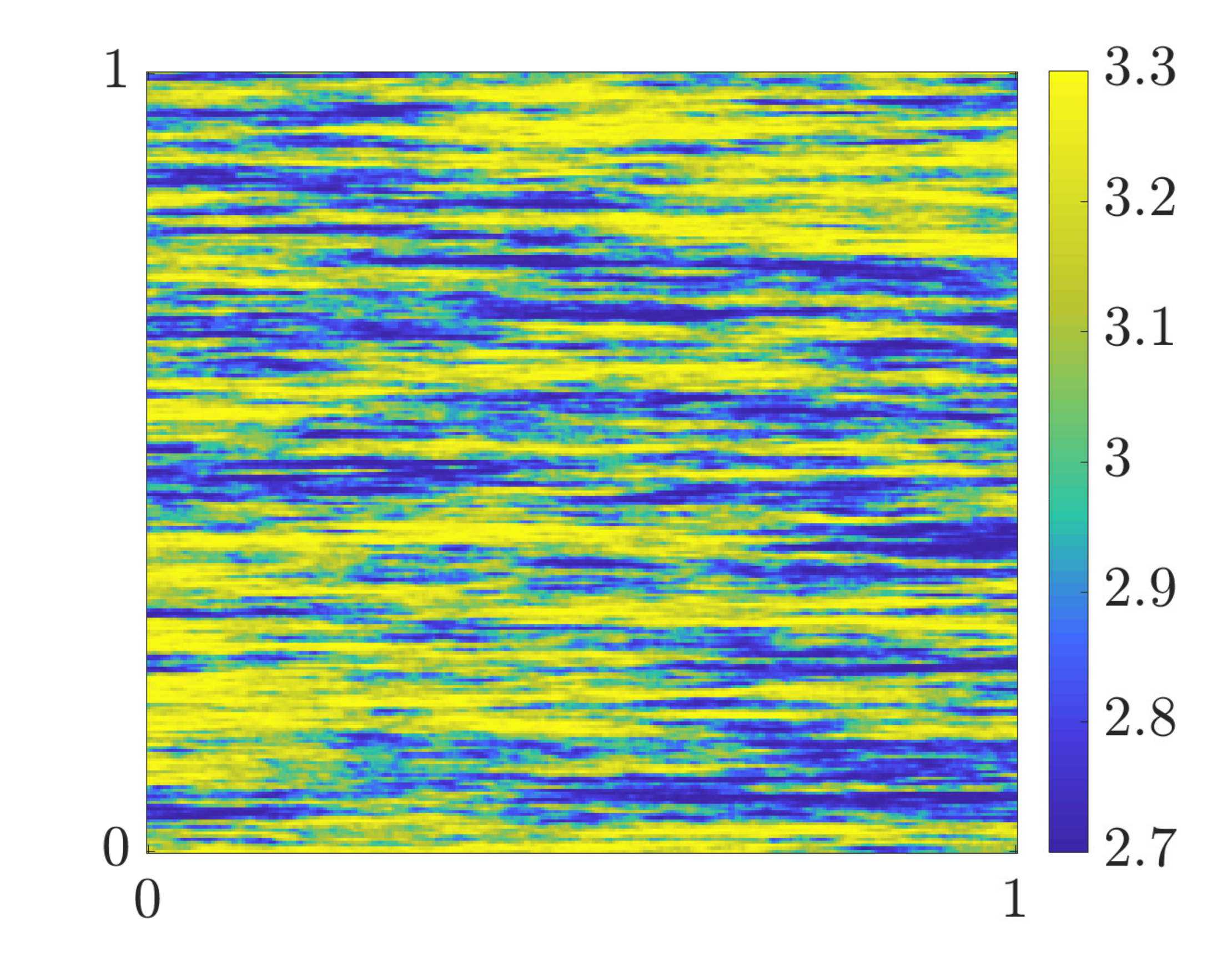}}
\subcaption{$\mathcal{U}(2.7 ,3.3,C_{\Phi_2})$}
\end{subfigure}
\caption{Samples of random fields generated using the isotropic Mat\'ern parameter $\Phi_1$ (top row) and the anisotropic parameter $\Phi_2$ (bottom row); and standard normal marginal distribution (left column) and uniform marginal distribution (right column). The notation $\mathcal{U}(2.7,3.3,C_{\Phi})$ represents a random field with uniform marginal distribution, $\mathcal{U}[2.7,3.3]$, with spatial correlation defined by $C_{\Phi}$.}\label{sampleRF}
\end{figure}

We will experiment here with both isotropic and anisotropic Mat\'ern covariance models. In Table \ref{Matern}, the two Mat\'ern parameters sets are listed, $\Phi_1$ corresponding to an isotropic model and $\Phi_2$ to an anisotropic model. In Figure \ref{sampleRF}, we present some samples of the random fields with a Gaussian and a uniform marginal distribution, for the two Mat\'ern parameters. We use $N_{PC} = 6$ in Equation~\eqref{gPC} for generating the random fields with uniform marginal distribution. Due to a small correlation length and low spatial regularity, the numerical solutions of the PDE with random coefficients based on $\Phi_2$ are comparatively more expensive to compute than those obtained with $\Phi_1$. A comprehensive study on the computational cost of solving elliptic PDEs with different Mat\'ern parameters can be found in \cite{prashant3}. We will study the effect of covariance functions on the performance of the solver for the Richards' equation.
\begin{table}[H]
\begin{center} 
\begin{tabular}{ccc}
\hline
${\Phi_1}$ & \hspace{1cm} & ${\Phi_2}$\\
\hline
(1.0,0.2,0.2,1)&\hspace{1cm}&$(0.5,0.1,0.01,1)$ \\
\hline
\end{tabular}
\end{center}
\caption{Two combinations of Mat\'ern parameters $\Phi=(\nu_c,{\lambda_{cx}},{\lambda_{cz}},\sigma_c^2)$ corresponding to isotropic ($\Phi_1$) and anisotropic ($\Phi_2$) random fields.}\label{Matern}
\end{table}
\section{Modified Picard iteration combined with the cell-centered multigrid method}\label{MPCCMG}

Algorithms based on the modified Picard iteration from Celia et al. \cite{Celia1990} are often employed as efficient iterative solution methods for the Richards' equation. These methods are relatively easy to implement, as they do not require the computation of Jacobians and they also have low storage requirements. Within each modified Picard iteration, a diffusion equation with variable coefficients needs to be solved. For this, we propose to utilize the \emph{cell-centered multigrid} (CCMG) for heterogeneous diffusion coefficients, as proposed in \cite{prashant3,WESSELING198885,MOLENAAR199625}. The CCMG algorithm is efficient as it is constructed with a simple set of transfer operators and it has been demonstrated to perform well for a large class of highly heterogeneous and also jumping diffusion coefficients \cite{prashant3}. % Later on, we will numerically show that the resulting modified Picard-CCMG solver is robust for a large range of soil parameters, making it ideal for multilevel Monte Carlo methods.

\subsection{Modified Picard iteration}

We briefly recall the fully-implicit Picard iteration for the mixed formulation of the Richards' equation from \cite{Celia1990}. With $\Delta t$ the time-step and for any integer $J>1$, we define a uniform temporal grid by $\{t^j = j\Delta t, j = 0, \ldots, J\}$. The iteration number within a time-step is denoted by an integer $k>0$. For simplicity, we use a simplified notation for $\theta^{j,k} = \theta(p^{j,k})$ and ${K}^{j,k} = {K}_sK_{rw}(p^{j,k})$. The backward Euler approximation of \eqref{mixedForm} is then written as

\begin{equation}
\label{Mixed_based}
\frac{\theta^{j+1,k+1} - \theta^j}{\Delta t} - \nabla \cdot {K}^{j+1,k} \nabla p^{j+1,k+1} - \frac{\partial K}{\partial z}^{j+1,k}=  {f}^{j+1}.
\end{equation}

The key idea of the modified Picard iteration is the use of a Taylor expansion for $\theta^{j+1,k+1}$ with respect to $p$, i.e. 

\begin{equation}\label{taylorThetap}
\theta^{j+1,k+1} = \theta^{j+1,k} + \frac{\partial \theta}{\partial p}^{j+1,k}(p^{j+1,k+1}-p^{j+1,k}) + \mathcal{O}\left(\delta p^2\right),
\end{equation}
 where the derivative $\frac{\partial \theta(p)}{\partial p}  = C(p)$ is analytically computed by using \eqref{sc_van}. By neglecting the higher-order terms in \eqref{taylorThetap} and substitution in \eqref{Mixed_based}, we get
\begin{equation}
\label{Mixed_based_2}
C(p^{j+1,k}) \frac{\delta p^{j+1,k}}{\Delta t}  + \frac{\theta^{j+1,k} - \theta^j}{\Delta t} - \nabla \cdot {K}^{j+1,k} \nabla  p^{j+1,k+1} - \frac{\partial K}{\partial z}^{j+1,k} =  {f}^{j+1},
\end{equation}
with  $\delta p^{j+1,k} = p^{j+1,k+1} - p^{j+1,k}$. The above equation can be expressed in the form:
\begin{equation}
\label{Mixed_based_3}
C(p^{j+1,k}) \frac{\delta p^{j+1,k}}{\Delta t}  - \nabla \cdot {K}^{j+1,k} \nabla \delta p^{j+1,k} = \nabla \cdot {K}^{j+1,k} \nabla p^{j+1,k} + \frac{\partial K}{\partial z}^{j+1,k} + {f}^{j+1} - \frac{\theta^{j+1,k} - \theta^j}{\Delta t}.
\end{equation}

The next pressure head iterate is obtained by the update $p^{j+1,k+1} = p^{j+1,k}+\delta p^{j+1,k}$. Notice that the left-hand side of the above equation is the residual associated with the Picard iteration, which should be equal to  zero for a converged solution. Therefore, one may use $||\delta p^{j+1,k} ||_{\infty}<\varepsilon_{PI}$ as a stopping criterion with $\varepsilon_{PI}>0$ as the tolerance for Picard iteration.  The pressure head at time $t^{j+1}$ is then given by $p^{j+1} = p^{j+1,k+1}$, with $k$ the total number of Picard iterations to converge to $\varepsilon_{PI}$. The iterative scheme \eqref{Mixed_based_3} is a general mixed-formulation Picard iteration, which results in perfect mass balance. 

\subsection{Cell-centered multigrid}

Focussing on the $k$-th Picard iteration \eqref{Mixed_based_3} at time $t^{j+1}$, the following elliptic PDE with variable coefficients is obtained, using simplified notation, 

\begin{eqnarray}\label{diffusion_problem}
\frac{\widetilde{C}}{\Delta t}\delta \widetilde{p} -\nabla \cdot \left(\widetilde{K} \nabla \delta \widetilde{p} \right) & = & {\tilde f} \hspace{0.25cm} \text{in}\quad {\cal D}, \\
\delta \widetilde{p} & =&  0 \quad \text{in} \quad {\Gamma_D}\cup{\Gamma_N},\nonumber
%\delta \widetilde{p} & = & g_N \hspace{0.1cm}\text{in}\quad {\Gamma_N},
\end{eqnarray}
with the known quantities
\begin{equation}
\widetilde{C} = C(p^{j+1,k}),\quad \widetilde{K} = K^{j+1,k}\quad\text{and}\quad \widetilde{f} =  \nabla \cdot {K}^{j+1,k} \nabla p^{j+1,k} + \frac{\partial K}{\partial z}^{j+1,k} + {f}^{j+1} - \frac{\theta^{j+1,k} - \theta^j}{\Delta t},\nonumber
\end{equation}
and the unknown $\delta \widetilde{p} = \delta p^{j+1,k}$. To discretize the above problem, we use a {\em cell-centered finite volume scheme} for which the hydraulic conductivity at the cell-face is based on the harmonic averaging of the hydraulic conductivities from the adjacent cells, derived by the continuity of fluxes \cite{WESSELING198885,MOLENAAR199625}. 

For the discretization of \eqref{diffusion_problem}, a uniform grid $\mathcal{D}_h$ on a unit square domain with the same mesh width $h = 1/M, M \in \mathbb{N}$ in both directions, 
\begin{equation}
\label{Grid_fine}
\mathcal{D}_h = \left\{(x_{i_1},z_{i_2}); x_{i_1} = \left(i_1-\frac12\right) h,z_{i_2} = \left(i_2-\frac12\right) h, i_1,i_2 = 1,\ldots,M\right\}, 
\end{equation}
is considered. For each \emph{interior cell} (edges do not lie on a boundary) with center $(x_{i_1},z_{i_2})$, denoted by $\mathcal{D}_h^{i_1,i_2}$, we obtain a five-point scheme,

\begin{equation}
c^h_{i_1,i_2} \delta \widetilde{p}_{i_1,i_2} + w^h_{i_1,i_2} \delta \widetilde{p}_{i_1-1,i_2} + e^h_{i_1,i_2} \delta \widetilde{p}_{i_1+1,i_2} + s^h_{i_1,i_2} \delta \widetilde{p}_{i_1,i_2-1} + n^h_{i_1,i_2} \delta \widetilde{p}_{i_1,i_2+1} = \widetilde{f}^h_{i_1,i_2},
\label{five_stencil}
\end{equation}
with 
\begin{eqnarray}
w^h_{i_1,i_2} &=& -\frac{2}{h^2} \frac{\widetilde{K}_{i_1,i_2} \widetilde{K}_{i_1-1,i_2}}{\widetilde{K}_{i_1,i_2} + \widetilde{K}_{i_1-1,i_2}}, \quad
e^h_{i_1,i_2} = -\frac{2}{h^2} \frac{\widetilde{K}_{i_1,i_2} \widetilde{K}_{i_1+1,i_2}}{\widetilde{K}_{i_1,i_2} + \widetilde{K}_{i_1+1,i_2}}, 	\nonumber \\
s^h_{i_1,i_2} &=& -\frac{2}{h^2} \frac{\widetilde{K}_{i_1,i_2} \widetilde{K}_{i_1,i_2-1}}{\widetilde{K}_{i_1,i_2} + \widetilde{K}_{i_1,i_2-1}}, \quad
n^h_{i_1,i_2} = -\frac{2}{h^2} \frac{\widetilde{K}_{i_1,i_2} \widetilde{K}_{i_1,i_2+1}}{\widetilde{K}_{i_1,i_2} + \widetilde{K}_{i_1,i_2+1}}, \nonumber \\
c^h_{i_1,i_2} &=& -(w^h_{i_1,i_2} + e^h_{i_1,i_2} + n^h_{i_1,i_2}+ s^h_{i_1,i_2}) +\frac{\widetilde{C}_{i_1,i_2}}{\Delta t} \nonumber,
\end{eqnarray}
where, for instance, $\widetilde{K}_{i_1,i_2}$ is the diffusion coefficient associated with cell $\mathcal{D}_h^{i_1,i_2}$ and the source term $\widetilde{f}^h_{i_1,i_2}$ is an approximation of $\widetilde{f}$ in that cell. This scheme is modified appropriately for cells close to the boundary. 

\begin{figure}[H]
\includegraphics[trim={0cm 5cm 0cm 6.5cm},clip,scale=0.62]{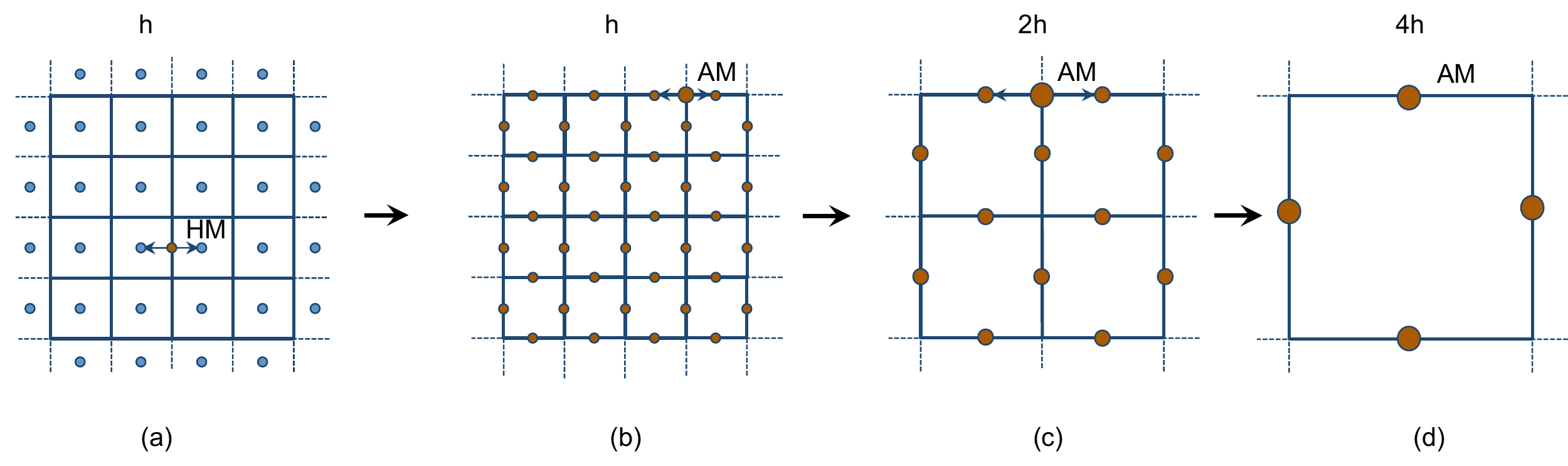}
\caption{Schematic representation of coefficient upscaling in the multigrid hierarchy (h-2h-4h). (a) Coefficient values at cell-centres (blue dots) (b) Values at face-centres (red dots) obtained from the harmonic mean (HM) of coefficients from two adjacent cell-centres. (c) Values at face-centres (bigger red dots) of the 2h-grid are based on arithmetic mean (AM) of coefficients from face-centres of the h-grid. (d) Values at face-centres (biggest red dots) of the 4h-grid are based on the arithmetic mean of the values of the coefficients from face-centres of the 2h-grid.}\label{diffUpscale}
\end{figure}

Next, we describe the multigrid method for solving the linear system arising from the above discretization. The multigrid hierarchy is based on uniform grid coarsening, i.e. the cell-width is doubled in each coarsening step in each direction.  As the smoothing method, we use the lexicographic Gauss-Seidel iteration, and as the transfer operators between the fine and coarse grids a simple piece-wise constant prolongation operator, $P_{2h}^h$, is applied and its scaled adjoint is used as the restriction operator $R_h^{2h}$ on the cell-centered grid. In classical stencil notation, these are written as,
\begin{equation}\label{eq:transfer}
P_{2h}^h = 
\left]
\begin{array}{ccc}
1 &  & 1 \\ 
& \star & \\ 
1 & & 1
\end{array}
\right[_{2h}^{h}, \ \hspace{2cm}
R_h^{2h}=\frac{1}{4}\left[ \begin{array}{ccc}
1 &  & 1\\  &\star&  \\ 1 & & 1
\end{array}\right]^{2h}_h,
\end{equation}
respectively, where $\star$ denotes the position of the cell center. The coarse grid operator is obtained via a direct discretization of the PDE operator on the coarse grid. For this discretization on a coarser grid, we need to appropriately define the diffusion coefficients on the coarse cell edges. The technique to define the suitable diffusion coefficient on a coarse cell edge is graphically described in Figure \ref{diffUpscale} and its caption. In \cite{prashant3}, the $W(2,2)$-cycle was found to be a very robust and efficient multigrid cycling strategy, and, therefore, we also employ this cycle in our experiments. The number of multigrid iterations is based on the stopping criterion, $||\mathcal{L}_h\delta \tilde{p}_h - \tilde{f}_h||_{\infty}<\varepsilon_{MG}$, where $\mathcal{L}_h$ denotes the linear operator after the discretization of Equation~\eqref{diffusion_problem} and $\varepsilon_{MG}>0$.
We consider the modified Picard method in this article as it is widely adopted, although many modifications have been proposed to improve its robustness. For instance, the authors in \cite{MILLER2006525} studied a spatio-temporal adaptive solution method to improve the numerical stability of the modified Picard iteration. Another interesting improvement was proposed in \cite{LOTT2012}, where an Anderson acceleration was applied to improve the robustness and computational cost for the standard Picard iteration scheme. These improvements can easily be extended to the modified Picard-CCMG solver studied here. Also, there are a number of effective solution approaches based on Newton's method, see for e.g. \cite{Steffen,JONES2001763,Juncu2012}. These methods exhibit a quadratic convergence rate but are very sensitive to initial solution approximations. 
\subsection{Performance of the modified Picard-CCMG solver}\label{numexp1}
We study the performance of the modified Picard-CCMG solver for a range of values of the parameters $\alpha,n$ and the effect of the heterogeneity of the hydraulic conductivity on the performance of the solver. For this we consider an \emph{infiltration problem} \cite{Juncu2012,Caputo2008} on a two-dimensional computational domain $\mathcal{D}=(0,1)^2$. The initial and boundary conditions are prescribed as follows:
\begin{eqnarray}\label{travellingwave}
p(x,z,0) = -0.4(1-\exp(-80z)),\hspace{1cm}p(x,1,t) &=& -0.4,\\ \nonumber
p(x,0,t) = 0.1,\hspace{1cm} \frac{\partial p}{\partial x}\bigg |_{x=0,1} &=& 0.
%K_s^{(bl)} = 0.2 \text{ [m/h]},\hspace{1cm} \theta_s^{(bl)} &=&0.50,\hspace{1cm} \theta_r^{(bl)} = 0.05.
\end{eqnarray}
The right-hand side is assumed to be zero, and we consider the a final time $T_{final} =0.1\text{ [h]}$. In Table \ref{alpha_n}, we provide a list of 20 values for $\alpha$ and $n$, used in the experiments. In total, we test 400 pairings of $\alpha$ and $n$. Parameters $\theta_s = 0.50$ and $\theta_r  = 0.05$ are fixed, as they do not pose any problems for the convergence rate of the solver. The samples of hydraulic conductivity are generated according to \eqref{RandK}, with $K_s^{(bl)} =  0.2 \text{ [m/h]}$ and the covariance is based on the two Mat\'ern parameters from Table \ref{Matern}.  
\begin{table}[H]
\label{fourMat}
\begin{center} 
\begin{tabular}{ccccc|ccccc}
\hline
$\alpha$&&&&&$n$&&&&\\
\hline
0.2 &0.4 &0.6&0.8&1.0&1.1&1.2&1.3&1.4&1.5\\
1.2 &1.4 &1.6&1.8&2.0&1.6&1.7&1.8&1.9&2.0\\
2.2 &2.4 &2.6&2.8&3.0&2.2&2.4&2.6&2.8&3.0\\
3.2 &3.4 &3.6&3.8&4.0&3.2&3.4&3.6&3.8&4.0\\
\hline
\end{tabular}
\end{center}
\caption{Set of $\alpha,n$ values used to for benchmarking the modified Picard-CCMG solver.}\label{alpha_n}
\end{table}
A similar test was performed in \cite{LOTT2012} for a deterministic steady-state flow governed by the Richards' equation. We perform our experiments in a probabilistic framework. For a given pair $\alpha,n$, we generate 64 random hydraulic conductivity fields and solve Richards' equation with conditions given in \eqref{travellingwave} for each sample. This is done as the number of multigrid iterations varies depending on the random realization of the hydraulic conductivity field. The cost of solving one instance of stochastic Richards' equation is expressed in terms of the total number of multigrid $W(2,2)-$cycles needed to solve the time-dependent problem. Here, by total number of $W(2,2)-$ cycles means the sum of multigrid iterations needed to reach $T_{final}$. For all experiments, we set the tolerances $\varepsilon_{PI},\varepsilon_{MG} = 10^{-5}$. The solution method was terminated with failure when the maximum number of nonlinear iterations (set to 50) was exceeded at any time-step.

\begin{figure}
\begin{subfigure}[b]{0.49\textwidth}
\includegraphics[trim={1.5cm 0cm 2cm 1cm},clip, scale = 0.2]{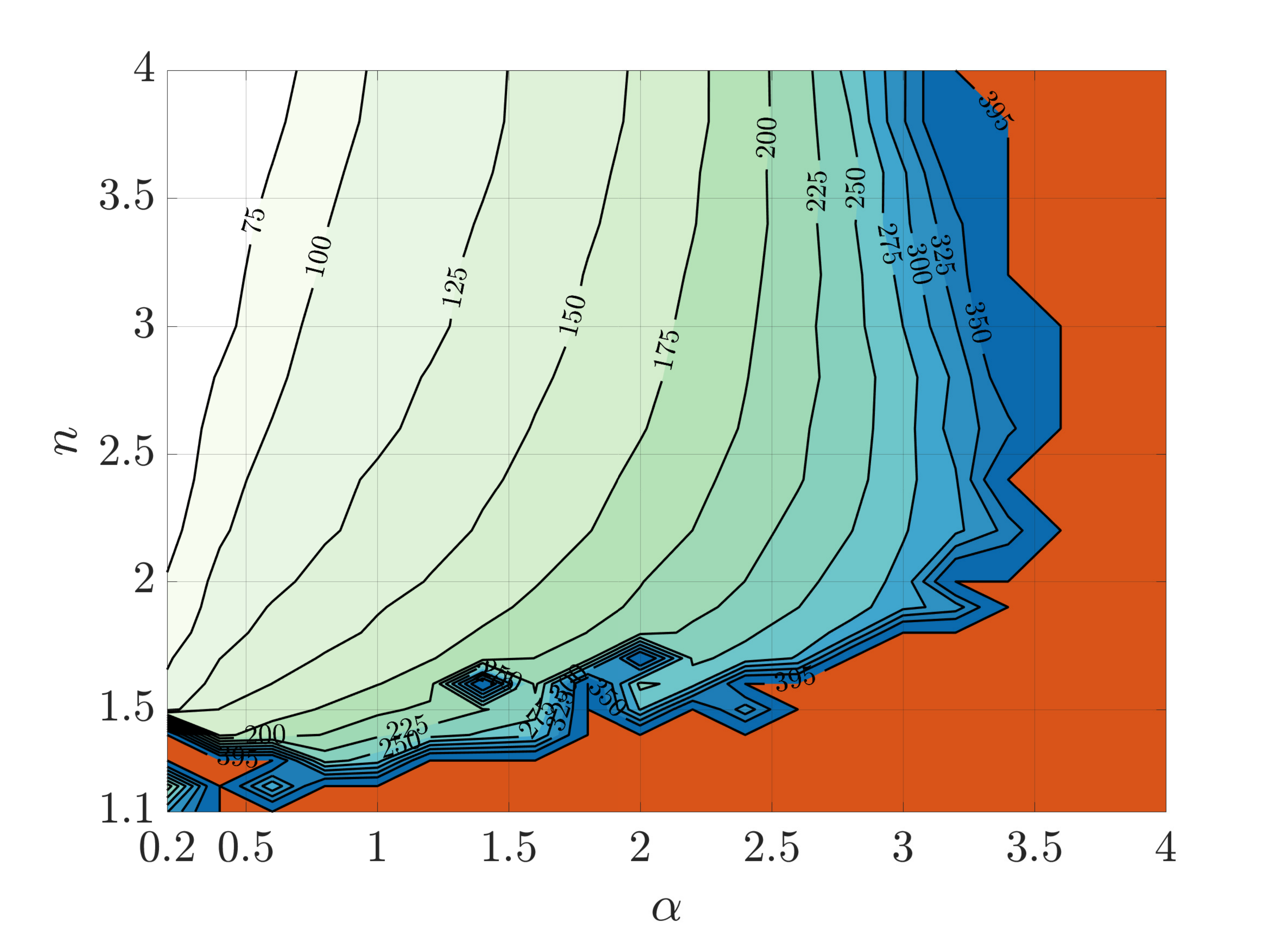}
\subcaption{$h=1/32,\Delta t = 1/64$}
\end{subfigure}
\begin{subfigure}[b]{0.49\textwidth}
\includegraphics[trim={1.5cm 0cm 2cm 1cm},clip, scale = 0.2]{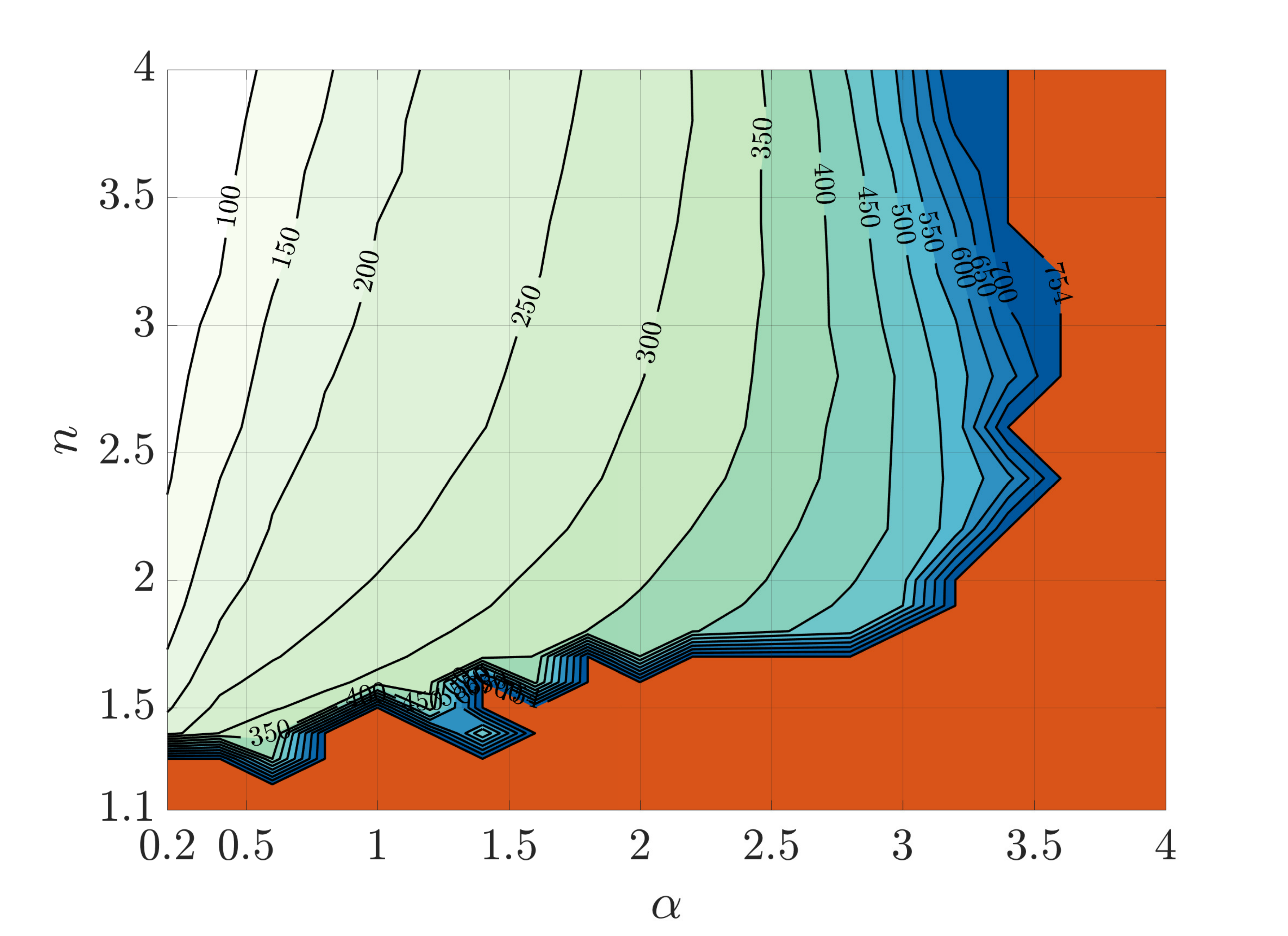}
\subcaption{$h=1/32,\Delta t = 1/128$}
\end{subfigure}
\begin{subfigure}[b]{0.49\textwidth}
\includegraphics[trim={1.5cm 0cm 2cm 1cm},clip, scale = 0.2]{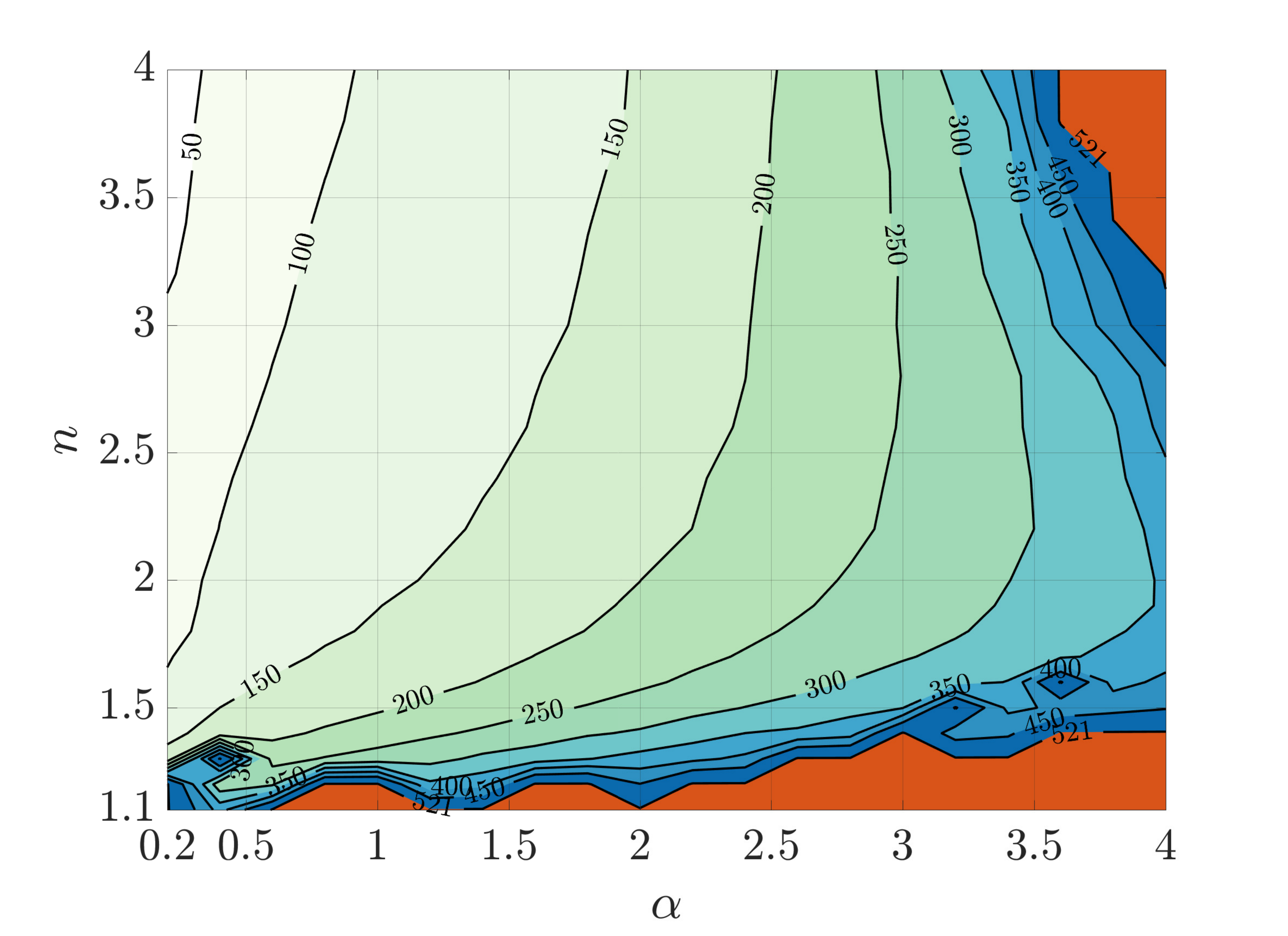}
\subcaption{$h=1/64,\Delta t = 1/64$}
\end{subfigure}
\begin{subfigure}[b]{0.49\textwidth}
\includegraphics[trim={1.5cm 0cm 2cm 1cm},clip, scale = 0.2]{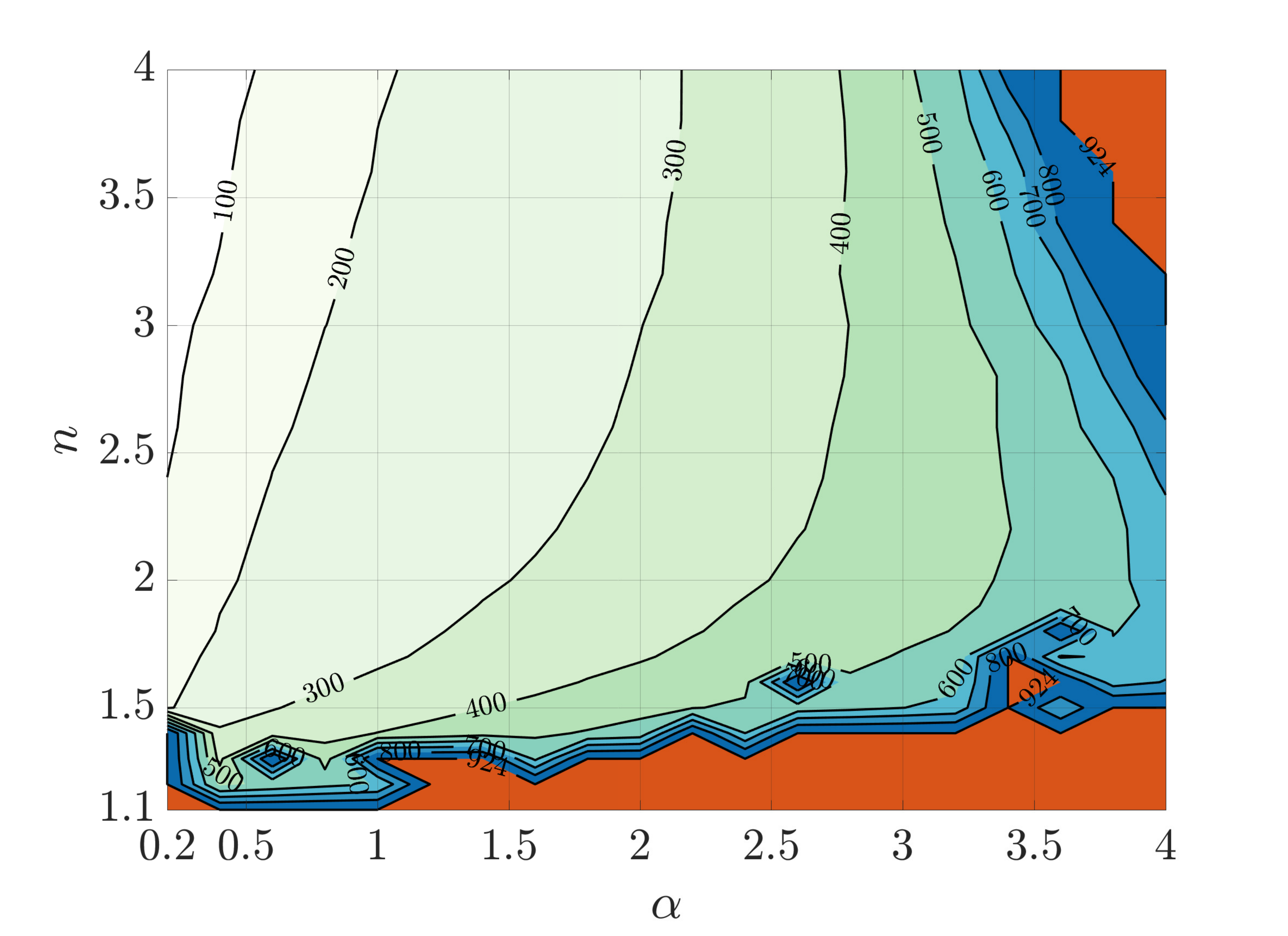}
\subcaption{$h=1/64,\Delta t = 1/128$}
\end{subfigure}
\caption{Contour plots of  the average number of multigrid iterations needed to solve the infiltration problem using the modified Picard-CCMG solver for an isotropic hydraulic conductivity field generated using $\Phi_1$.}\label{HM_Phi1}
\end{figure}

\begin{figure}
\begin{subfigure}[b]{0.49\textwidth}
\includegraphics[trim={1.5cm 0cm 1.5cm 1cm},clip,scale = 0.2]{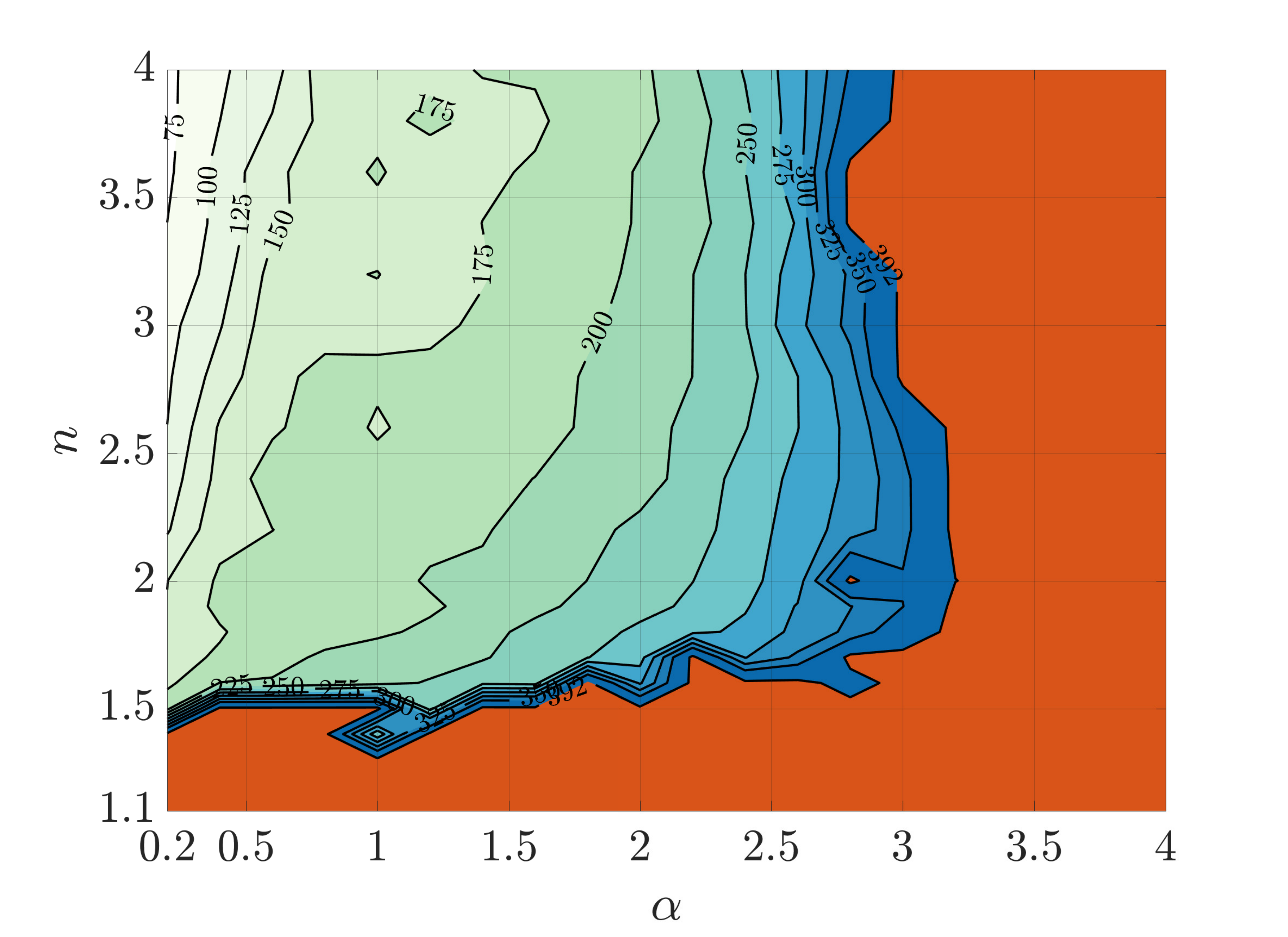}
\subcaption{$h=1/32,\Delta t = 1/64$}
\end{subfigure}
\begin{subfigure}[b]{0.49\textwidth}
\includegraphics[trim={1.5cm 0cm 2cm 1cm},clip,scale = 0.2]{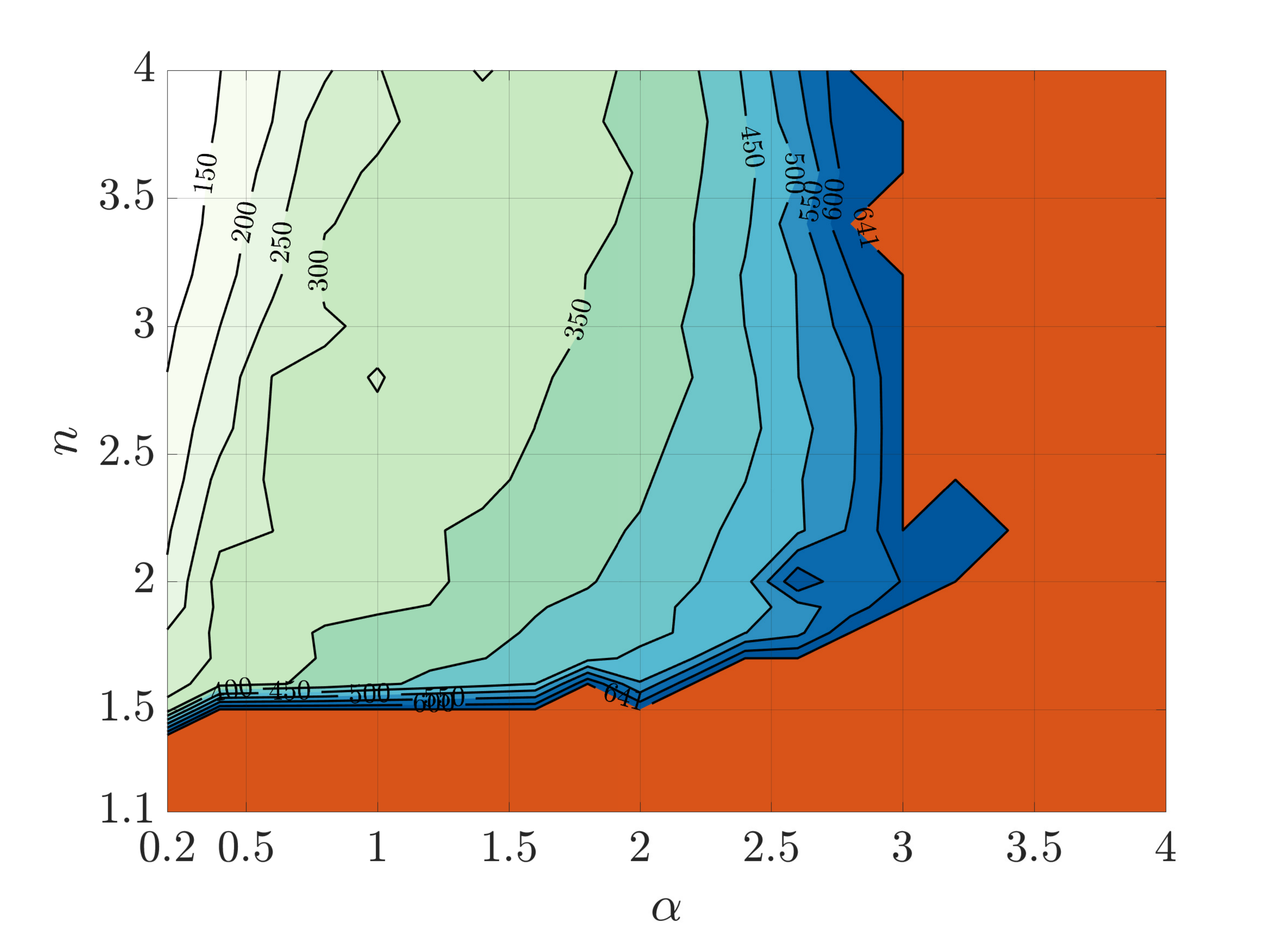}
\subcaption{$h=1/32,\Delta t = 1/128$}
\end{subfigure}
\begin{subfigure}[b]{0.49\textwidth}
\includegraphics[trim={1.5cm 0cm 2cm 1cm},clip,scale = 0.2]{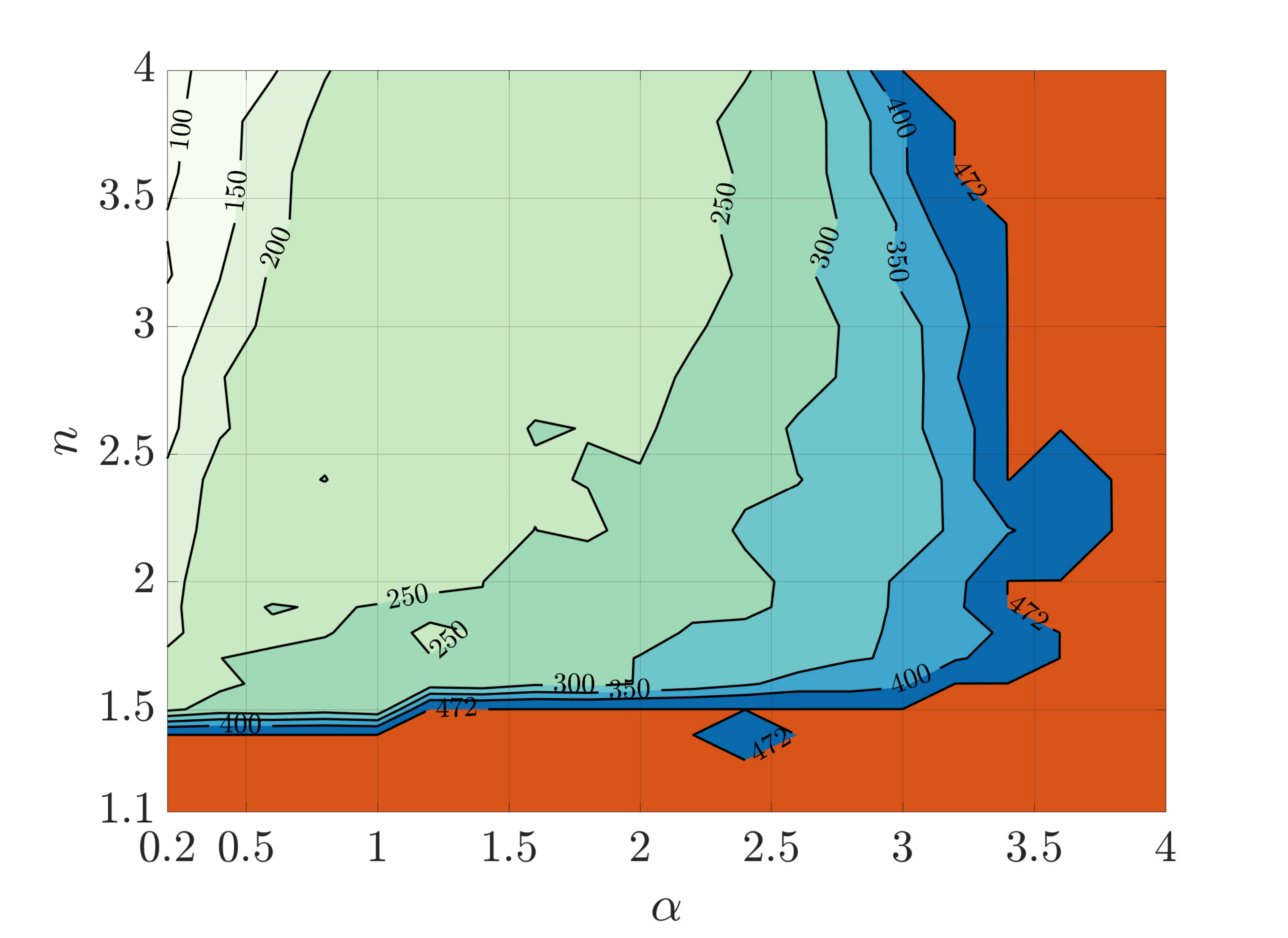}
\subcaption{$h=1/64,\Delta t = 1/64$}
\end{subfigure}
\begin{subfigure}[b]{0.49\textwidth}
\includegraphics[trim={1.5cm 0cm 2cm 1cm},clip,scale = 0.2]{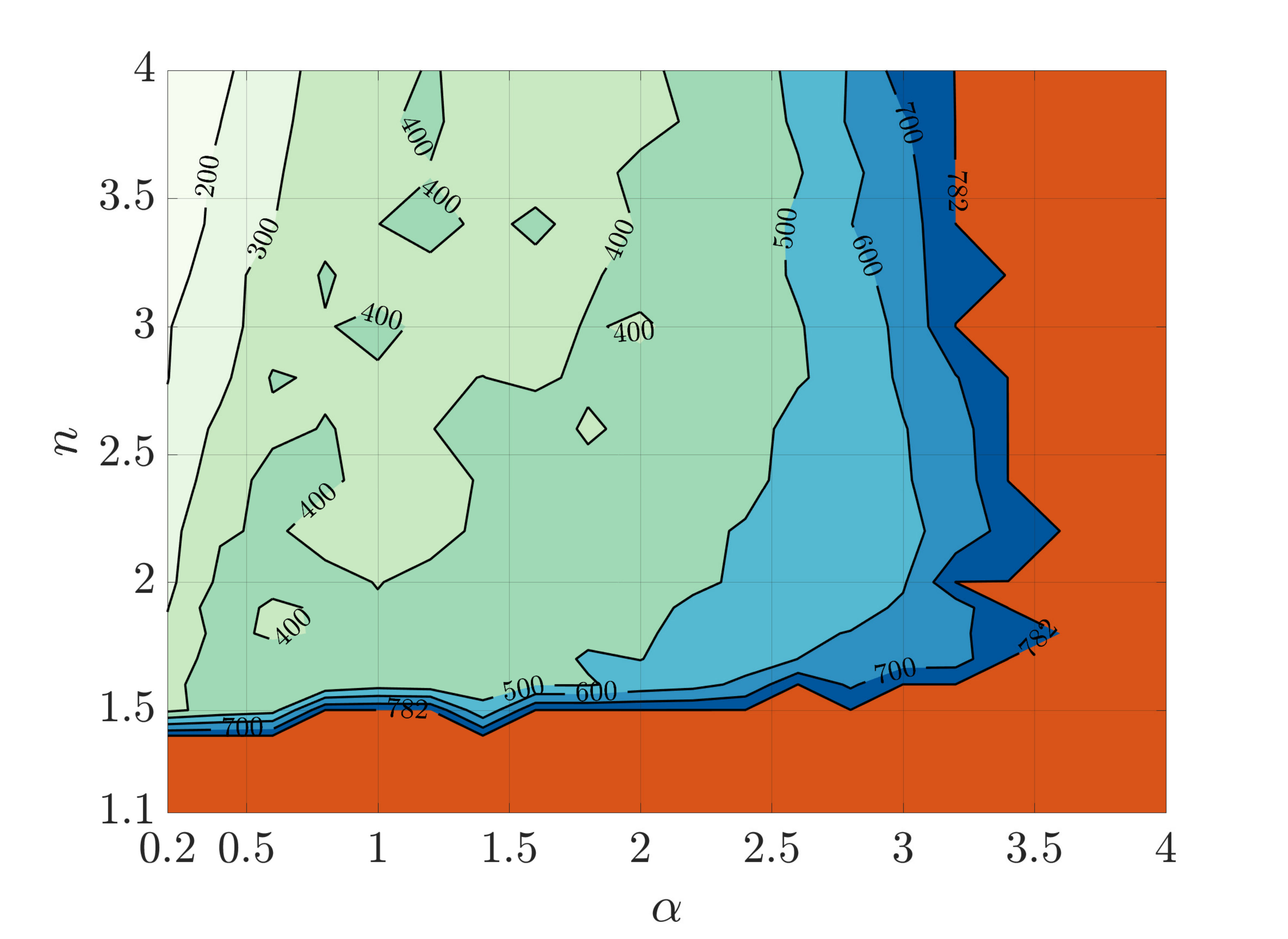}
\subcaption{$h=1/64,\Delta t = 1/128$}
\end{subfigure}
\caption{Contour plots of the average number of multigrid iterations needed to solve \eqref{travellingwave} using the modified Picard-CCMG solver for an anisotropic hydraulic conductivity field generated using $\Phi_2$.}\label{HM_Phi2}
\end{figure}

In Figures \ref{HM_Phi1}-\ref{HM_Phi2}, we show the average cost (average of 64 random realizations of the hydraulic conductivity) for solving the stochastic Richards' equation for four different combinations of spatial and temporal grid sizes and for two Mat\'ern parameter sets, $\Phi_1$ and $\Phi_2$, listed in Table \ref{Matern}. The region in red in the figures denotes the $(\alpha,n)$ values for which the modified Picard-CCMG solver failed to converge at least once out of the 64 samples. 

Based on numerical experiments, the performance of the modified Picard-CCMG solver for the stochastic Richards' equation can be summarized as follows:
\begin{itemize}
\item In general, the cost increases by decreasing $n$ and increasing $\alpha$. The cost of the solver rises steeply for $n<1.5$ and $\alpha>3.0$, and the cost increment with respect to the decrease in the value of $n$ is more pronounced compared to the increase in $\alpha$. 
\item While a spatio-temporal mesh refinement improves the robustness with respect to $\alpha$ and $n$, the improvement is less pronounced for $n$ and may require a very fine mesh as $n\rightarrow 1$.
\item For a given spatio-temporal mesh, the modified Picard-CCMG solver is less robust and more expensive for anisotropic hydraulic conductivity compared to the isotropic case. A similar $(\alpha,n)$-robustness can be achieved for anisotropic cases by using a sufficiently refined mesh.
\end{itemize}

\begin{figure}
\begin{subfigure}[b]{0.49\textwidth}
\includegraphics[trim={2.2cm 0cm 0cm 0cm},clip,scale = 0.2]{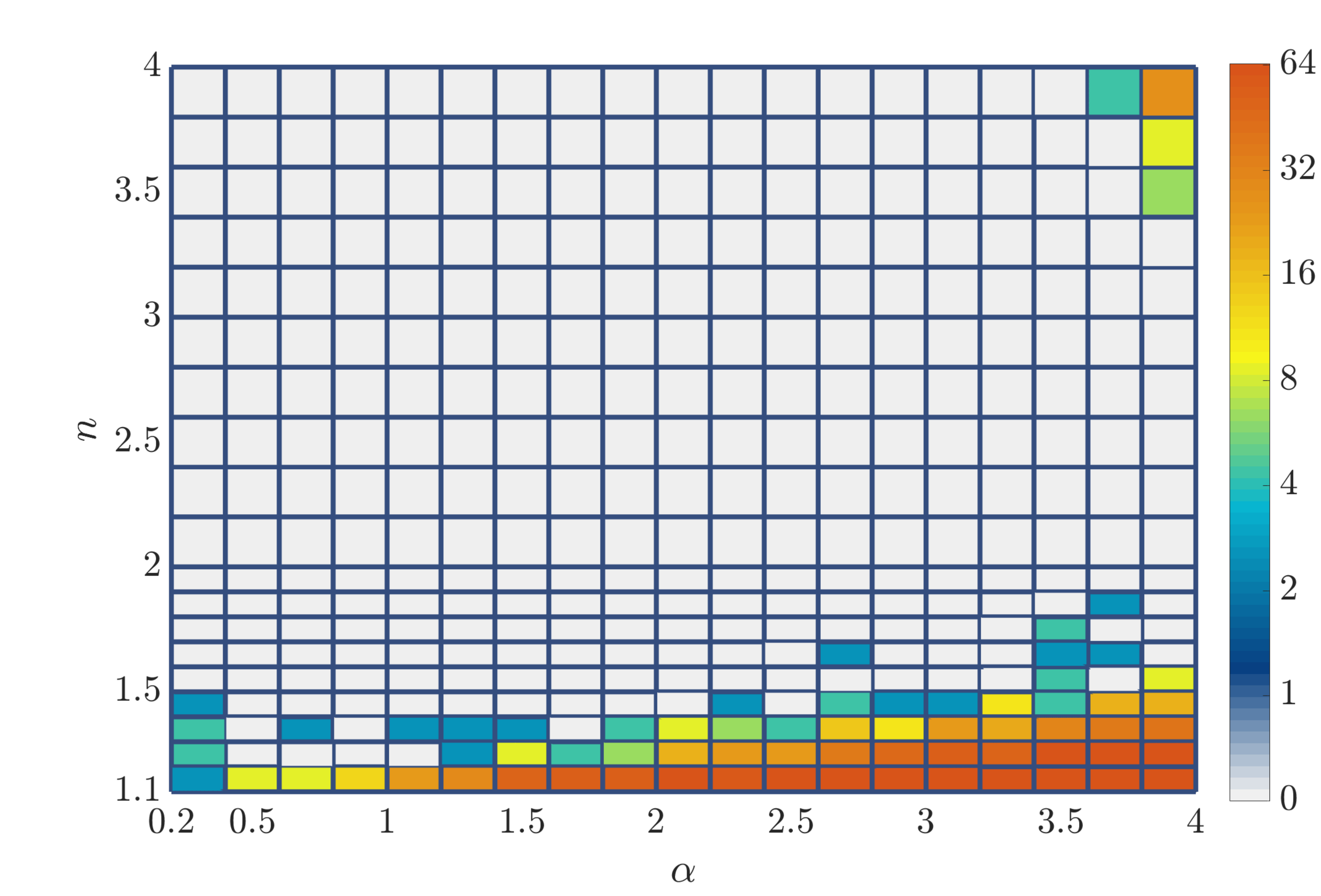}
\subcaption{$h=1/64,\Delta t = 1/128,\Phi_1$}
\end{subfigure}
\begin{subfigure}[b]{0.49\textwidth}
\includegraphics[trim={2.2cm 0cm 0cm 0cm},clip,scale = 0.2]{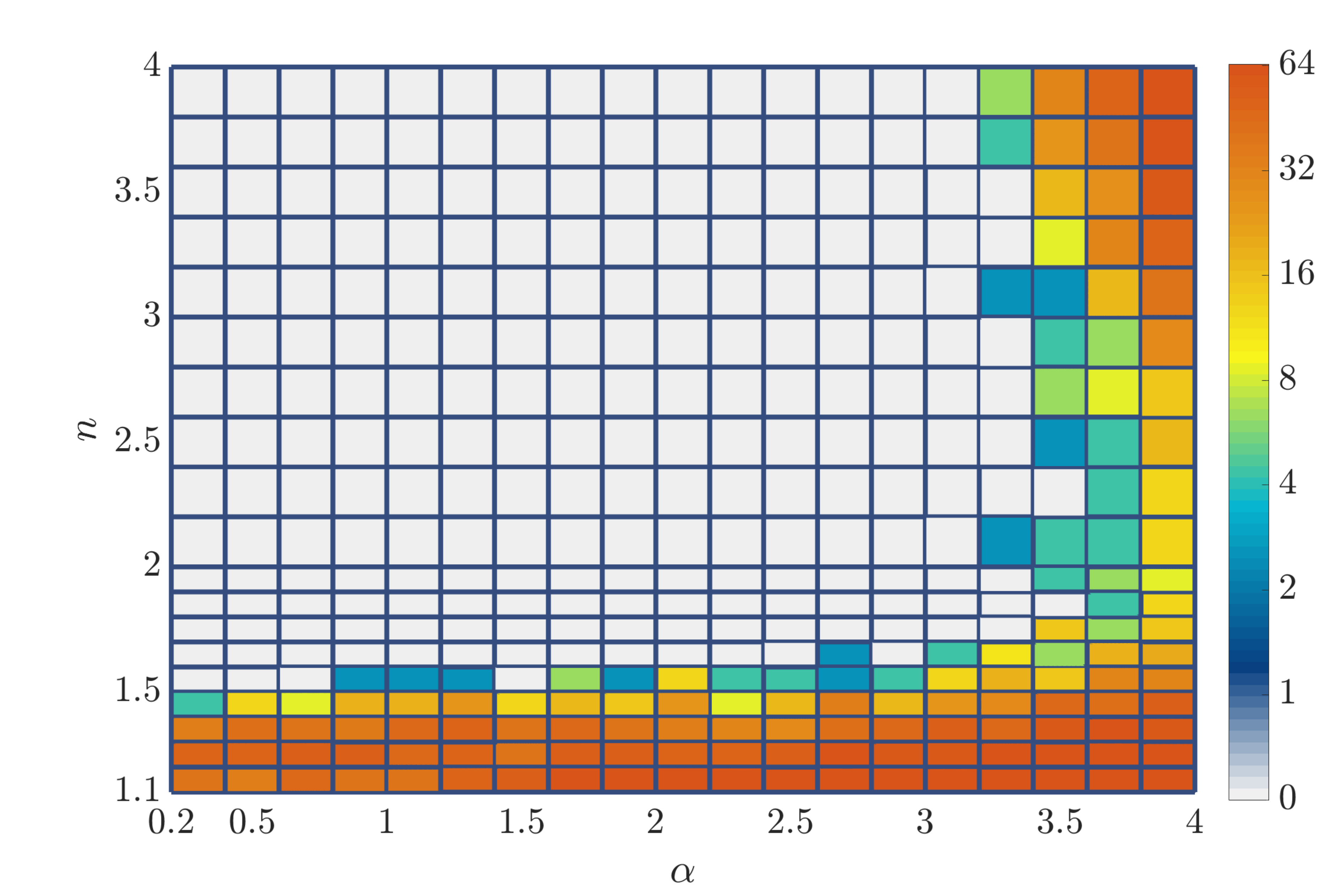}
\subcaption{$h=1/64,\Delta t = 1/128,\Phi_2$}
\end{subfigure}
\caption{Counting the number of samples (out of 64), for which the modified Picard-CCMG solver does not converge.}\label{failures}
\end{figure}

The standard deviation contours for the cost show a similar behavior as the average cost contour and we observe a large standard deviation for the cost when $n<1.5$ and $\alpha>3.0$. In Figure \ref{failures}, we present the number of samples (out of 64 samples), for which the solver did not converge for $\Phi_1$ and $\Phi_2$. For almost all samples convergence failed with $\alpha$ close to 4.0 and $n$ close to 1.1.

A few remarks are in order. We point out that the $(\alpha,n)$-cost map may vary depending on the type of boundary and initial conditions as well as on $T_{final}$. For instance, in the above experiments, an initially wet profile for the porous media was considered. We expect the performance of the modified Picard-CCMG solver to vary for problems in which infiltration takes place into an initially dry media and the convergence rates may depend on the values of $\theta_r$ and $\theta_s$ (see e.g. \cite{ZHA2017}). Furthermore, the robustness of the solver will also depend on the properties of the hydraulic conductivity field such as on the degree of heterogeneity and anisotropy.  These topics will be actively explored in the future work.

\section{Multilevel Monte Carlo with parametric continuation}\label{nonMLMC}
We have observed in the preceding section that the total number of multigrid iterations increases rapidly with a decrease in the value of parameter $n$ and an increase in $\alpha$. We also noticed that the solver is less robust on a coarse spatio-temporal mesh. Therefore, when using the original MLMC estimator for a ``difficult'' $(\alpha,n)$ pair, a relatively fine spatio-temporal mesh will be required (and employed), even on the coarsest level of the MLMC hierarchy, resulting in an expensive estimator. To deal with this drawback, we propose an MLMC estimator based on the \emph{parametric continuation} technique. In this approach, we solve the original problem only on the finest level of the MLMC hierarchy and {\em simplify} the parameter settings dictating the nonlinearities as we work on coarser levels. This allows us to include a comparatively coarser spatio-temporal mesh compared to the original MLMC estimator as simpler problems are solved on coarser levels. 

This idea is motivated by continuation based multigrid solvers for nonlinear boundary value problems \cite{Bra77,brandt2003multigrid,Mittelmann1}. In the context of multigrid solvers, continuation is commonly applied in the FMG-FAS (Full MultiGrid- Full Approximation Scheme) algorithm. In these algorithms,  the continuation process is integrated in the FMG hierarchy, where the coarse grid solves the simplest problem and is used as a good first approximation for the next grid with a slightly more complicated problem. This process is repeated until the finest grid is reached where the target problem is solved.  Although the continuation strategy works well for a large class of nonlinear problems, there is no guarantee that the simpler problem is close enough to the next difficult problem. One can use \emph{bifurcation diagrams} to understand the solution dependence on nonlinearity dictating parameters. These diagrams can also reveal multiple branches and bifurcation points, where the solution differs greatly even if there is a slight perturbation in the parameter value. In such cases, an \emph{arclength} procedure \cite{arclength} can be applied to determine the appropriate perturbation size.

\subsection{MLMC estimator}

To explain the MLMC estimator, we consider the  pressure head field at some final time $T_{final}$ as the QoI. Further, we define a spatio-temporal hierarchy of grid levels $\{\mathcal{D}_{\ell},\mathcal{T}_{\ell}\}^L_{\ell=0}$ using
\begin{equation}\label{hmax_ch6}
h_\ell  = \Delta t_{\ell} = \mathcal{O}(s^{-\ell}h_0),
\end{equation}
where  $h_0$ is the cell-width on the coarsest mesh $\mathcal{D}_0$ and $s>0$ represents a grid refinement factor. We further define a hierarchy of parameter sets, $\{\boldsymbol{\Theta}_{\ell}\}_{\ell=0}^L$, where $\boldsymbol{\Theta}_L$ is the parameter set corresponding to the target (strongly nonlinear) problem to be solved. For instance, we can define a parametric hierarchy using the set of van Genuchten parameters, e.g. $\boldsymbol{\Theta}_{\ell} = \{\alpha_\ell^{(bl)},n_\ell^{(bl)}\}$. The approximation of the pressure head on the level $\ell$ at  $T_{final}$ is denoted by $p_{h_\ell,\boldsymbol{\Theta}_\ell}$. Using the linearity of the expectation operator, one can define the expected value of the pressure head on the finest level, $L$, with the original parameter set, $\boldsymbol{\Theta}_{L}$, by the following telescopic sum:
\begin{equation}\label{telescope_ch6}
\mathbb{E}[p_{h_L,\boldsymbol{\Theta}_L}] = \mathbb{E}[p_{h_0,\boldsymbol{\Theta}_0}] + \sum_{\ell =1}^L \mathbb{E}[p_{h_\ell,\boldsymbol{\Theta}_\ell} - p_{h_{{\ell-1}},\boldsymbol{\Theta}_{\ell-1}}].
\end{equation}
Note that for $\boldsymbol{\Theta}_0=\boldsymbol{\Theta}_1=...=\boldsymbol{\Theta}_L$, we have the standard MLMC estimator which solves the same problem on all levels. In terms of computational effort, it is cheaper to approximate $\mathbb{E}[p_{h_0,\boldsymbol{\Theta}_0}]$ by a standard Monte Carlo estimator, as the samples are computed on a coarse spatio-temporal mesh based on an ``easy nonlinear parameter set $\mathbf{\Theta}_0$''. Furthermore, the correction term, $\mathbb{E}[p_{h_\ell,\boldsymbol{\Theta}_\ell} - p_{h_{\ell-1},\boldsymbol{\Theta}_{\ell-1}}]$, can be accurately determined using only a few samples as the level-dependent variance, $\mathcal{V}_\ell := \mathbb{V}[p_{h_\ell,\boldsymbol{\Theta}_\ell} - p_{h_{\ell-1},\boldsymbol{\Theta}_{\ell-1}}]$, is typically small, since the random variables  $p_{h_\ell,\boldsymbol{\Theta}_\ell}$ and $ p_{h_{\ell-1},\boldsymbol{\Theta}_{\ell-1}}$ are positively correlated. Note that the correlation will depend on the grid parameters $h_\ell$ and $h_{\ell-1}$ as well as on the difference between the nonlinear parameters $\boldsymbol{\Theta}_\ell$ and $\boldsymbol{\Theta}_{\ell-1}$. We will elaborate on this later on.

%More precisely, the value of the level-dependent variance depends on the correlation between $p_{h_\ell,\boldsymbol{\Theta}_\ell}$ and $ p_{h_{\ell-1},\boldsymbol{\Theta}_{\ell-1}}$, i.e.
%\begin{equation}
%\mathcal{V}_\ell = \mathbb{V}[p_{h_\ell,\boldsymbol{\Theta}_\ell}]+\mathbb{V}[p_{h_{\ell-1},\boldsymbol{\Theta}_{\ell-1}}]-2\text{Cov}(p_{h_\ell,\boldsymbol{\Theta}_\ell},p_{h_{\ell-1},\boldsymbol{\Theta}_{\ell-1}}).
%\end{equation}

%In general, the correlation between $p_{h_\ell,\boldsymbol{\Theta}_\ell}$ and $ p_{h_{\ell-1},\boldsymbol{\Theta}_{\ell-1}}$ deteriorates as the difference between $h_\ell$ and $h_{\ell-1}$ increase. Further, the correlation also depends on the parameters $\boldsymbol{\Theta}_\ell$ and $\boldsymbol{\Theta}_{\ell-1}$ associated with the nonlinearity.  As, For simplicity, we assume

Each of the expectations in the MLMC estimator \eqref{telescope_ch6} can be independently computed using the standard MC simulation. We define a multilevel estimator, $\mathpzc{E}^{ML}_{L}[p_{h_L,\boldsymbol{\Theta}_L}]$, constructed using a sum of $L+1$ MC estimators: 
\begin{align}\label{MLMCestimator_ch6}
\mathbb{E}[p_{h_L,\boldsymbol{\Theta}_L}] \approx \mathpzc{E}^{ML}_L[p_{h_L,\boldsymbol{\Theta}_L}] : =&\sum^L_{\ell=0} \mathpzc{E}^{MC}_{N_\ell}[p_{h_\ell,\boldsymbol{\Theta}_\ell} - p_{h_{\ell-1},\boldsymbol{\Theta}_{\ell-1}}] ,
\end{align}
where $\mathpzc{E}^{MC}_{N_\ell}[p_{h_\ell,\boldsymbol{\Theta}_\ell} - p_{h_{\ell-1},\boldsymbol{\Theta}_{\ell-1}}]$ is the standard MC estimator obtained by averaging $N_\ell$ independent, identically distributed (i.i.d.) samples as
\begin{equation}\label{SMC_ch6}
\mathpzc{E}^{MC}_{N_\ell}[p_{h_\ell,\boldsymbol{\Theta}_\ell} - p_{h_{\ell-1},\boldsymbol{\Theta}_{\ell-1}}]:=\left(\frac{1}{N_\ell}\sum^{N_\ell}_{i=1} (p_{h_\ell,\boldsymbol{\Theta}_\ell}(\omega_i) - p_{h_{\ell-1},\boldsymbol{\Theta}_{\ell-1}}(\omega_i))\right).
\end{equation}
with $\omega_i$ denoting an event in the stochastic domain $\Omega$ and $p_{h_{-1},\boldsymbol{\Theta}_{-1}}=0$. It is expected that the number of MLMC samples $N_\ell\in\mathbb{N}$ forms a decreasing sequence for increasing $\ell$. In order to keep the variance of the correction terms small, the MC samples, $p_{h_\ell,\boldsymbol{\Theta}_\ell}(\omega_i) - p_{h_{\ell-1},\boldsymbol{\Theta}_{\ell-1}}(\omega_i)$, should be based on the same random input $\omega_i$ for the simulation on two consecutive levels $\ell$ and $\ell-1$. 

Similarly, a multilevel estimator for the variance of the pressure head, $\mathbb{V}[p_{h_L,\boldsymbol{\Theta}_L}]$, can be defined as
\begin{equation}\label{ML_var_ch6}
\mathbb{V}[p_{h_L,\boldsymbol{\Theta}_L}] \approx \mathcal{V}_L^{ML}[p_{h_L,\boldsymbol{\Theta}_L}] := \sum_{\ell=0}^L \mathcal{V}^{MC}_{N_\ell}[p_{h_\ell,\boldsymbol{\Theta}_\ell}]-\mathcal{V}^{MC}_{N_{\ell}}[p_{h_{\ell-1},\boldsymbol{\Theta}_{\ell-1}}],
\end{equation}
where the variance $\mathcal{V}^{MC}_{N_\ell}[p_{h_\ell,\boldsymbol{\Theta}_\ell}]$ is computed as 
\begin{equation}\label{MC_var_ch6}
\mathcal{V}^{MC}_{N_\ell}[p_{h_\ell,\boldsymbol{\Theta}_\ell}] \approx \frac{1}{N_\ell-1}\sum_{i=1}^{N_\ell} \left(p_{h_\ell,\boldsymbol{\Theta}_\ell}(\omega_i) -\mathpzc{E}_{N_{\ell}}^{MC}[p_{h_\ell,\boldsymbol{\Theta}_\ell}]\right)^2.
\end{equation} 
Again, the computational savings for the variance estimator \eqref{ML_var_ch6} are obtained by computing individual variances $\mathcal{V}^{MC}_{N_\ell}[p_{h_\ell,\boldsymbol{\Theta}_\ell}]$ and $\mathcal{V}^{MC}_{N_{\ell}}[p_{h_{\ell-1},\boldsymbol{\Theta}_{\ell-1}}]$ using the same random inputs $\{\omega_i\}_{i=1}^{N_\ell}$. The above variance estimator can be seen as an extension of the standard multilevel variance estimator proposed in \cite{bierig}. 

For the multilevel estimators, an appropriate spatial interpolation procedure is required to combine expectations from all levels. Typically, the polynomial order of the interpolation scheme should be equal to or higher than the order of the discretization to avoid any additional dominant source of error. In some more detail, when using the estimator \eqref{MLMCestimator_ch6} to compute $\mathpzc{E}^{ML}_L[p_{h_L,\boldsymbol{\Theta}_L}]$, we begin by computing $\mathpzc{E}_{N_0}^{MC}[p_{h_0,\boldsymbol{\Theta}_0}]$ on the coarsest grid $\mathcal{D}_0$. This quantity is then interpolated to the next finer grid $\mathcal{D}_1$ and is added to the correction term $\mathpzc{E}^{MC}_{N_1}[p_{h_1,\boldsymbol{\Theta}_1} - p_{h_0,\boldsymbol{\Theta}_0}]$ resulting in a two-level estimate $\mathpzc{E}_{1}^{ML}[p_{h_1,\boldsymbol{\Theta}_1}]$. This is again interpolated to the next grid level $\mathcal{D}_2$ and added to the next correction term $\mathpzc{E}^{MC}_{N_2}[p_{h_2,\boldsymbol{\Theta}_2} - p_{h_1,\boldsymbol{\Theta}_1}]$. This process is repeated until the final level is reached. 

\subsection{Accuracy of MLMC estimator}
Throughout this paper, we use the $L^2-$ based norm for the error analysis of the multilevel Monte Carlo estimator. We assume that the pressure considered belongs to the functional space $\Lomd$ corresponding to the space of square-integrable measurable functions $p:\Omega\rightarrow L^2(\mathcal{D})$ for a previously defined probability space $(\Omega,\mathbb{F},\mathbb{P})$. These spaces are equipped with the norm
\begin{equation}\label{eq:RMSEnorm}
\lnorm p(\mathbf{x},T,\omega)\rnorm_{\Lomd} := \mathbb{E}\left[\lnorm p(\mathbf{x},T,\omega)\rnorm^2_{L^2(\mathcal{D})} \right]^{\tfrac12} =\left(\int_{\Omega}\lnorm p(\mathbf{x},T,\omega)\rnorm^2_{L^2(\mathcal{D})} \text{d}\mathbb{P} \right)^{\tfrac12}.
\end{equation}
The mean-square error (MSE) in  $\mathpzc{E}^{ML}_{L}[p_{h_L,\boldsymbol{\Theta}_L}]$ can then be expressed as the sum of the discretization and the sampling errors as

\begin{align}\label{MSEMLMC_ch6}
\lnorm\mathbb{E}[p_{\boldsymbol{\Theta}_L}]-  \mathpzc{E}^{ML}_{L}[p_{h_L,\boldsymbol{\Theta}_L}]\rnorm_{\Lomd}^2
\leq & \lnorm\mathbb{E}[p_{\boldsymbol{\Theta}_L}]  -   \mathbb{E}[p_{h_L,\boldsymbol{\Theta}_L}]\rnorm_{\Ld}^2 +\lnorm\mathbb{E}[p_{h_L,\boldsymbol{\Theta}_L}]-  \mathpzc{E}^{ML}_L[p_{h_L,\boldsymbol{\Theta}_L}] \rnorm_{\Lomd}^2.
%=& (C_1h_L^{\alpha})^2 + \sum^L_{\ell=0}\frac{\mathcal{V}_\ell}{N_\ell},
\end{align} 
Both errors in the MLMC estimator can be dealt with separately. The discretization error can be quantified as:
\begin{equation}\label{numErr}
\lnorm\mathbb{E}[p_{\boldsymbol{\Theta}_L}]  -   \mathbb{E}[p_{h_L,\boldsymbol{\Theta}_L}]\rnorm_{\Ld}\leq c_0h_L^{a},\qquad a>0,
\end{equation}
where $c_0$ is a constant independent of $h_L$ but depending on the parameter set ${\boldsymbol{\Theta}_L}$. The rate $a$ typically depends on the regularity of the PDE and the accuracy of the discretization.  
The next task is to bound the sampling errors. As the MLMC estimator $\mathpzc{E}^{ML}_L[p_{h_L,\boldsymbol{\Theta}_L}]$ is composed of $L+1$ independent MC estimators, the sampling error in the MLMC estimator is just the sum of sampling errors from the individual MC estimators. Therefore,
\begin{equation}\label{eq:MLMCVar}
\lnorm\mathbb{E}[p_{h_L,\boldsymbol{\Theta}_L}]-  \mathpzc{E}^{ML}_L[p_{h_L,\boldsymbol{\Theta}_L}] \rnorm_{\Lomd}^2=\sum^L_{\ell=0}\frac{\lnorm\mathcal{V}_\ell\rnorm_{\Ld}}{N_\ell},
\end{equation}
see \cite{mishra2012sparse,MULLER2013685} for a proof. Obtaining a bound on the level-variance $\lnorm\mathcal{V}_\ell\rnorm_{\Ld}$ is more involved due to its dependence on the grid size $h_\ell$ as well as on the nonlinearity parameter set $\boldsymbol{\Theta}_{\ell}$. We numerically estimate it by
\begin{eqnarray}\label{level_var_ch6}
\lnorm\mathcal{V}_\ell\rnorm_{\Ld} &=& \lnorm\mathbb{V}[p_{h_\ell,\boldsymbol{\Theta}_\ell} - p_{h_{\ell-1},\boldsymbol{\Theta}_{\ell-1}}]\rnorm_{\Ld}\nonumber \\
&\approx& \frac{1}{N_\ell-1}\sum_{i=1}^{N_\ell}\int_\mathcal{D}\Bigg( \mathpzc{E}^{MC}_{N_\ell}[p_{h_\ell,\boldsymbol{\Theta}_\ell} - p_{h_{\ell-1},\boldsymbol{\Theta}_{\ell-1}}] -  (p_{h_\ell,\boldsymbol{\Theta}_\ell}(\omega_i) - p_{h_{\ell-1},\boldsymbol{\Theta}_{\ell-1}}(\omega_i)) \Bigg)^2.
\end{eqnarray}
To achieve a tolerance of $\varepsilon$, one needs to ensure that
\begin{equation}\label{errMLMC}
\lnorm\mathbb{E}[p_{\boldsymbol{\Theta}_L}]-  \mathpzc{E}^{ML}_{L}[p_{h_L,\boldsymbol{\Theta}_L}]\rnorm_{\Lomd}^2 \leq (c_0h_L^{a})^2 + \sum^L_{\ell=0}\frac{\lnorm\mathcal{V}_\ell\rnorm_{\Ld}}{N_\ell}< \varepsilon^2.
\end{equation}
The total cost of the MLMC estimator can be expressed as $\mathcal{W}^{ML}_{L} =\sum^L_{\ell=0}N_\ell \mathcal{W}_{\ell}$, where $\mathcal{W}_{\ell} = \mathcal{O}(h_{\ell}^{-\gamma})$ corresponds to the cost of computing one sample on level $\ell$. For time-dependent problems, the rate $\gamma \geq d+1$, with $d$ the number of spatial dimensions.  As proposed in \cite{MLMC1,MLMC2}, the number of samples at different levels is typically derived by minimizing the total cost such that the sampling error of the MLMC estimator reduces below $\varepsilon^2$, i.e.,
\begin{equation}\label{MLMCopt_ch6}
 \text{min}\left(\sum_{\ell=0}^L N_\ell \mathcal{W}_\ell\right)\quad\text{s.t}\quad\sum^L_{\ell=0}\frac{\lnorm\mathcal{V}_\ell\rnorm_{\Ld}}{N_\ell} = \varepsilon^2.
\end{equation}
Using the standard Lagrange multiplier approach \cite{MLMC1}, gives us 
\begin{equation}\label{MLMCSAMP_ch6}
N_\ell =\varepsilon^{-2}\left(\sum_{\ell=0}^L\sqrt{\lnorm\mathcal{V}_\ell\rnorm_{\Ld}\mathcal{W}_\ell}\right)\sqrt{\frac{\lnorm\mathcal{V}_\ell\rnorm_{\Ld}}{\mathcal{W}_\ell}},
\end{equation}
and hence the total cost to obtain a tolerance of $\varepsilon$ is given by
\begin{equation}
\mathcal{W}^{ML}_L(\varepsilon)=\sum_{\ell=0}^LN_\ell\mathcal{W}_\ell=\varepsilon^{-2}\left(\sum_{\ell=0}^L\sqrt{\lnorm\mathcal{V}_\ell\rnorm_{\Ld}\mathcal{W}_\ell}\right)^2.
\end{equation}
In the above formula, the product ${\lnorm\mathcal{V}_\ell\rnorm_{\Ld}\mathcal{W}_\ell}$ determines the cost contribution from any level $\ell$. For instance, if the product decays with increasing $\ell$, the dominant cost comes from the coarsest level whereas if the product grows with $\ell$, the dominant contribution comes from the finest level. 

\begin{remark}
The optimal number of samples given in \eqref{MLMCSAMP_ch6} is based on a pre-defined hierarchy of parameters $\{\boldsymbol{\Theta}_{\ell}\}_{\ell=0}^L$. A more general approach is to find $N_\ell$ along with the parameter set $\{\boldsymbol{\Theta}_{\ell}\}_{\ell=0}^L$ for which the total cost of the MLMC estimator is minimum. Solving such optimization problem analytically is non-trivial. Furthermore, numerically obtaining the best values for $\boldsymbol{\Theta}_{\ell}$ can also be highly expensive. In the numerical experiments section, we will discuss some heuristics that can be applied to find $\boldsymbol{\Theta}_{\ell}$.
\end{remark}
%One of the principle issues in designing the continuation based MLMC estimator is deciding the step size for the parameter $\boldsymbol{\Theta}$.
\subsubsection{MLMC algorithm with parametric continuation}

To compute the estimator $\mathpzc{E}^{ML}_L[p_{h_L,\boldsymbol{\Theta}_L}]$, the standard MLMC algorithm from \cite{MLMC1,MLMC2} cannot be directly employed as it requires solving the same problem on all grid levels.  Here, we describe a {\em modified version} of the standard MLMC technique to compute $\mathpzc{E}^{ML}_L[p_{h_L,\boldsymbol{\Theta}_L}]$. This algorithm assumes that the total number of levels in the MLMC hierarchy and the values of the nonlinearity parameters $\boldsymbol{\Theta}_\ell$ for all levels are known in advance. The algorithm can be described by the following steps:
\begin{algorithm}[H]
\caption{$PC\_MLMC$ algorithm}\label{PC_MLMC}
\begin{algorithmic}[1]
\State Fix the tolerance $\varepsilon$, $\mathcal{D}_\ell,\mathcal{T}_\ell,\boldsymbol{\Theta}_\ell$ and warm-up samples $N_\ell^*$ for $\ell = 0,1,2,...,L$.
\State Compute quantities $\mathpzc{E}^{MC}_{N_\ell}[p_{h_\ell,\boldsymbol{\Theta}_\ell} - p_{h_{\ell-1},\boldsymbol{\Theta}_{\ell-1}}]$ and $\lnorm\mathcal{V}_\ell\rnorm_{\Ld}$ using samples $N_\ell  = N^*_\ell$ for all levels.
\State Update $N_\ell$ using the formula \eqref{MLMCSAMP_ch6} for all levels.
\State Compute additional samples and update $\mathpzc{E}^{MC}_{N_\ell}[p_{h_\ell,\boldsymbol{\Theta}_\ell} - p_{h_{\ell-1},\boldsymbol{\Theta}_{\ell-1}}]$ and $\lnorm\mathcal{V}_\ell\rnorm_{\Ld}$ for all levels.
\State Perform steps 3-4 until no additional samples are needed on any level.
\end{algorithmic}
\end{algorithm}
In the above algorithm, the value of $N^*_{\ell}$ should not be set too high, especially not on the finest level, in order to avoid oversampling. Further, the cost per sample $\mathcal{W}_\ell$ can also be estimated ``on-the-fly'' by averaging the CPU times from the computation of warm-up samples.
\section{Numerical experiments}\label{numexp2}
We evaluate the performance of the new MLMC estimator and study the improvements with respect to the standard MLMC estimator.  For all the experiments, we use the infiltration problem with conditions given  in \eqref{travellingwave}, however, with $T_{final}=0.2$ [h] (in hours) and  the two Mat\'ern covariance parameters from Table \ref{Matern}. We employ a geometric hierarchy of spatio-temporal grids with refinement factor $s=2$ in \eqref{hmax_ch6} and we use $h_\ell = \Delta t_\ell$. For all experiments, the following baseline values are prescribed, $K_s^{(bl)} = 0.2$ [m/h] \text{(in metres/hour)}, $\theta_s^{(bl)} = 0.5$ and $\theta_r^{(bl)} = 0.05$; different baseline values for $\alpha^{(bl)}$ and $n^{(bl)}$ are studied. The uncertainty in the soil parameters is defined according to the values presented in Table \ref{soilparams}. The sampling and upscaling procedure for a Gaussian random field is described in \ref{appendix_A1}. The sampling of random fields with uniform marginal is described in Section \ref{stochRE}.
\begin{table}[H]
\begin{center}
\begin{tabular}{|c|c|}
\hline
Quantity & Uncertainty\hspace{1.2cm} \\
\hline
$Z(\mathbf{x},\omega)$ & $\mathcal{N}(0,C_\Phi)\hspace{1.6cm}$\\
$\varepsilon_{\theta_s}(\mathbf{x},\omega) $ & $\mathcal{U}(-0.05,0.05,C_{\Phi})\hspace{0.3cm}$   \\
$\varepsilon_{\theta_r}(\mathbf{x},\omega)$   & $\mathcal{U}(-0.005,0.005,C_{\Phi})$ \\
$\varepsilon_{\alpha}(\mathbf{x},\omega)$     & $\mathcal{U}(-0.2,0.2,C_{\Phi})\hspace{0.7cm}$    \\
$\varepsilon_{n}(\mathbf{x},\omega)$             & $\mathcal{U}(-0.05,0.05,C_{\Phi})\hspace{0.3cm}$ \\
\hline
\end{tabular}
\end{center}
\caption{Description of uncertainty for different soil parameters.}\label{soilparams}
\end{table}
Note that the above stochastic model is extremely high-dimensional as it comprises five independent random fields. For each random field the degree of freedom is equal to the number of grid points in the sampling mesh. The dimensionality can be reduced using the KL-expansion method \eqref{KLexp_ch6}, however as we use random fields with small correlation lengths, we will still need to use a very large number of KL-modes for an accurate representation of these random fields. 
\subsection{Convergence of discretization bias}
\begin{figure}[H]
\begin{subfigure}[b]{0.49\textwidth}
\includegraphics[trim={0cm 0cm 0cm 0cm},clip,scale = 0.285]{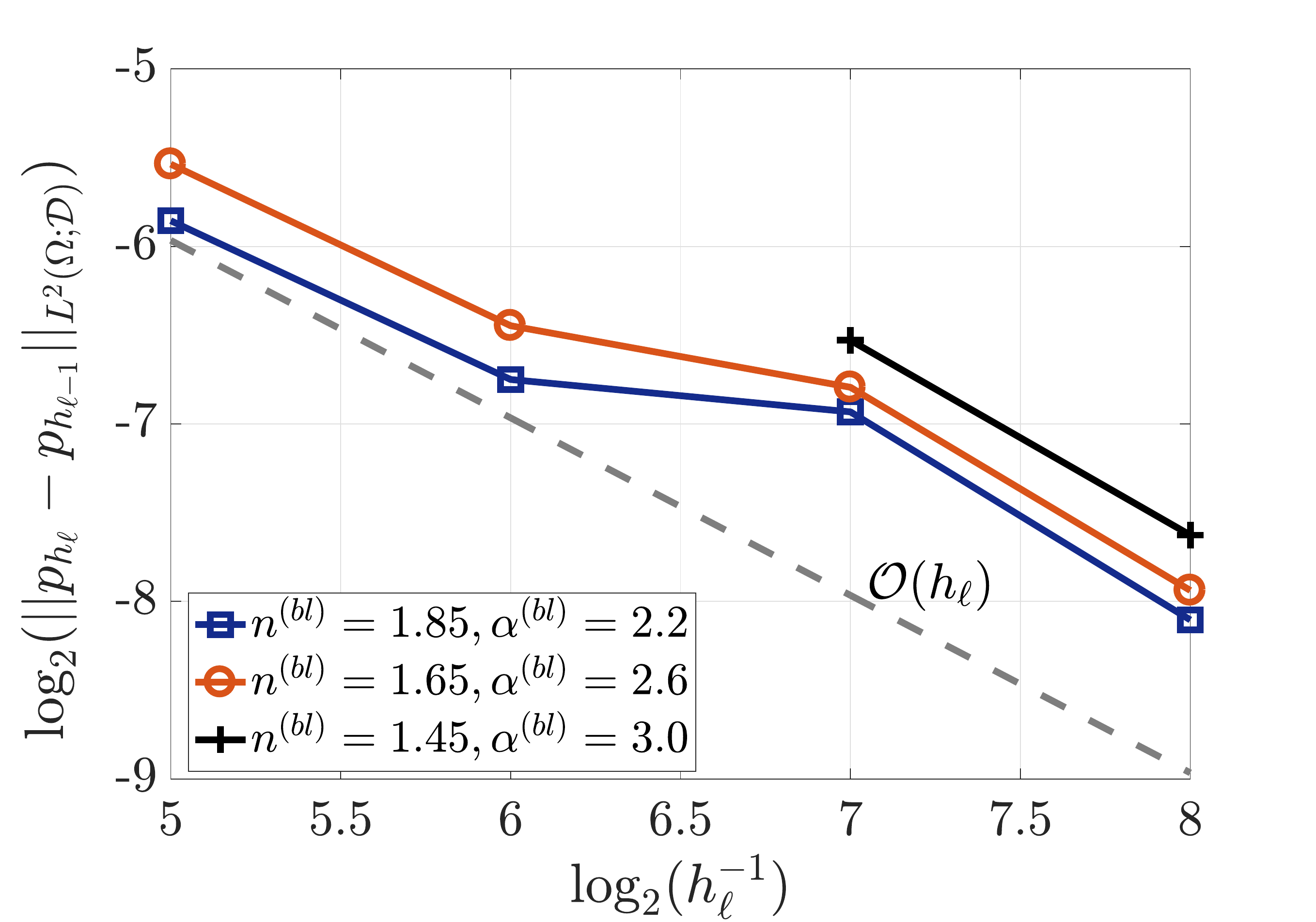}
\subcaption{$\Phi_1$}
\end{subfigure}
\begin{subfigure}[b]{0.49\textwidth}
\includegraphics[trim={0cm 0cm 0cm 0cm},clip,scale = 0.28]{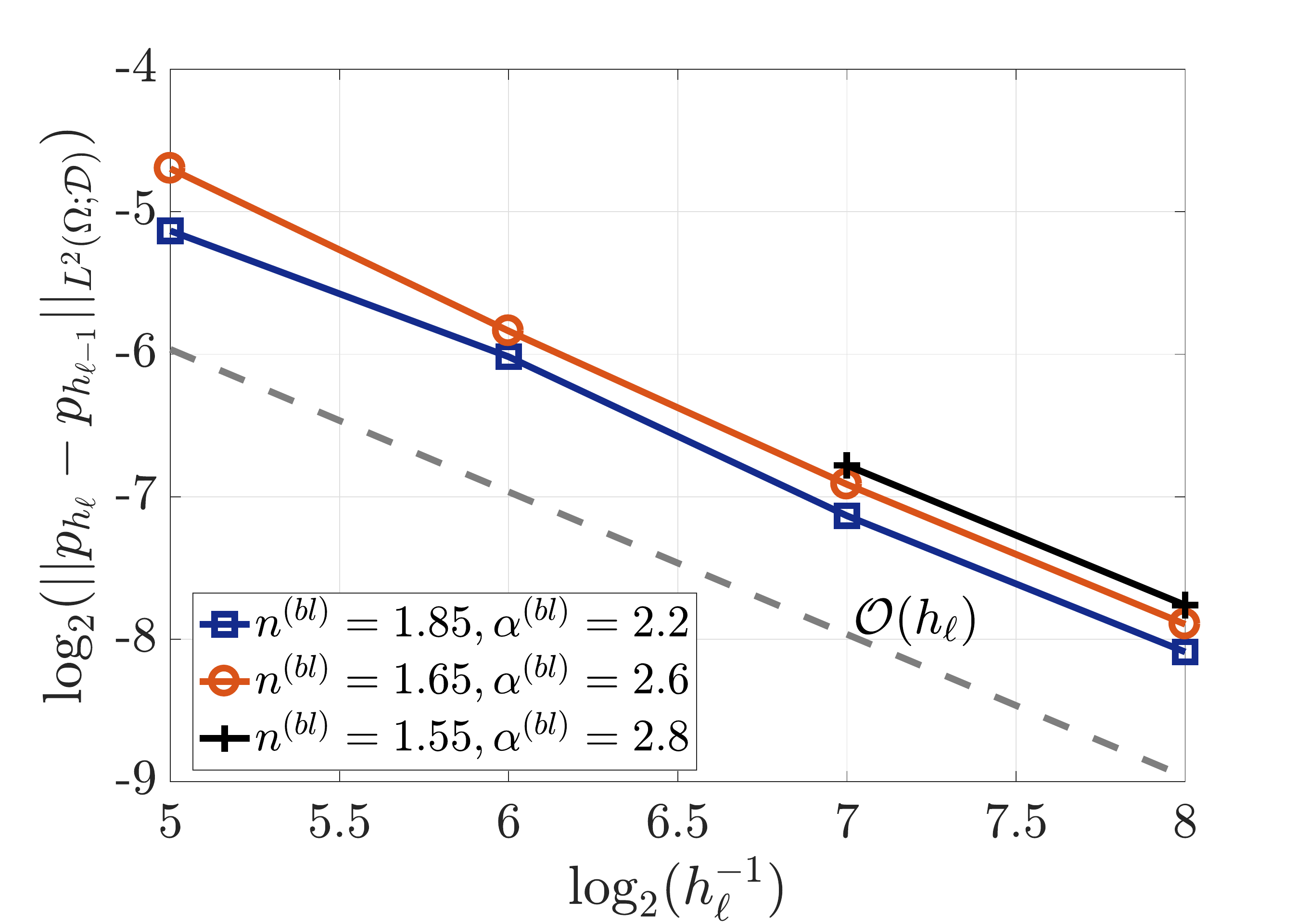}
\subcaption{$\Phi_2$}
\end{subfigure}
\caption{Convergence of discretization error with respect to mesh refinement for different baseline values of $\alpha^{(bl)}$ and $n^{(bl)}$.}\label{FV_err_ch6}
\end{figure}

We begin by analyzing the reduction  of the discretization error $||p_{h_\ell} - p_{h_{\ell-1}}||_{L^2(\Omega;\mathcal{D})}$ with respect to mesh refinement for different baseline values of $\alpha^{(bl)}$ and $n^{(bl)}$. The relative error is  used to bound the exact discretization bias as 
\begin{equation}\label{discErr}
||p - p_{h_{\ell}}||_{L^2(\Omega;\mathcal{D})}\leq \frac{||p_{h_\ell} - p_{h_{\ell-1}}||_{L^2(\Omega;\mathcal{D})}}{s^a -1},
\end{equation}
where $a$ is the convergence rate defined in \eqref{numErr}. The relative errors for, $\Phi_1$ and $\Phi_2$ are presented in the left and right pictures in Figure \ref{FV_err_ch6}, respectively. For both cases a convergence rate close to first-order is observed, i.e. $a\approx1$.  The convergence rate typically depends on the order of the spatio-temporal discretization scheme as well as on the smoothness parameter $\nu_c$ in the covariance function. In fact, the dominant error comes from the first-order accurate backward Euler time discretization. The magnitude of the error grows with increasing $\alpha^{(bl)}$ and reduces with increasing $n^{(bl)}$ values. Note that for the most difficult cases, $n^{(bl)} = 1.45, \alpha^{(bl)} =3.0$ for $\Phi_1$ and $n^{(bl)} = 1.55, \alpha^{(bl)} =2.8$ for $\Phi_2$, the convergent solutions are obtained from $h_\ell =1/64$ onwards.
\subsection{MLMC simulation}
Here, we describe the algorithm to compute the multilevel estimator $\mathpzc{E}^{ML}_L[p_{h_L,\boldsymbol{\Theta}_L}]$. We perform the MLMC simulations for two test cases based on $\Phi_1$ and $\Phi_2$, respectively. The original problem for $\Phi_1$ uses $\mathbf{\Theta}_L = \{\alpha_L^{(bl)},n_L^{(bl)}\} = \{3.0,1.45\}$ and for $\Phi_2$ the original problem is based on $\mathbf{\Theta}_L = \{2.8,1.55\}$. 

\begin{figure}
\begin{subfigure}[b]{0.32\textwidth}
\includegraphics[trim={0cm 0cm 0cm 0cm},clip,scale = 0.205]{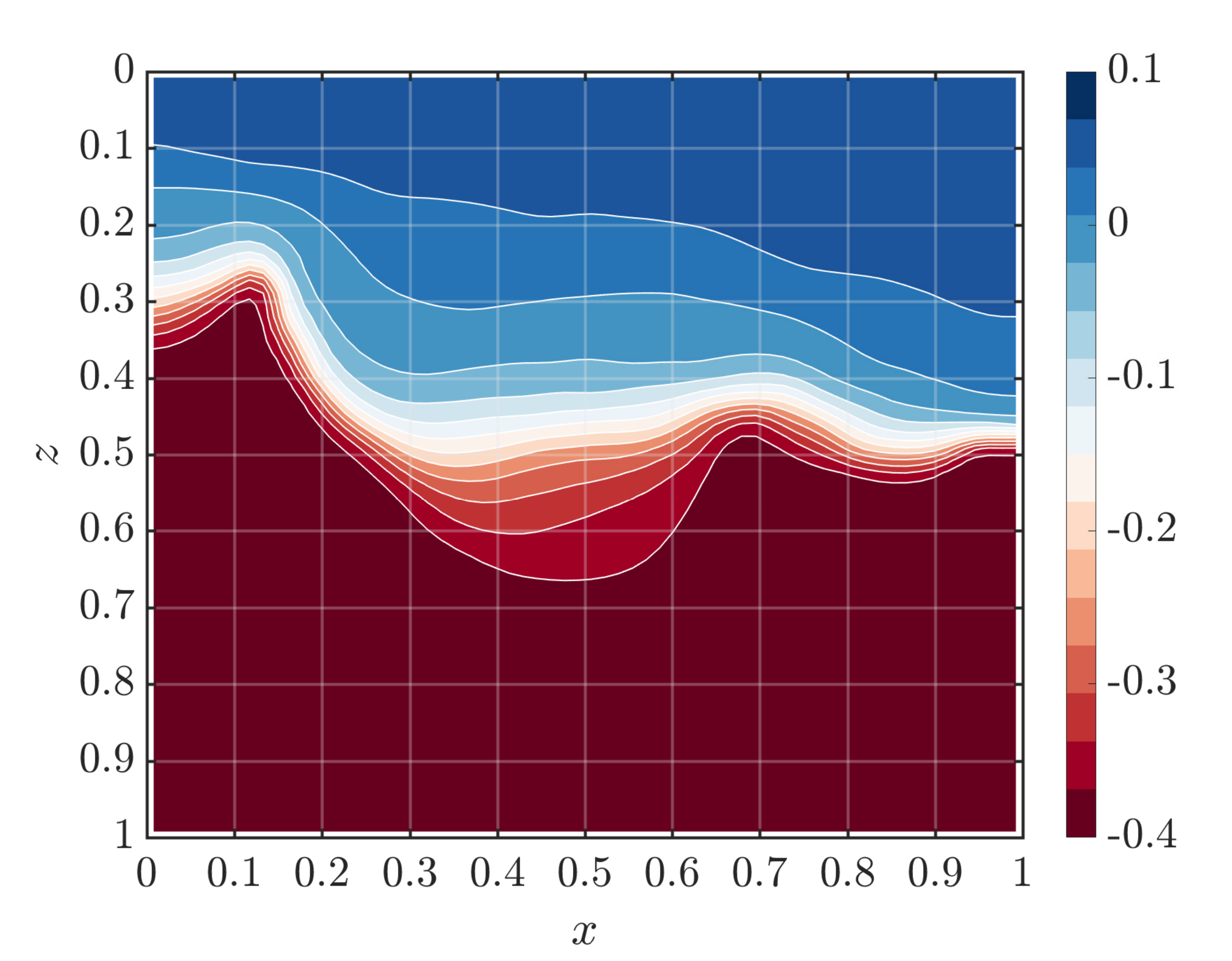}
\subcaption{$n^{(bl)} = 1.45,\alpha^{(bl)}=3.0$}
\end{subfigure}
\begin{subfigure}[b]{0.32\textwidth}
\includegraphics[trim={0cm 0cm 0cm 0cm},clip,scale = 0.205]{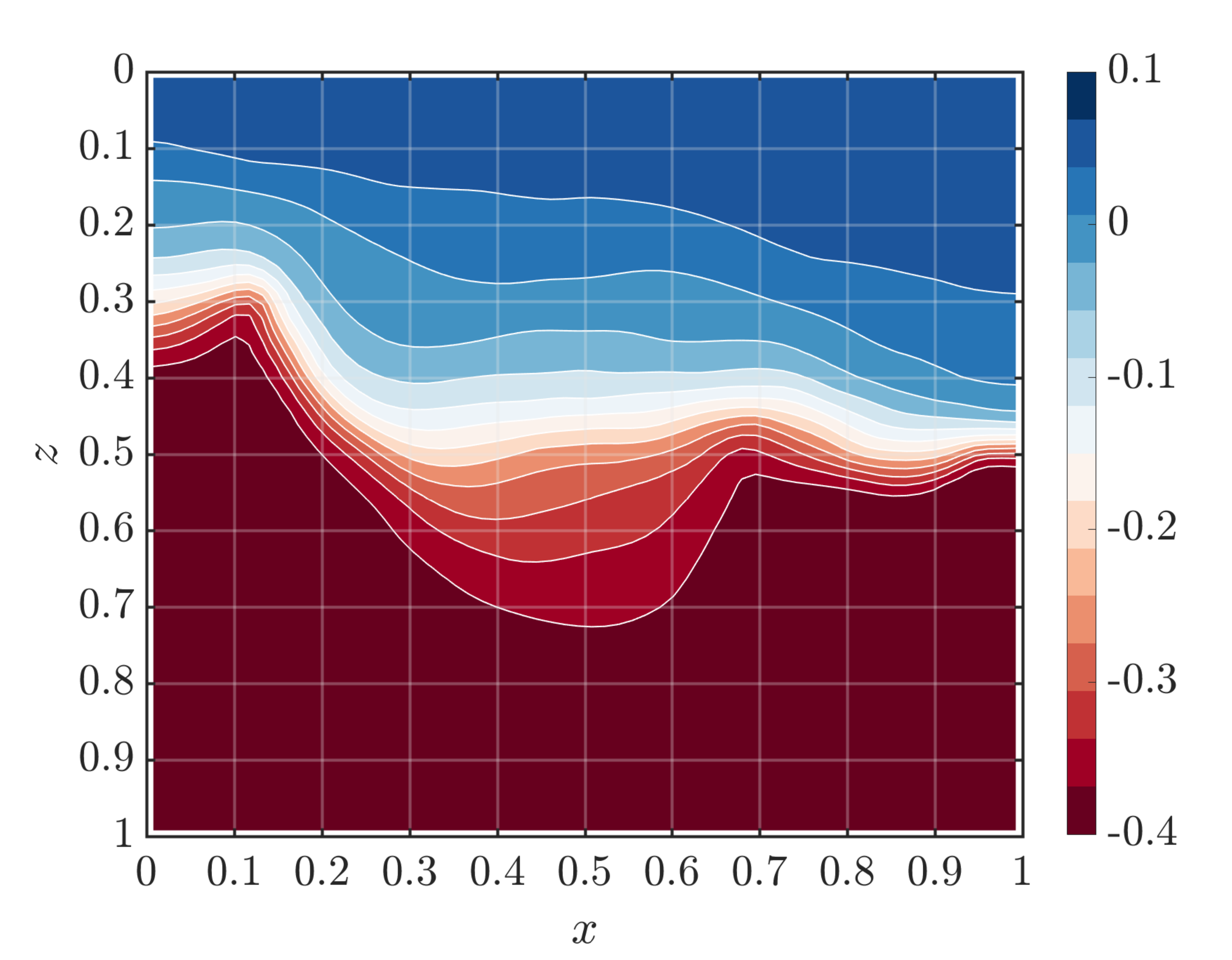}
\subcaption{$n^{(bl)} = 1.65,\alpha^{(bl)}=2.6$}
\end{subfigure}
\begin{subfigure}[b]{0.33\textwidth}
\includegraphics[trim={0cm 0cm 0cm 0cm},clip,scale = 0.205]{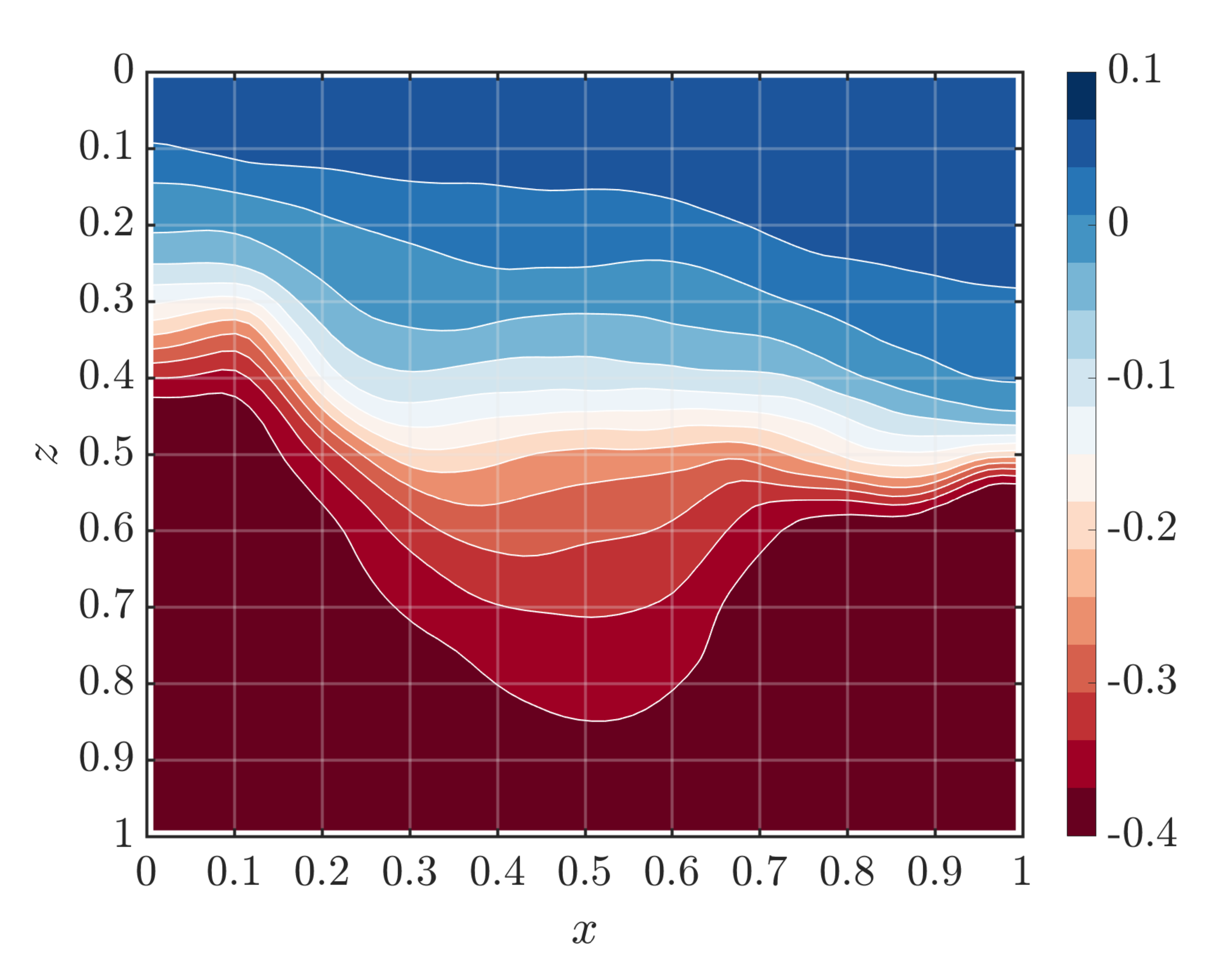}
\subcaption{$n^{(bl)} = 1.85,\alpha^{(bl)}=2.2$}
\end{subfigure}
\caption{Comparison of pressure head fields at $T_{final} = 0.2$ [h] for different baseline values of the parameters, $(n^{(bl)},\alpha^{(bl)})$ but with the same random fields $Z,\varepsilon_{\theta_s},\varepsilon_{\theta_r}, \varepsilon_{\alpha},\varepsilon_n$. Solutions are based on $h=\Delta t = 1/64$ and the Mat\'ern parameter set $\Phi_1$.}\label{nn_alpha_comp}
\end{figure}

\begin{figure}
\centering
\includegraphics[trim={0cm 0cm 0cm 0cm},clip,scale = 0.22]{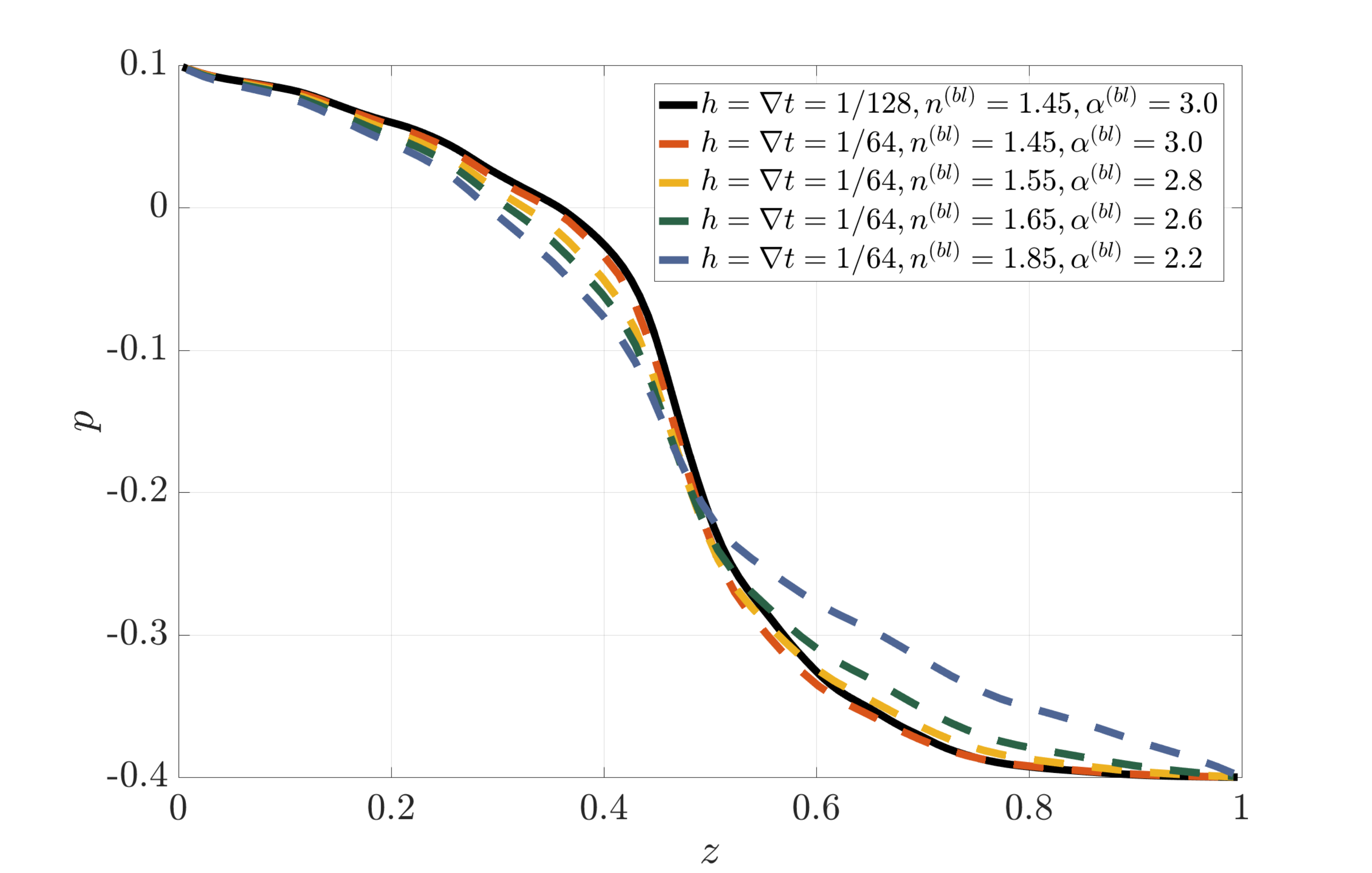}
\caption{Comparison of cross sections of the pressure heads from Figure \ref{nn_alpha_comp} at $x=0.5$.}\label{n_alpha_comp2}
\end{figure}

We first investigate the correlations for the pressure head profiles when the baseline values for $\alpha^{(bl)}$ and $n^{(bl)}$ are varied however employing the same random fields. In Figure \ref{nn_alpha_comp}, we compare three pressure head solutions with different baseline values, and with the same random fields, $Z(\mathbf{x},\omega),\varepsilon_{\theta_s}(\mathbf{x},\omega),\varepsilon_{\theta_r}(\mathbf{x},\omega), \varepsilon_{\alpha}(\mathbf{x},\omega)$ and $\varepsilon_n(\mathbf{x},\omega)$ (see Section \ref{stochRE}). Clearly, the pressure head profile becomes more \emph{diffusive} when ``easier'' parameters are prescribed. We also compare the cross sections of the pressure head profiles at $x = 0.5$ in Figure \ref{n_alpha_comp2}. For reference, we use the solution on the fine grid $h=\Delta t=1/128$ (black solid line) and compare it with different pairs of $n^{(bl)}$ and $\alpha^{(bl)}$ values on the next coarse grid $h=\Delta t=1/64$.  The profiles with the same $(n^{(bl)},\alpha^{(bl)})$-values are very close and the deviation increases as the two parameters are set to ``easier'' values. Thus, we can conclude that the correlation decays as the difference between the baseline values of the nonlinear parameters widens.

\begin{figure}
\begin{subfigure}[b]{0.49\textwidth}
\includegraphics[trim={0cm 0cm 0cm 0cm},clip,scale = 0.285]{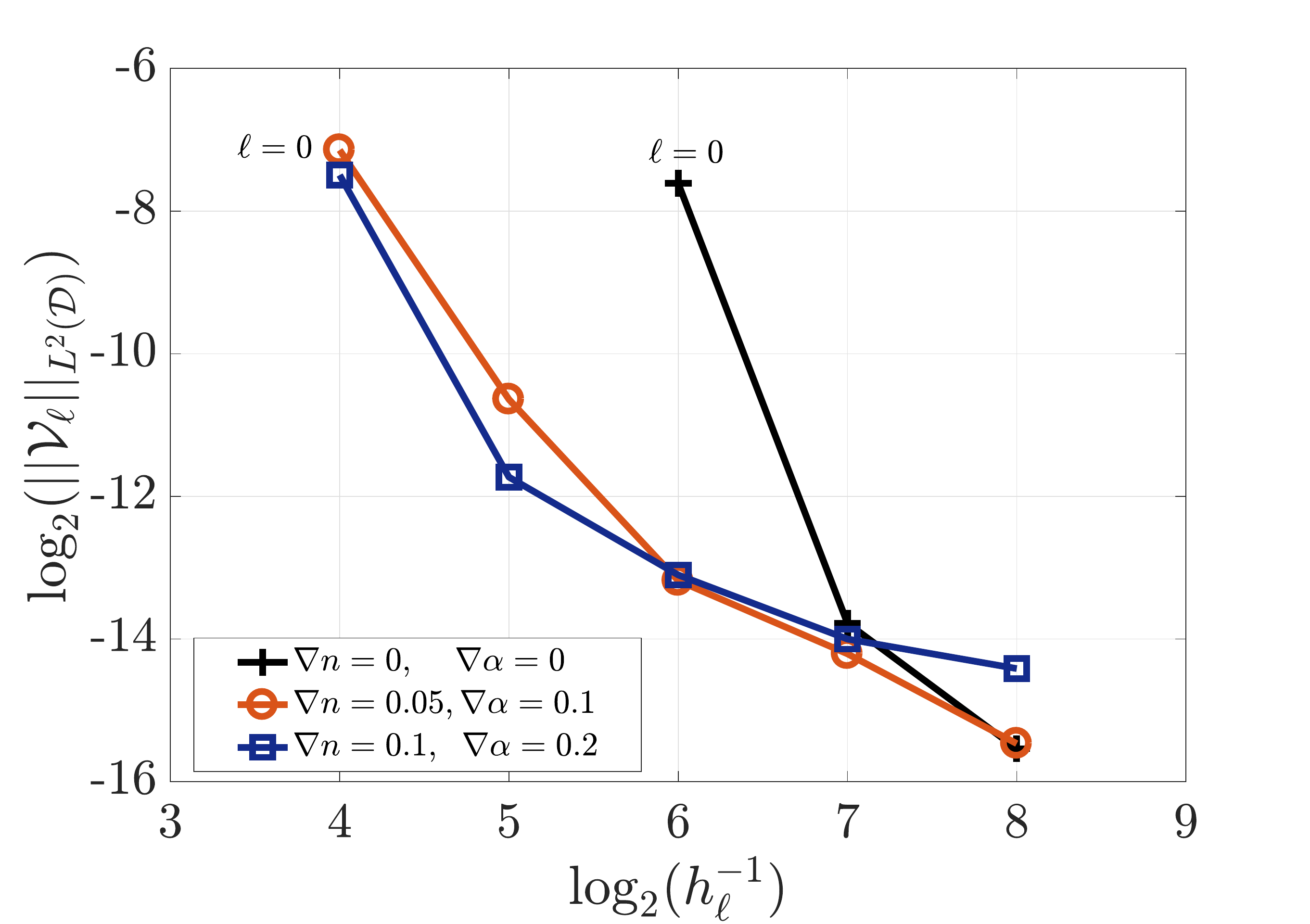}
\subcaption{$h_L=1/256,\Phi_1,\alpha^{(bl)}_L=3.0,n^{(bl)}_L=1.45$}
\end{subfigure}
\begin{subfigure}[b]{0.49\textwidth}
\includegraphics[trim={0cm 0cm 0cm 0cm},clip,scale = 0.28]{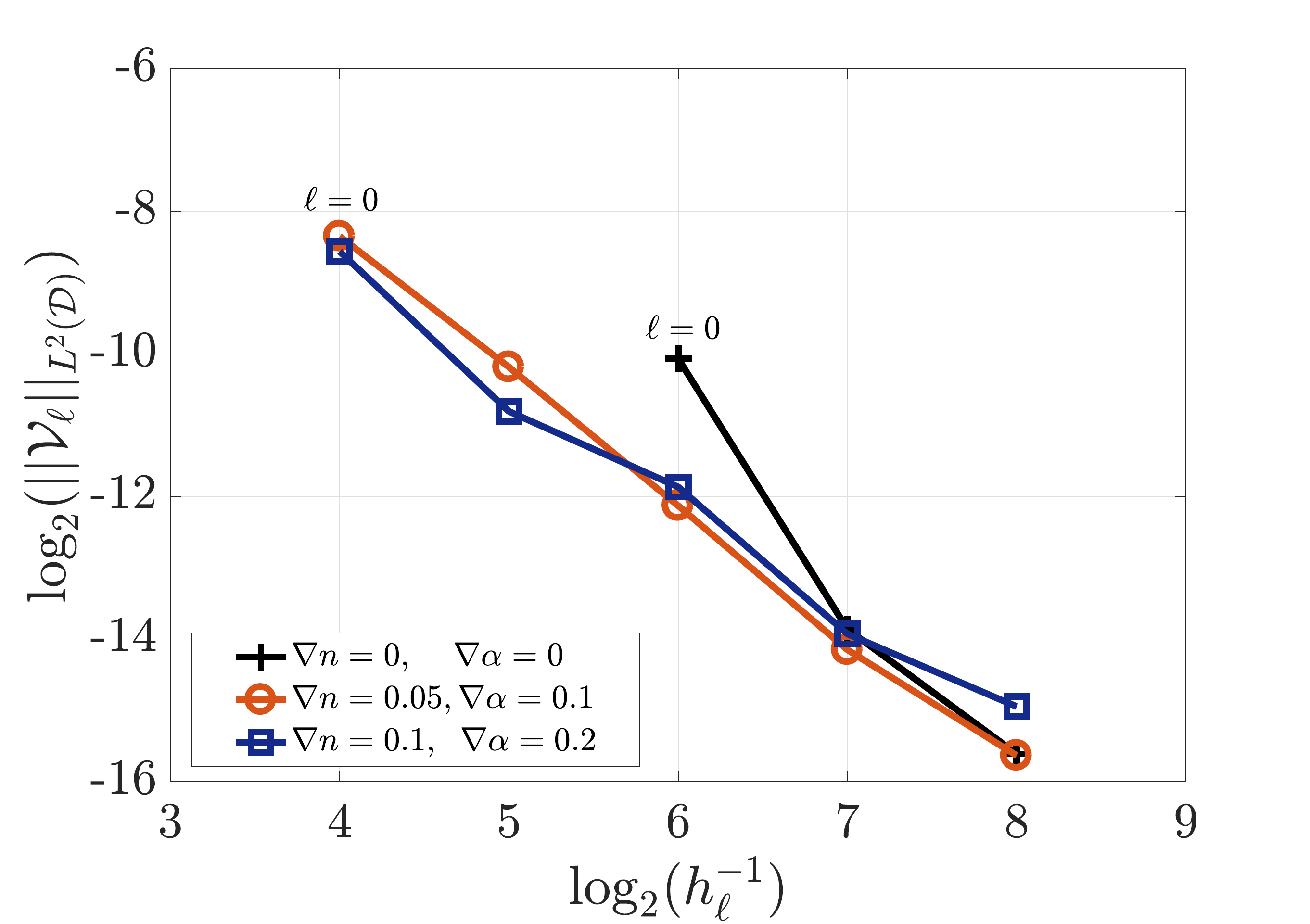}
\subcaption{$h_L=1/256,\Phi_2, \alpha^{(bl)}_L=2.8,n^{(bl)}_L=1.55$}
\end{subfigure}
\caption{Comparison of level-dependent variances using different pairs of $\nabla \alpha$,$\nabla n$ for the two Mat\'ern parameters.}\label{MLMC_var}
\end{figure}

Next, we study the behavior of the level-dependent variance $\lnorm\mathcal{V}_\ell\rnorm_{\Ld}$ when using the parametric continuation approach. For this we define the so-called {\em parametric continuation variables}, $\nabla \alpha = \alpha_{\ell} - \alpha_{\ell-1}$ and $\nabla n = n_{\ell-1} - n_{\ell}$, with the purpose to reduce the nonlinearity when processing coarse grids. In Figure \ref{MLMC_var}, we plot $\lnorm\mathcal{V}_\ell\rnorm_{\Ld}$ computed using \eqref{level_var_ch6} for different $(\nabla \alpha$,  $\nabla n)$ pairs for the two Mat\'ern parameter sets. The original problem is solved with $h_L = 1/256$. The parameter sets $\mathbf{\Theta}_L$ and $\mathbf{\Theta}_\ell$, for $\ell = L-1,L-2,...,0$, are obtained by employing $\nabla \alpha $ and $\nabla n$. The black line represents the variance when same problem is solved on all levels, i.e. $\nabla \alpha=\nabla n=0$, corresponding to the original MLMC estimator. Using the original approach, we can only process three levels in the MLMC hierarchy. The red and blue lines in the figure correspond to the variance which is computed using $\nabla \alpha=0.05,\nabla n=0.1$ and  $\nabla \alpha=0.1,\nabla n=0.2$, respectively.  For these two cases, we can incorporate a larger number of  coarse levels, up to $h_0=1/16$, as milder nonlinear problems are solved on these coarse levels. Furthermore, for levels $\ell<L$ the variance is smaller, compared to the case without continuation (where $\nabla \alpha=\nabla n=0$) which will result in lower number of samples on these levels. Here, we wish to highlight the fact that choosing optimal values for $\nabla\alpha$ and $\nabla n$ is important. For example, when $\alpha=0.1,\nabla n=0.2$, the variance on the finest level increases as compared to the variance found with the original MLMC approach. Due to this, an increasing number of samples will be needed on the finest level, compared to the original MLMC estimator, which is undesirable as it may result in an expensive estimator. On the other hand, for a smaller perturbation $\nabla\alpha=0.05,\nabla n=0.1$, the magnitude of the variance on the finest level is similar to that of the original MLMC estimator. Therefore, the number of samples on the finest level will be more-or-less similar to the original MLMC estimator. 

Alternatively, one can avoid high variance on the finest level by using zero perturbations on the two finest levels, i.e. $\nabla \alpha_L=\nabla \alpha_{L-1} = 0$ and $\nabla n_L=\nabla n_{L-1} = 0$ and choosing non-zero perturbations on next coarser levels from $L-2$ onwards.  This way, we solve the original problem on the two finest levels. 

\begin{figure}
\begin{subfigure}[b]{0.49\textwidth}
\includegraphics[trim={0cm 0cm 0cm 0cm},clip,scale = 0.285]{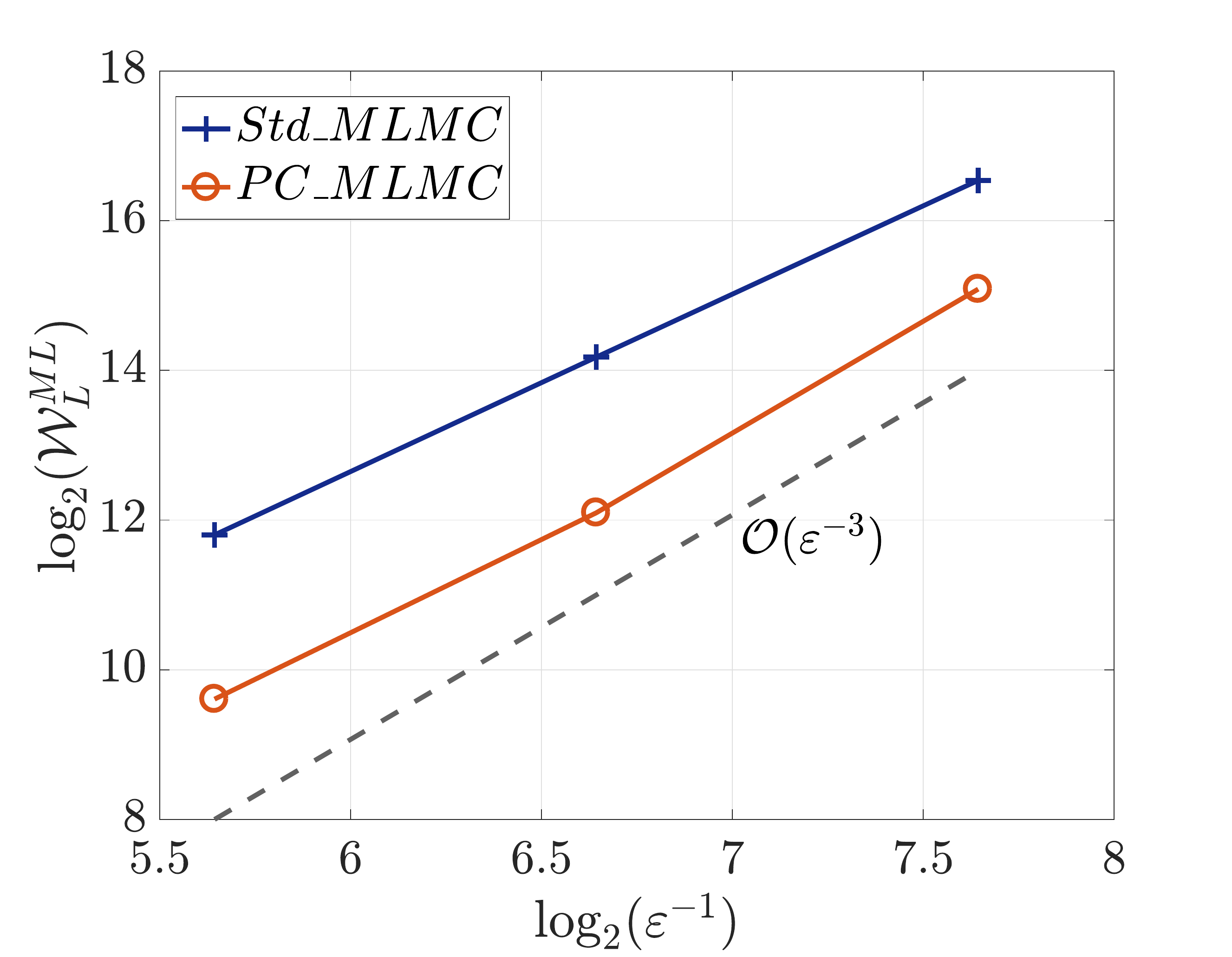}
\subcaption{$\Phi_1,\alpha^{(bl)}_L=3.0,n^{(bl)}_L=1.45$}
\end{subfigure}
\begin{subfigure}[b]{0.49\textwidth}
\includegraphics[trim={0cm 0cm 0cm 0cm},clip,scale = 0.28]{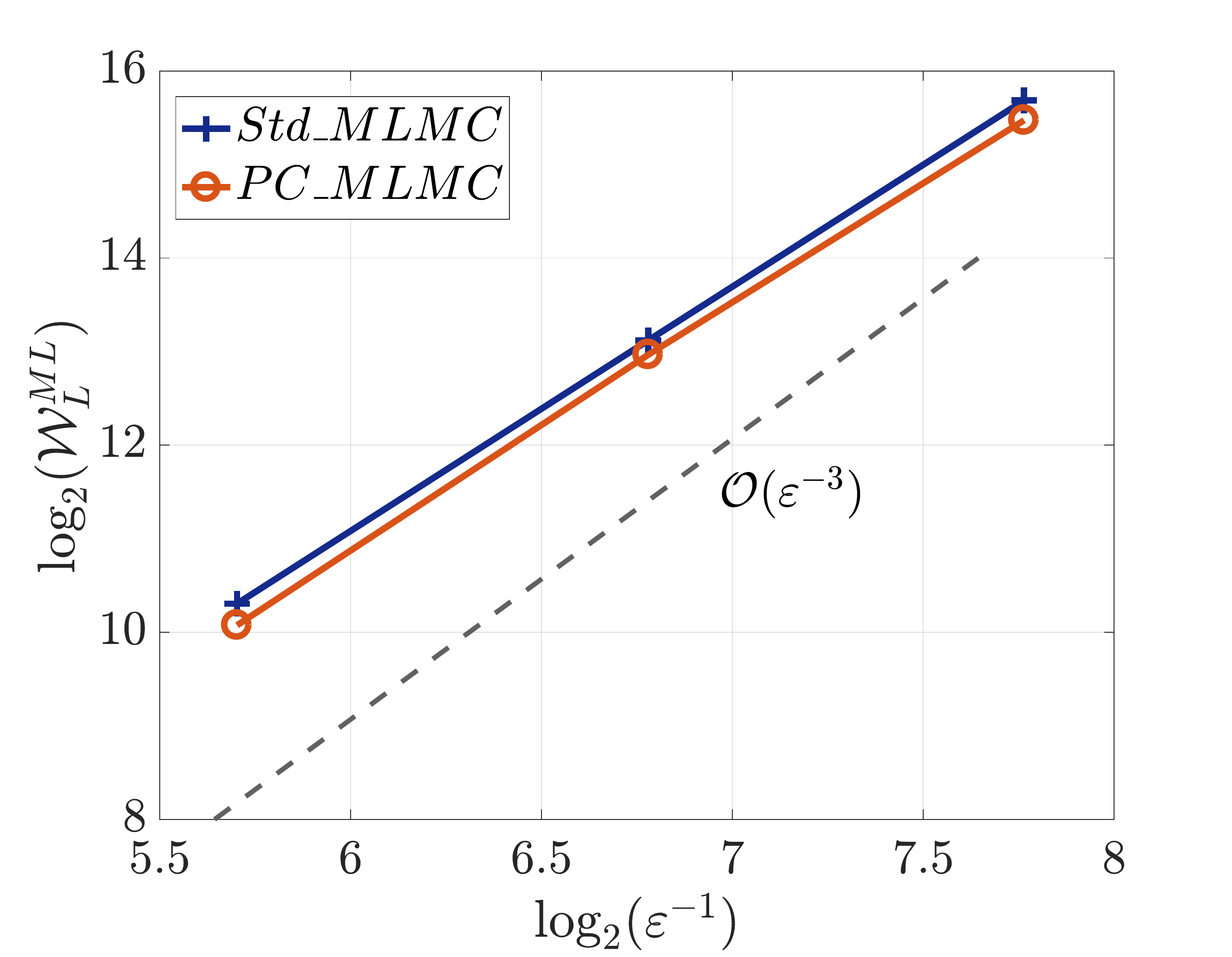}
\subcaption{$\Phi_2, \alpha^{(bl)}_L=2.8,n^{(bl)}_L=1.55$}
\end{subfigure}
\caption{Comparison of CPU times $\mathcal{W}_L^{ML}$ (sec) for two different estimators.}\label{cost_compare_pc_std}
\end{figure}

Now we apply the parametric continuation based MLMC estimator denoted by $PC\_MLMC$, to compute the mean and variance of the pressure head field and also perform comparisons with respect to the standard MLMC estimator which is denoted by $Std\_MLMC$. For this, we use the two previously discussed test cases: isotropic covariance $\Phi_1$ with baseline values $n^{(bl)} = 1.45,\alpha^{(bl)}=3.0$ and anisotropic covariance $\Phi_2$ with baseline values $n^{(bl)} = 1.55,\alpha^{(bl)}=2.8$. For simplicity, we use the continuation variables  $\nabla \alpha=0.05,\nabla n=0.1$ for both the isotropic and anisotropic cases. We compare the number of samples needed on different grids for three values of the tolerances. For the $Std\_MLMC$ estimator, the coarsest possible level is $h_0=\Delta t_0=1/64$, whereas for the $PC\_MLMC$, we use $h_0=\Delta t_0 =1/16$. As the discretization error shows a first-order decay (see Figure \ref{FV_err_ch6}), we set the tolerance $\varepsilon = \mathcal{O}(h_L)$. We use  Algorithm \ref{PC_MLMC} to reduce the sampling error to $\varepsilon$. In Table \ref{phi1_ch6}, the two estimators for the isotropic Mat\'ern parameter $\Phi_1$ are compared.  Due to the sample optimization strategy \eqref{MLMCSAMP_ch6},  a large number of samples is shifted to coarser grids when using the $PC\_MLMC$ estimator. Furthermore, a fewer number of samples are required for the $PC\_MLMC$ estimator compared to the $Std\_MLMC$, even on the finest level.  This is due the fact that the sum $\sum_{\ell=0}^L\sqrt{\lnorm\mathcal{V}_\ell\rnorm_{\Ld}\mathcal{W}_\ell}$ for the $PC\_MLMC$ estimator is slightly smaller than for $Std\_MLMC$. Moreover, a large computational gain is induced by the reduction in the number of samples on grid $h_\ell=1/64$, for instance, for $\varepsilon=0.005$, the number of samples reduced from 438 to 35 when using the parametric continuation. In Figure \ref{cost_compare_pc_std} (a) the CPU times for the two estimators are also compared. We observe a speed-up of about a factor of three for $\varepsilon=0.005$. 

A similar test is performed for the anisotropic problem. The number of samples for different tolerances are provided in Table \ref{phi2_ch6} and the CPU times in Figure \ref{cost_compare_pc_std} (b). Again some improvement in computation times are observed, although the gain is not as pronounced as for the first problem. This is due to the fact that the second case uses simpler baseline values $n^{(bl)} = 1.55,\alpha^{(bl)}=2.8$ and the cost reduction with parameter simplification is not very rapid. For the isotropic case with $n^{(bl)} = 1.45,\alpha^{(bl)}=3.0$, the cost decay is more rapid with parameter simplification. This is more evident from the cost map in Figure \ref{HM_Phi2}, where we see more dense contour lines around $n^{(bl)}= 1.45$. Therefore, the parametric continuation approach is very effective when a strongly nonlinear stochastic problem needs to be solved. 

We also wish to highlight the fact that both MLMC estimators are \emph{optimal} since the cost scales as $\mathcal{O}(\varepsilon^{-3})$, which is similar to the computational complexity of solving one deterministic problem on the finest grid, i.e. $\mathcal{O}(h_L^{-3})$ and $h_L=\mathcal{O}(\varepsilon)$.
\begin{table}[H]
\centering
\begin{tabular}{|c|c|c|c|c|c|c|}
\hline
\multirow{2}{*}{$h_\ell$} &
\multicolumn{2}{c|}{$N_{\ell}(h_L=1/64,\varepsilon = 0.02 )$} &
\multicolumn{2}{c|}{$N_{\ell}(h_L=1/128,\varepsilon = 0.01)$} &
\multicolumn{2}{c|}{$N_{\ell}(h_L=1/256,\varepsilon = 0.005)$} \\
\cline{2-7}
    & $Std\_MLMC$ & $PC\_MLMC$ & $Std\_MLMC$ &$ PC\_MLMC$ & $Std\_MLMC$ &$PC\_MLMC$ \\
    \hline
    1/16 &$-$ &115& $-$&459 & $-$ & 1833\\
    \hline
    1/32 &$-$ &11&$-$ &44 & $-$ & 176\\
    \hline
    1/64 &28 &3& 110&9 &438 &35\\
    \hline
    1/128&$-$ & $-$&5 & 2&18 &8\\
    \hline
    1/256 &$-$&$-$&$-$ &$-$ &4& 3\\
    \hline
\end{tabular}
\caption{Comparison of number of samples needed to achieve tolerances $\varepsilon$ using the standard MLMC ($Std\_MLMC$) and parametric continuation MLMC ($PC\_MLMC$) estimators for $\Phi_1,n_L^{(bl)} = 1.45,\alpha_L^{(bl)}=3.0$. Entries with symbol $(-)$ indicate zero samples needed for that grid.}\label{phi1_ch6}
\end{table}
\begin{table}[H]
\centering
  \begin{tabular}{|c|c|c|c|c|c|c|}
    \hline
    \multirow{2}{*}{$h_\ell$} &
      \multicolumn{2}{c|}{$N_{\ell}(h_L=1/64,\varepsilon = 0.0184)$} &
      \multicolumn{2}{c|}{$N_{\ell}(h_L=1/128,\varepsilon = 0.0092)$} &
      \multicolumn{2}{c|}{$N_{\ell}(h_L=1/256,\varepsilon = 0.0046)$} \\
      \cline{2-7}
    & $Std\_MLMC$ & $PC\_MLMC$ & $Std\_MLMC$ &$ PC\_MLMC$ & $Std\_MLMC$ &$PC\_MLMC$ \\
    \hline
    1/16 &$-$ &96 &$-$ &427 &$-$ &1659 \\
    \hline
    1/32 &$-$ &18 &$-$ &77 &$-$ & 301\\
    \hline
    1/64 & 9&4 &36 &16 & 140&60 \\
    \hline
    1/128 &$-$ &$-$ &4 & 4&15 &12 \\
    \hline
    1/256 &$-$&$-$&$-$&$-$&3 & 3\\
    \hline
  \end{tabular}
\caption{Comparison of number of samples needed to achieve tolerances $\varepsilon$ using the standard MLMC ($Std\_MLMC$) and parametric continuation MLMC ($PC\_MLMC$) estimators for $\Phi_2,n^{(bl)} = 1.55,\alpha^{(bl)}=2.8$.}\label{phi2_ch6}
\end{table}

\begin{figure}[h]
\begin{subfigure}[b]{0.33\textwidth}
\includegraphics[trim={0cm 0cm 0cm 0cm},clip,scale = 0.175]{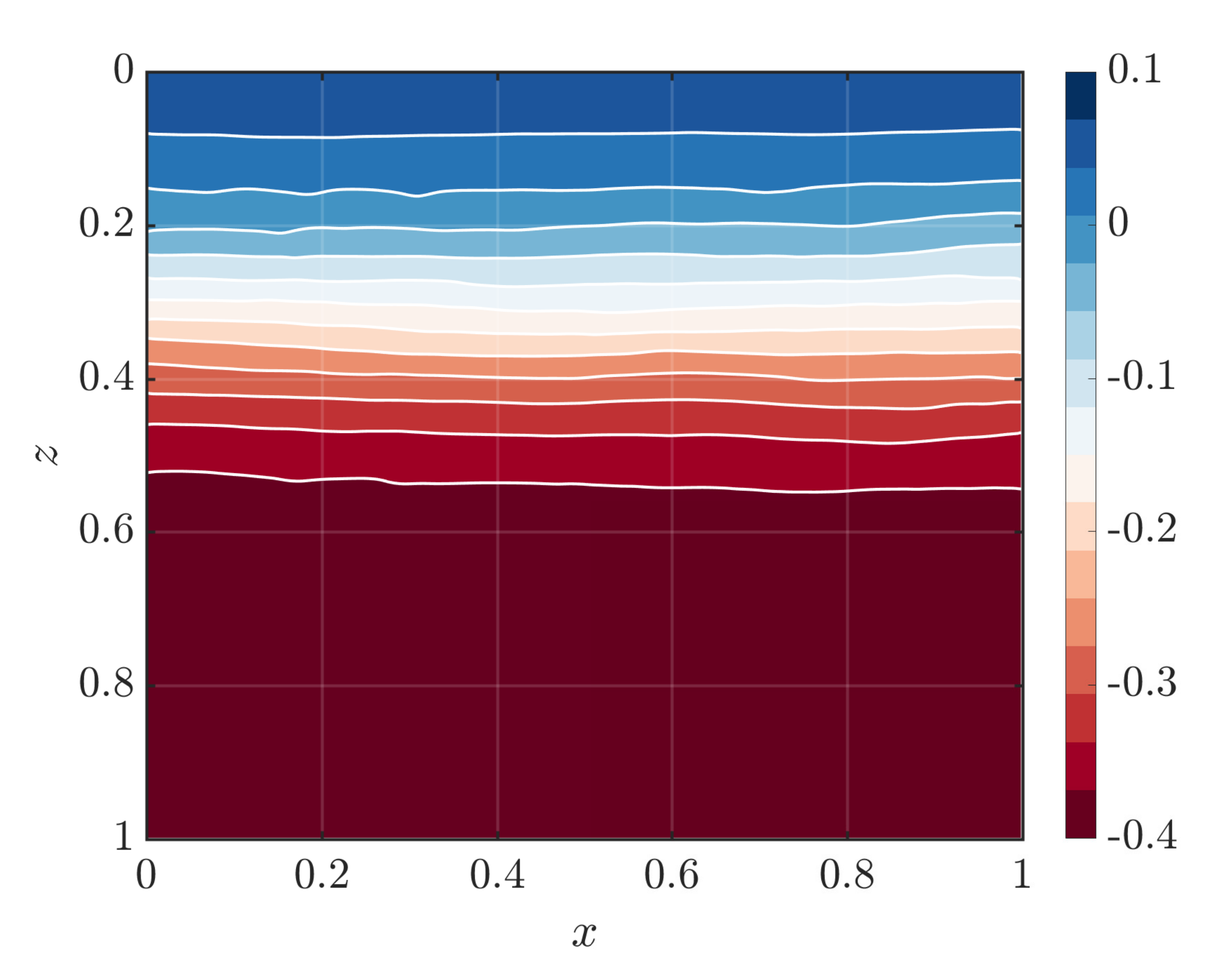}
\subcaption{$Std\_MLMC$}
\end{subfigure}
\begin{subfigure}[b]{0.33\textwidth}
\includegraphics[trim={0cm 0cm 0cm 0cm},clip,scale = 0.175]{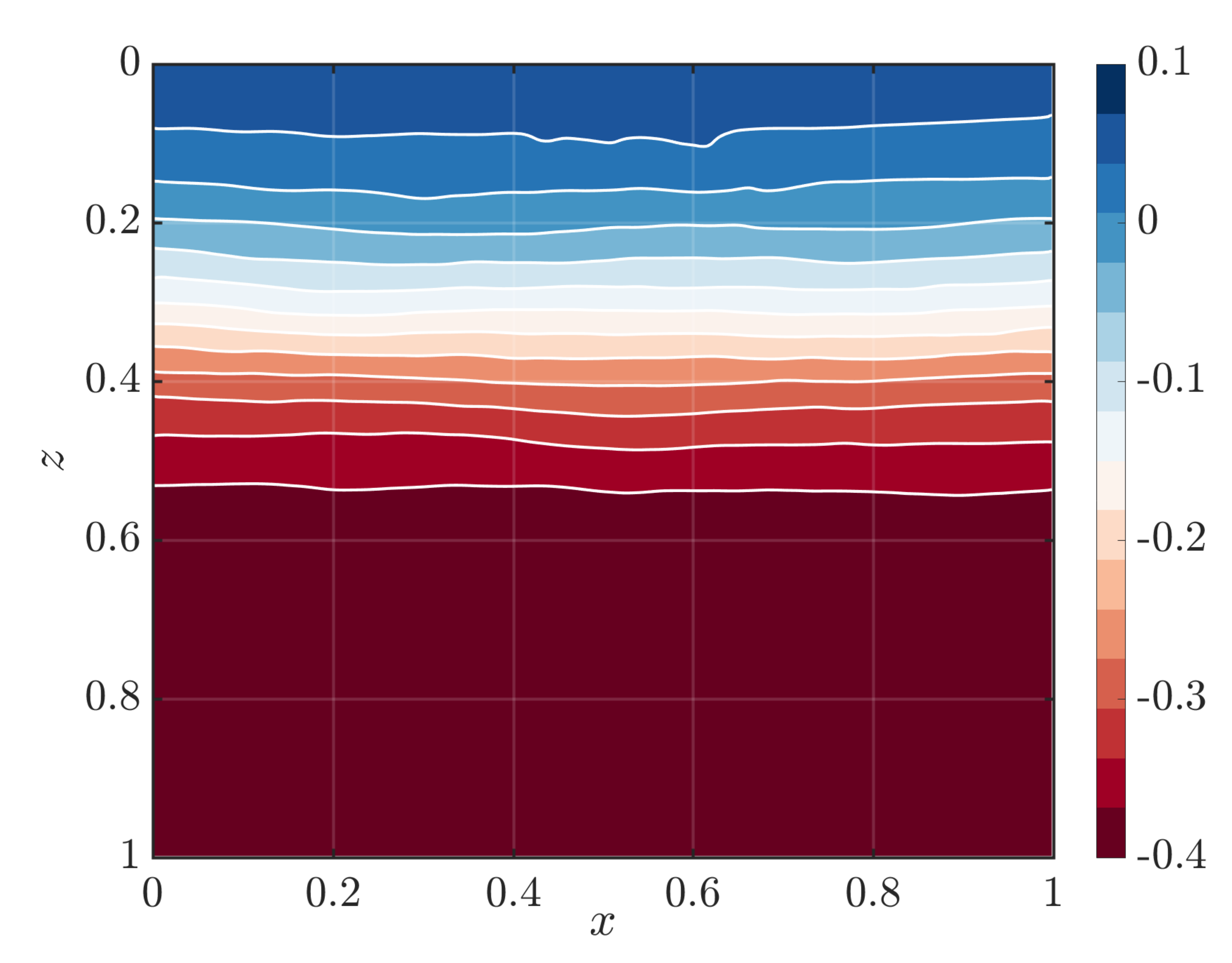}
\subcaption{$PC\_MLMC$}
\end{subfigure}
\begin{subfigure}[b]{0.32\textwidth}
\includegraphics[trim={0cm 0cm 0cm 0cm},clip,scale = 0.175]{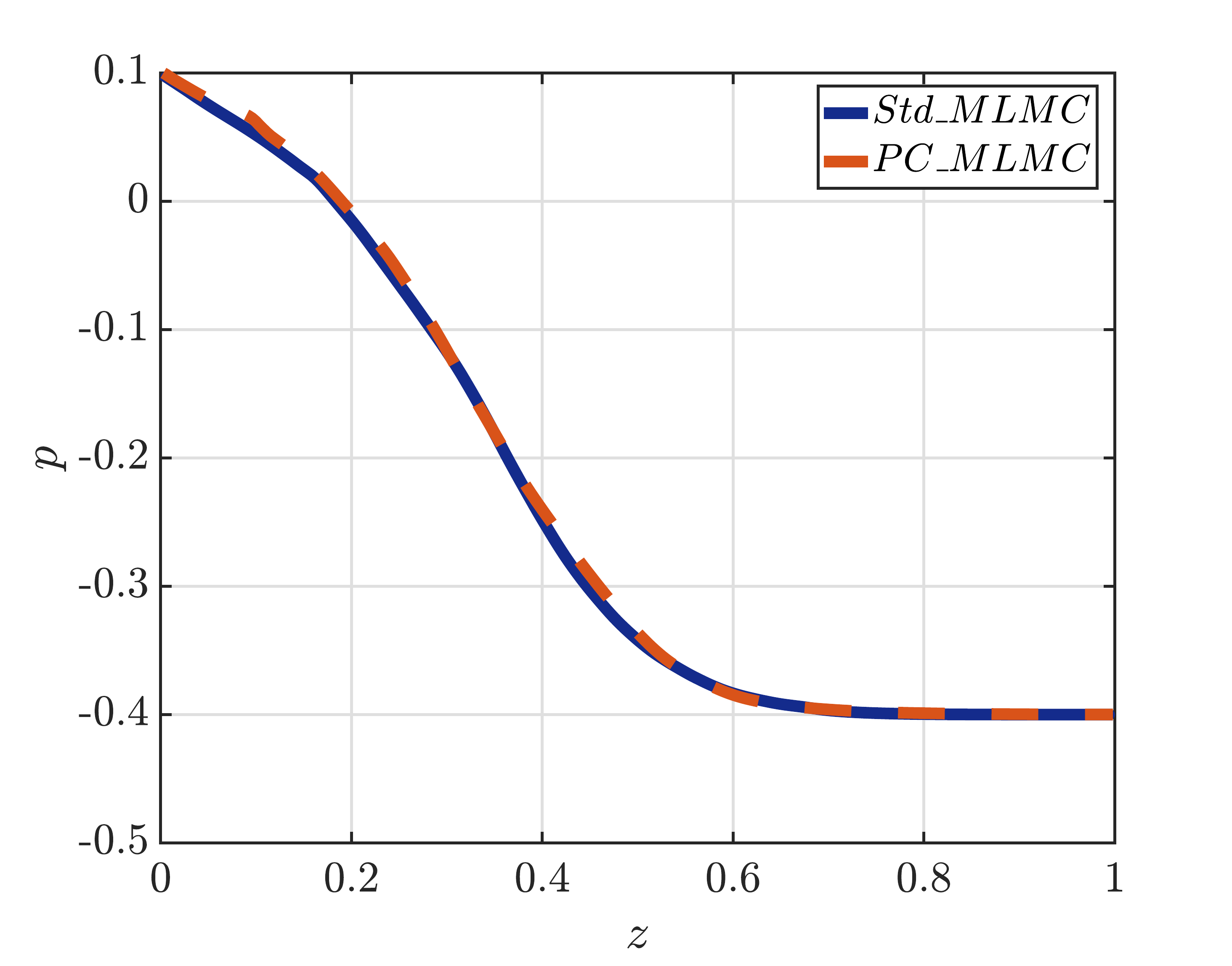}
\subcaption{Mean pressure head profile at $x=0.5$.}
\end{subfigure}
\caption{Comparison of mean pressure head field for $\Phi_1,\alpha^{(bl)}_L=3.0,n^{(bl)}_L=1.45, T_{final} = 0.2$ [h] computed using the two MLMC estimators with finest level $h_L=\Delta t_L =1/256$ and $\varepsilon=0.005$.}\label{Phi1_mean_ch6}
\end{figure}

\begin{figure}[h]
\begin{subfigure}[b]{0.33\textwidth}
\includegraphics[trim={0cm 0cm 0cm 0cm},clip,scale = 0.175]{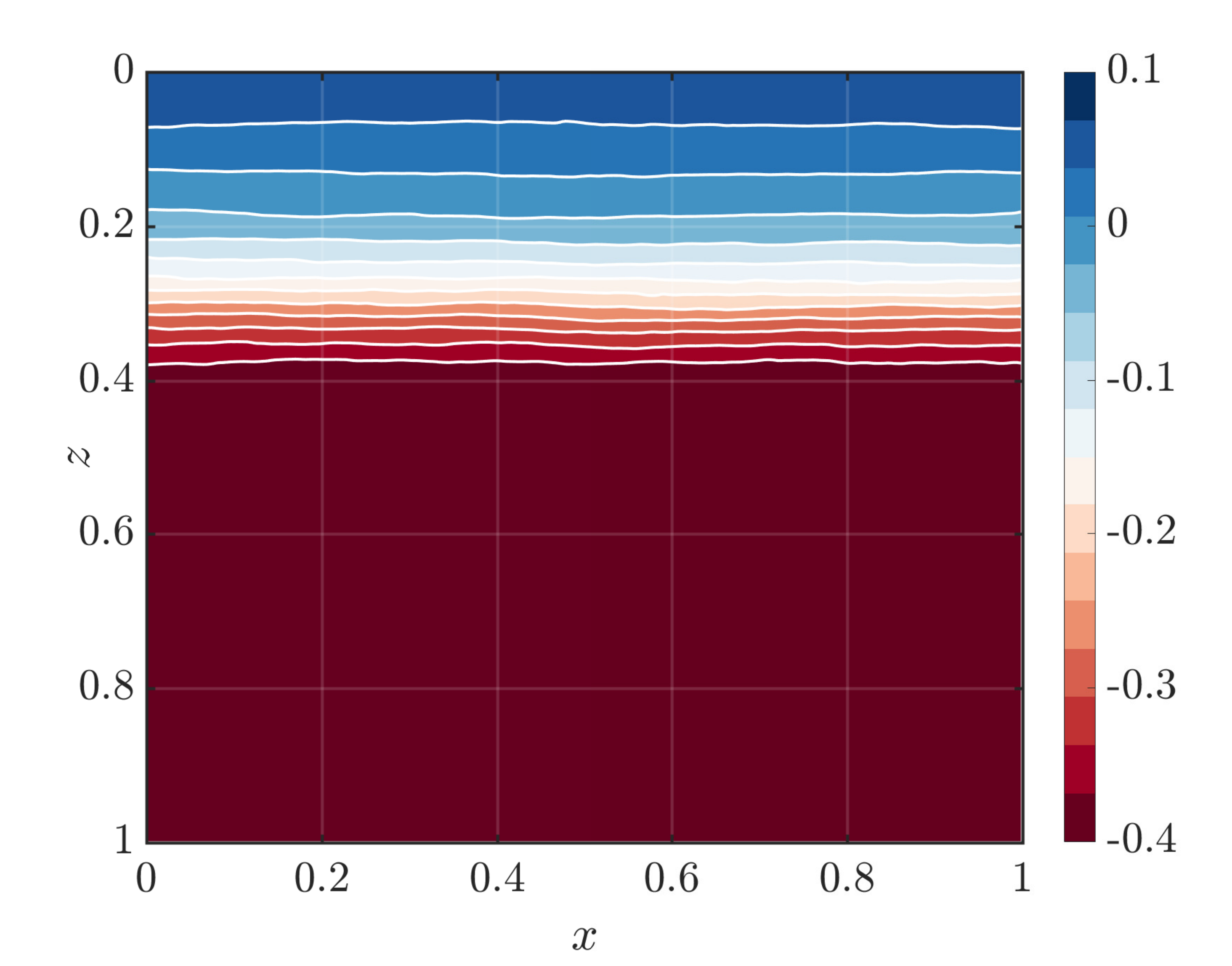}
\subcaption{$Std\_MLMC$}
\end{subfigure}
\begin{subfigure}[b]{0.33\textwidth}
\includegraphics[trim={0cm 0cm 0cm 0cm},clip,scale = 0.175]{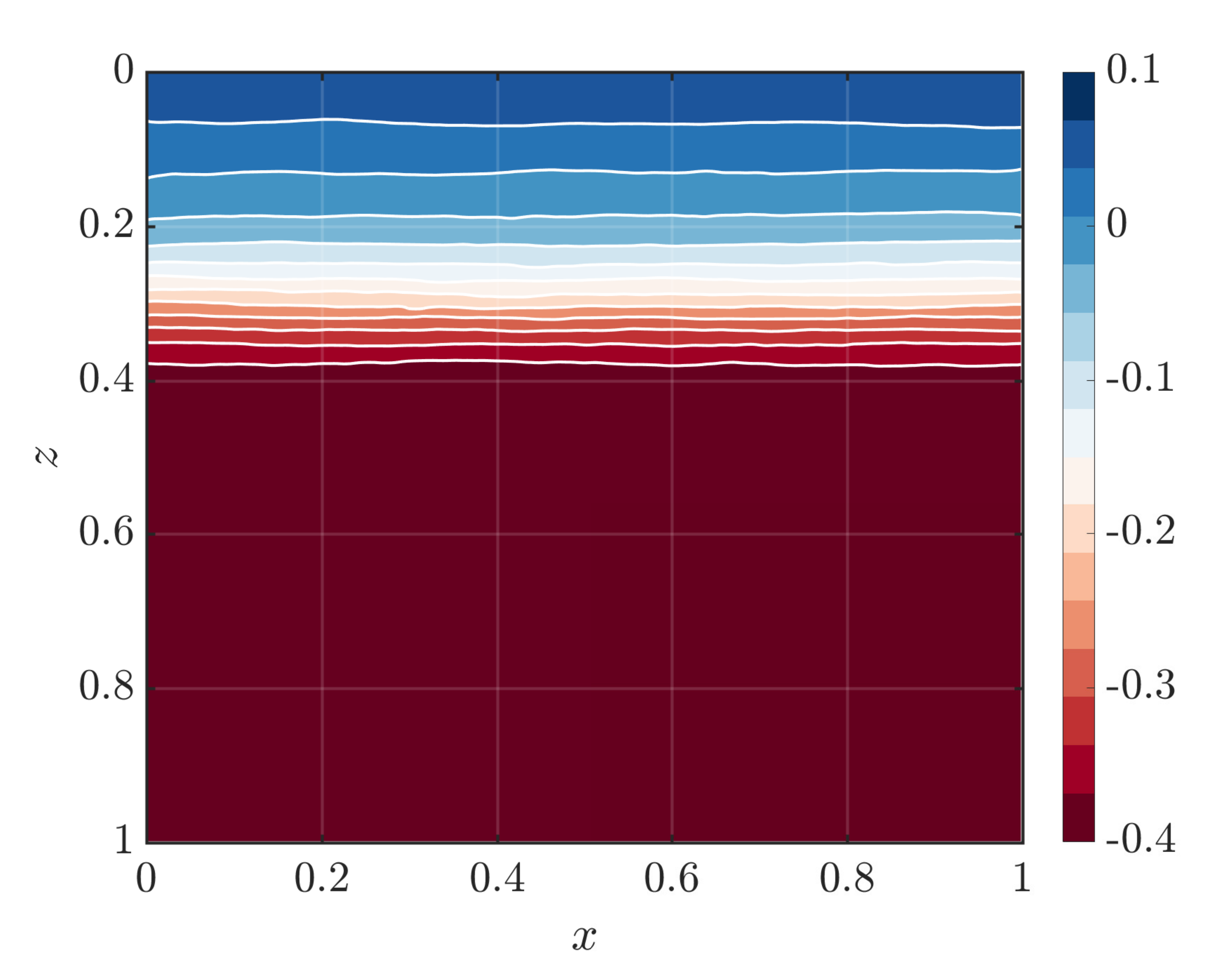}
\subcaption{$PC\_MLMC$}
\end{subfigure}
\begin{subfigure}[b]{0.32\textwidth}
\includegraphics[trim={0cm 0cm 0cm 0cm},clip,scale = 0.175]{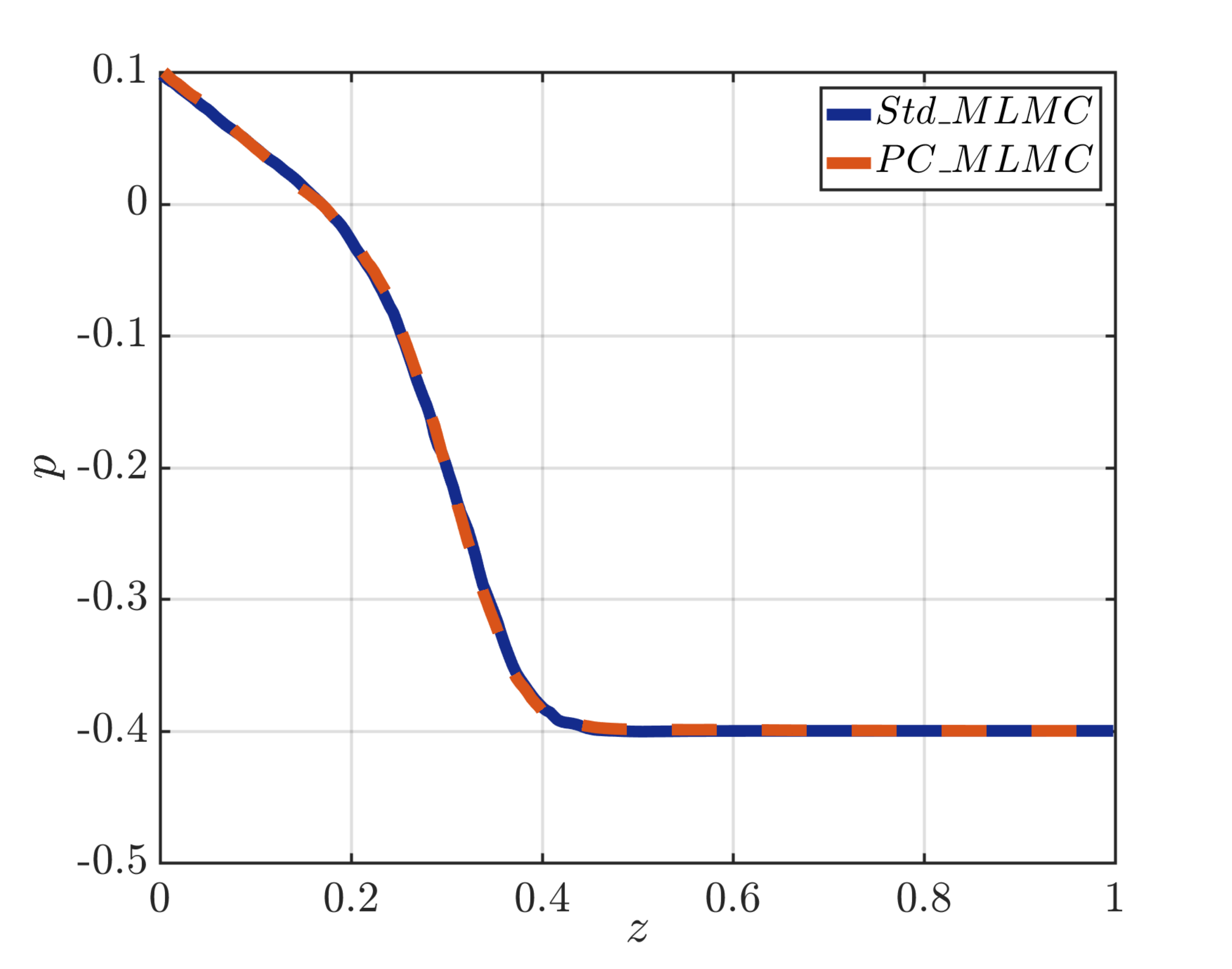}
\subcaption{Mean pressure head profile at $x=0.5$.}
\end{subfigure}
\caption{Comparison of mean pressure head field for $\Phi_2,\alpha^{(bl)}_L=2.8,n^{(bl)}_L=1.55, T_{final} = 0.2$ [h] computed using the two MLMC estimators with finest level $h_L=\Delta t_L =1/256$ and $\varepsilon=0.0046$.}\label{Phi2_mean_ch6}
\end{figure}

\begin{figure}
\begin{subfigure}[b]{0.33\textwidth}
\includegraphics[trim={0cm 0cm 0cm 0cm},clip,scale = 0.175]{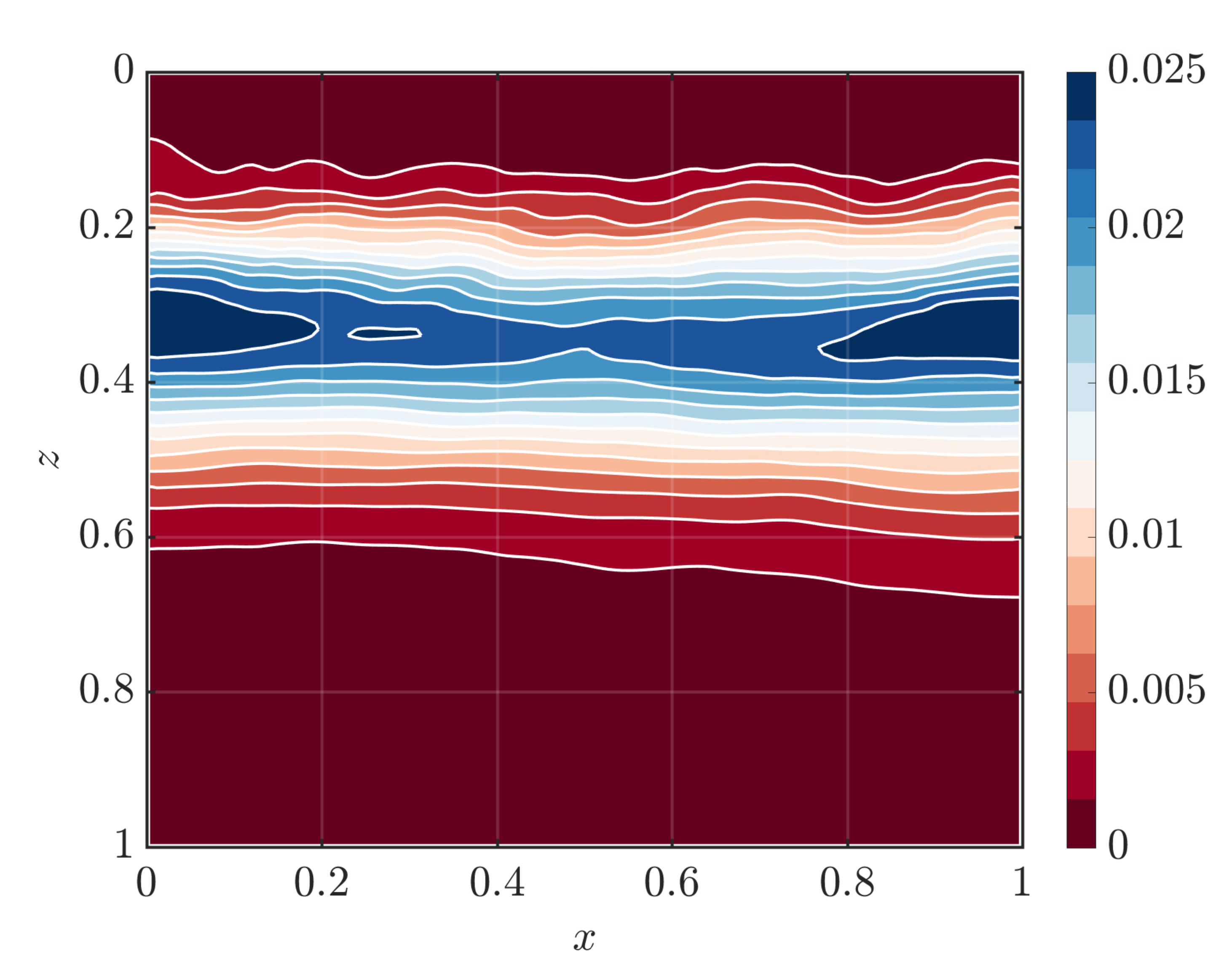}
\subcaption{$Std\_MLMC$}
\end{subfigure}
\begin{subfigure}[b]{0.33\textwidth}
\includegraphics[trim={0cm 0cm 0cm 0cm},clip,scale = 0.175]{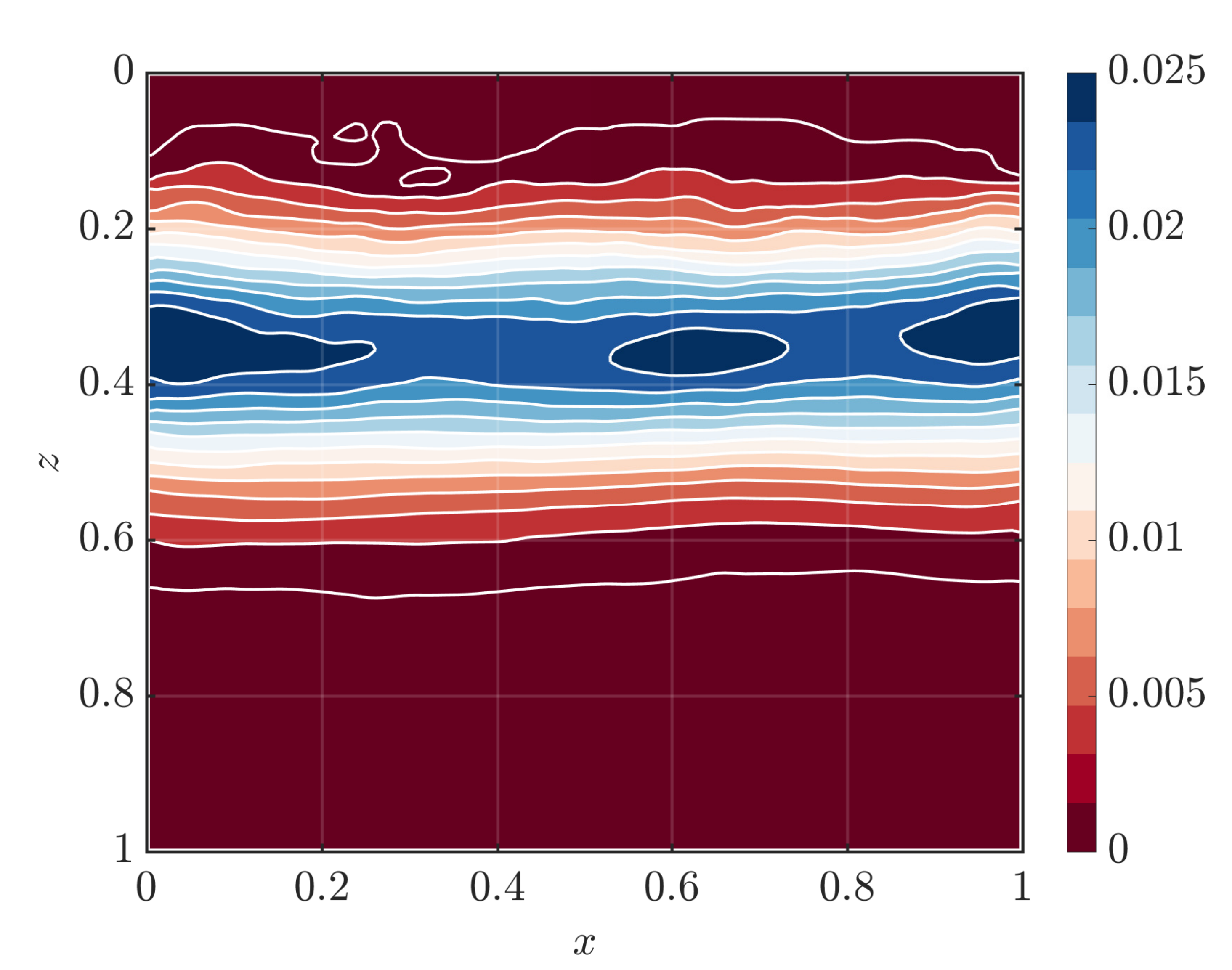}
\subcaption{$PC\_MLMC$}
\end{subfigure}
\begin{subfigure}[b]{0.32\textwidth}
\includegraphics[trim={0cm 0cm 0cm 0cm},clip,scale = 0.175]{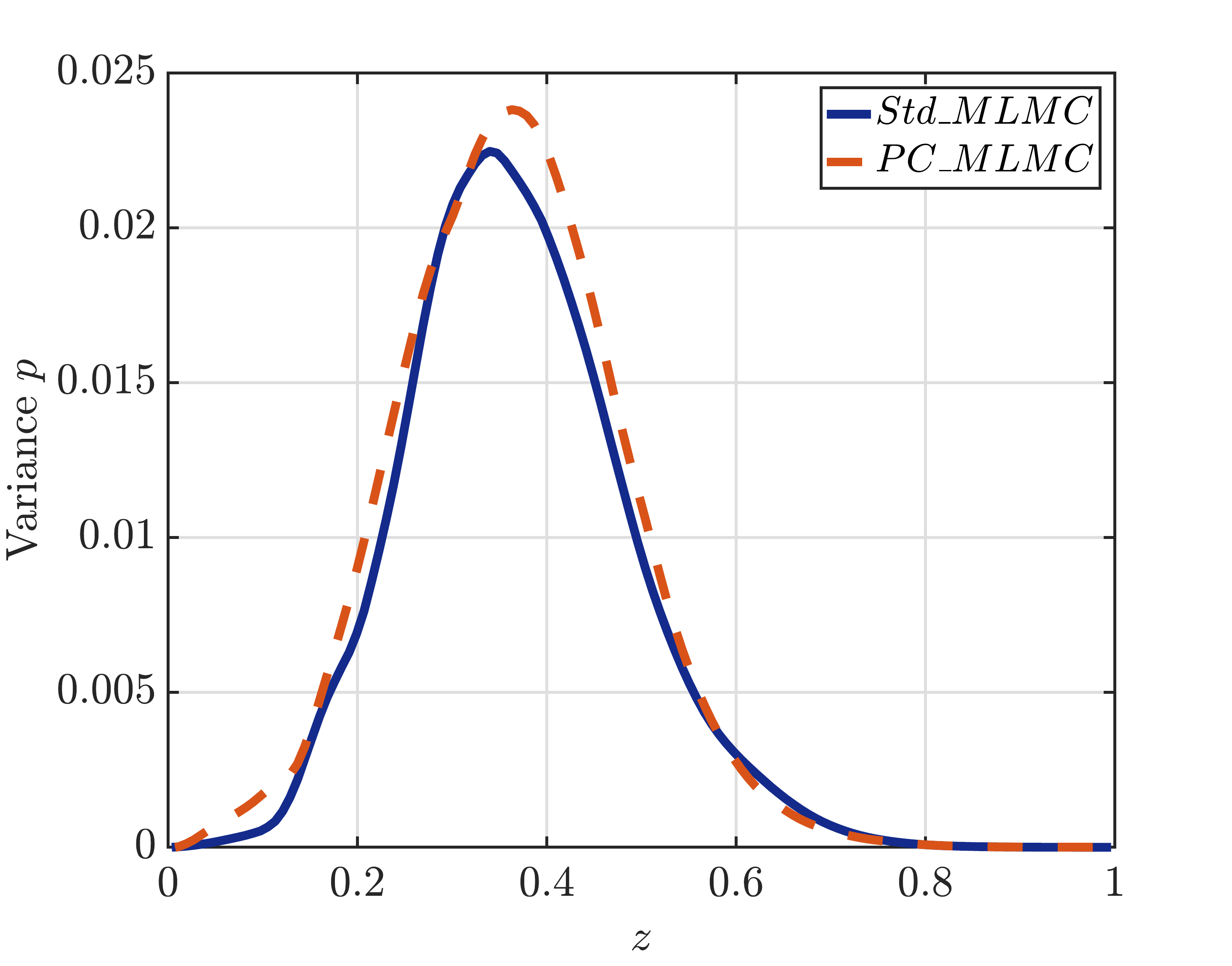}
\subcaption{Variance of pressure head profile at $x=0.5$.}
\end{subfigure}
\caption{Comparison of the variance of the pressure head field for $\Phi_1,\alpha^{(bl)}_L=3.0,n^{(bl)}_L=1.45, T_{final} = 0.2$ [h] computed using the two MLMC variance estimators with finest level $h_L=\Delta t_L =1/256$.}\label{Phi1_var_ch6}
\end{figure}

In the last part of this section, we validate the stochastic moments computed using the proposed estimator. It is expected that the mean pressure field computed using the two MLMC estimators should converge to a similar solution for a given tolerance. In Figure \ref{Phi1_mean_ch6}, the mean pressure head profile for the isotropic case is shown. It is computed using the number of samples from Table \ref{phi1_ch6}, with $\varepsilon =0.005$.  For a closer inspection, we also compare the mean pressure head profiles at $x=0.5$. Similarly, the mean profile for the anisotropic case is presented in Figure \ref{Phi2_mean_ch6}, using the number of samples from Table \ref{phi2_ch6} with $\varepsilon=0.0046$. We see good agreement between the mean profiles computed from the two estimators. The isotropic case exhibits a seemingly smoother transition from the saturated to the unsaturated zone compared to the anisotropic problem. In Figure \ref{Phi1_var_ch6} we also present the variance field for the isotropic test case, computed using the multilevel variance estimator $\mathcal{V}_L^{ML}[p_{h_L,\boldsymbol{\Theta}_L}] $ given in \eqref{ML_var_ch6}. The two variance fields are very similar, although some discrepancy in the magnitude is observed. This is due to the fact that the two variance fields are computed using the samples based on the error analysis of $\mathpzc{E}^{ML}_L[p_{h_L,\boldsymbol{\Theta}_L}]$ (from Table \ref{phi1_ch6}) and not on the error analysis of $\mathcal{V}_L^{ML}[p_{h_L,\boldsymbol{\Theta}_L}] $. Thus, the two variance estimates may have different tolerances resulting in this slight mismatch. Readers are referred to \cite{bierig} for a detailed error analysis of the multilevel variance estimator. 
%As the corrections in the multilevel variance estimator requires taking difference between variances computed from two consecutive levels, we may encounter negative variance at some locations. 

The results from the two estimators also showed good agreement with the plain Monte Carlo solutions performed on the grid $h_L=1/128$. This is done in order to verify that a proper upscaling of the random fields on coarser levels is carried out while using the MLMC estimator. These results are omitted for the sake of brevity.

\section{Conclusion}\label{conclu6}
In this work, an efficient uncertainty propagation method for a high-dimensional stochastic extension of Richards' equation was proposed. All the soil parameters were treated as unknown and modeled as random fields with appropriate marginal distributions. We also studied a modified Picard iteration and cell-centered multigrid method for solving the nonlinear systems with heterogeneous coefficients. We found that the combined solver is robust for a wide parameter range and the performance further improves with spatio-temporal refinements. This combination of solvers is general, therefore, its robustness can be further improved by incorporating adaptive time stepping or by combining with other advanced techniques, for instance, using the Anderson acceleration proposed in \cite{LOTT2012}.

For computing the statistical moments of the solution of Richards' equation, a parametric continuation technique based multilevel Monte Carlo estimator was proposed. This estimator is very practical for this problem, as it requires solving the strongly nonlinear problem only on the finest level, where the solver is robust, and uses simpler nonlinear problems on the coarse grid levels for a variance reduction. For the stochastic Richards' equation, the proposed estimator  is more prominent regarding the computational gains compared to the standard MLMC method if  the problem is strongly nonlinear.  In general, this estimator is also applicable to other parameter dependent nonlinear PDEs. One of the research problems that needs to be addressed is finding a computationally viable way of obtaining optimal step sizes for the nonlinear parameters used in continuation. This problem will be actively investigated in future work.
\appendix
%\subsection*{Appendix A1}\label{appendix_A1}
%\begin{eqnarray}\label{MLMC_sampling}
%\lnorm\mathbb{E}[p_{h_L,\boldsymbol{\Theta}_L}] -  \mathpzc{E}^{ML}_L[p_{h_L,\boldsymbol{\Theta}_L}] \rnorm_{\Lomd}^2&=& \mathbb{E}\left[\lnorm\sum^L_{\ell=0}\mathbb{E}[p_{h_\ell,\boldsymbol{\Theta}_\ell}-p_{h_{\ell-1},\boldsymbol{\Theta}_{\ell-1}}]-\mathpzc{E}^{MC}_{N_\ell}[p_{h_\ell,\boldsymbol{\Theta}_\ell}-p_{h_{\ell-1},\boldsymbol{\Theta}_{\ell-1}}]\rnorm^2_{L^2(\mathcal{D})} \right]\nonumber\\
%& =&  \mathbb{E}\sum^L_{\ell=0}\left[\lnorm\mathbb{E}[p_{h_\ell,\boldsymbol{\Theta}_\ell}-p_{h_{\ell-1},\boldsymbol{\Theta}_{\ell-1}}]-\mathpzc{E}^{MC}_{N_\ell}[p_{h_\ell,\boldsymbol{\Theta}_\ell}-p_{h_{\ell-1},\boldsymbol{\Theta}_{\ell-1}}]\rnorm^2_{L^2(\mathcal{D})} \right]\nonumber\\
%&=&\sum^L_{\ell=0}\mathbb{E}\left[\lnorm\mathbb{E}[u_{h_\ell}-u_{h_{\ell-1}}]-\mathpzc{E}^{MC}_{N_\ell}[p_{h_\ell,\boldsymbol{\Theta}_\ell}-p_{h_{\ell-1},\boldsymbol{\Theta}_{\ell-1}}]\rnorm^2_{L^2(\mathcal{D})} \right]\nonumber\\
%&=&\sum^L_{\ell=0}\left[\frac{1}{N_\ell^2}\sum_{i=1}^{N_\ell}\mathbb{E}\left[\lnorm\mathbb{E}[p_{h_\ell,\boldsymbol{\Theta}_\ell}-p_{h_{\ell-1},\boldsymbol{\Theta}_{\ell-1}}]-[p_{h_\ell,\boldsymbol{\Theta}_\ell}(\omega_i) - p_{h_{\ell-1},\boldsymbol{\Theta}_{\ell-1}}(\omega_i)]\rnorm^2_{L^2(\mathcal{D})}\right] \right]\nonumber\\
%&= &\sum^L_{\ell=0}\frac{\lnorm\mathcal{V}_\ell\rnorm_{\Lomd}}{N_\ell},
%\end{eqnarray}
\section{Sampling and upscaling of Gaussian random fields}\label{appendix_A1}
In this section, we outline the procedure for sampling Gaussian random fields that is used to sample the hydraulic conductivity fields \eqref{RandK} and also the non-Gaussian soil parameters \eqref{gPC}. For this, we use the Fast Fourier Transform Moving Average (FFT-MA) algorithm from \cite{Ravalec2000}. Although this sampling method is similar to the Cholesky decomposition technique, the FFT-MA method achieves a faster factorization of a covariance matrix by making the computational domain periodic. The resulting covariance operator is also periodic and can be decomposed as a convolutional product. This allows us to compute the samples of the random fields using cheaper vector-vector products compared to the expensive matrix-vector operation when using Cholesky factorization. Next, we provide a brief description of FFT-MA method from \cite{Ravalec2000}.

When using the Cholesky factorization, the samples of correlated Gaussian random vectors $\mathbf{z}(\omega)$ can be obtained as:
\begin{equation}
\mathbf{C}_{\Phi}  = \mathbf{L}\mathbf{L}^T\quad \text{ and use }\quad\mathbf{z} = \mathbf{L}\mathbf{y},
\end{equation}
where $\mathbf{C}_{\Phi}$ is the covariance matrix constructed on some grid and $\mathbf{y}$ is a vector of i.i.d. samples from the standard normal distribution. The FFT-MA relies on the decomposition of the covariance function $C_\Phi(r)$ as a convolutional product of some function $S_\Phi(r)$ and its transpose $S'_\Phi(r)= S_\Phi(-r)$. We can express this decomposition as
\begin{equation}\label{convol1}
\mathbf{c} = \mathbf{s}* \mathbf{s}',
\end{equation}
where $\mathbf{c},\mathbf{s}$ are vectors obtained by evaluating $C_\Phi(r)$ and $S_\Phi(r)$, respectively, at grid points of the considered mesh. Moreover, the resulting vector $\mathbf{s}$ is also real, positive and symmetric and $\mathbf{s} = \mathbf{s}'$. Now, a correlated random vector $\mathbf{z}$ can be computed by using the convolution product
\begin{equation}\label{convol2}
\mathbf{z} = \mathbf{s}*\mathbf{y}.
\end{equation}
The FFT-MA method performs the above computations in the frequency domain. As the convolution product in spatial domain is equivalent to the component-wise product in the frequency domain, we can take a Fourier transform of \eqref{convol1} as
\begin{equation}\label{product1}
\mathcal{F}(\mathbf{c}) = \mathcal{F}(\mathbf{s})\cdot \mathcal{F}(\mathbf{s}) \implies \mathcal{F}(\mathbf{s}) = \sqrt{\mathcal{F}(\mathbf{c})},
\end{equation}
where $\mathcal{F}$ denotes the discrete FFT and $\cdot$ denotes component-wise multiplication. As the FFT operation requires a periodic signal, first we transform the vector $\mathbf{c}$ into a periodic signal, which is also real, positive and symmetric. For more details on the practical aspects of this transformation see \cite{mishra2016multi,RF2,RF4,QMC}. Here the component-wise square-root operation does not pose any problems as the power spectrum $\mathcal{F}(\mathbf{c})$ is real and positive. Further, the convolution product in \eqref{convol2} can be expressed as a vector-vector product in the frequency domain as
\begin{equation}\label{product2}
\mathcal{F}(\mathbf{z}) = \mathcal{F}(\mathbf{s}*\mathbf{y})=\mathcal{F}(\mathbf{s})\cdot\mathcal{F}(\mathbf{y}).
\end{equation}
An inverse fast Fourier transform is finally applied to synthesize the samples for Gaussian random fields
\begin{equation}
\mathbf{z} = \mathcal{F}^{-1}(\mathcal{F}(\mathbf{s})\cdot\mathcal{F}(\mathbf{y})).
\end{equation} 
It is pointed out that due to the periodicity in the covariance vector $\mathbf{c}$, the resulting random field $\mathbf{z}$ is also periodic. Therefore, we only retain the part of the vector that corresponds to the physical domain and the remaining part is discarded. Also note that it takes two FFT evaluations to obtain one sample of $\mathbf{z}$ (ignoring the FFT operation in \eqref{product1} that is performed just once). For a given mesh, the cost of sampling random fields is negligible compared to the cost of solving the nonlinear PDE using the modified Picard-CCMG solver. 

Next, we describe the upscaling procedure for the random fields from grid $\mathcal{D}_\ell$ to $\mathcal{D}_{\ell-1}$. For clarity, we denote the above vectors computed on mesh $\mathcal{D}_\ell$ with subscript $\ell$, for example,  $\mathbf{z}_\ell, \mathbf{s}_\ell, \mathbf{y}_\ell$. While estimating the correction term $\mathpzc{E}^{MC}_{N_\ell}[p_{h_\ell,\boldsymbol{\Theta}_\ell} - p_{h_{\ell-1},\boldsymbol{\Theta}_{\ell-1}}]$ in the telescopic sum \eqref{MLMCestimator_ch6}, the approximations $p_{h_\ell,\boldsymbol{\Theta}_\ell}(\omega_i)$ and $p_{h_{\ell-1},\boldsymbol{\Theta}_{\ell-1}}(\omega_i)$ need to be positively correlated such that the level  dependent variance $\lnorm\mathcal{V}_\ell\rnorm_{\Ld}$ is small (see Eq. \eqref{level_var_ch6}). This is typically achieved by first sampling the fine grid Gaussian random field $\mathbf{z}_\ell$ to compute $p_{h_\ell,\boldsymbol{\Theta}_\ell}(\omega_i)$ and using an upscaled version $\mathbf{z}_{\ell-1}$ for $p_{h_{\ell-1},\boldsymbol{\Theta}_{\ell-1}}(\omega_i)$. While performing such upscaling of random fields, it is important to ensure that the telescopic sum \eqref{MLMCestimator_ch6} is not violated. In other words, the expectation of the random variable $p_{h_\ell,\boldsymbol{\Theta}_\ell}$ when estimating $\mathbb{E}[p_{h_\ell,\boldsymbol{\Theta}_\ell} - p_{h_{\ell-1},\boldsymbol{\Theta}_{\ell-1}}]$ and $\mathbb{E}[p_{h_{\ell+1},\boldsymbol{\Theta}_{\ell+1}} - p_{h_{\ell},\boldsymbol{\Theta}_{\ell}}]$ should be the same, i.e.
\begin{equation}\label{fine_coarse}
\mathbb{E}[p_{h_\ell,\boldsymbol{\Theta}_\ell}]^{(coarse)} = \mathbb{E}[p_{h_\ell,\boldsymbol{\Theta}_\ell}]^{(fine)}, \quad\text{for} \quad \ell = \{0,1,...,L-1\}.
\end{equation}
Using spatial averaging for obtaining an upscaled version may result in a modified covariance structure on the coarser levels, violating \eqref{fine_coarse}. This issue can be avoided by using the \emph{covariance upscaling} \cite{mishra2016multi} that employs the spectral generator on two consecutive grids using the same normally distributed vector $\mathbf{y}_\ell$. When using the FFT-MA algorithm, the vector $\mathbf{y}_\ell$ is associated with respective grid points, coarser realizations of the fine grid Gaussian random field $\mathbf{z}_\ell$ can be obtained by using multi-dimensional averaging of the vector $\mathbf{y}_\ell$. For instance, in two dimensions for the cell-centred mesh,
\begin{equation}\label{averaging}
\mathbf{y}^{i,j}_{\ell-1} = \frac{1}{2}(\mathbf{y}^{2i-1,2j-1}_{\ell}+\mathbf{y}^{2i-1,2j}_{\ell}+\mathbf{y}^{2i,2j-1}_{\ell}+\mathbf{y}^{2i,2j}_{\ell}),
\end{equation}
where $(i,j)$ is the cell index for the mesh $\mathcal{D}_{\ell-1}$. The scaling by a factor 2 is needed to obtain a standard normal distribution for the averaged quantity $\mathbf{y}^{i,j}_{\ell-1}$. The coarse random field can now be simply assembled as
\begin{equation}\label{cov_upscale}
\mathbf{z}_{\ell-1} = \mathcal{F}^{-1}(\mathcal{F}(\mathbf{s}_{\ell-1})\cdot\mathcal{F}(\mathbf{y}_{\ell-1})).
\end{equation}
This process can be recursively applied to generate upscaled random fields on next coarser scales. As the averaging in \eqref{averaging} smooths out high frequencies, the upscaled version $\mathbf{z}_{\ell-1}$ will also be  smoother compared to $\mathbf{z}_{\ell}$. These upscaled Gaussian random fields can be utilized to generate upscaled non-Gaussian fields using $\eqref{gPC}$.
\bibliographystyle{elsarticle-num}

\end{document}